\documentclass[12pt,a4paper]{report}
\usepackage[english]{babel}
\usepackage{amsmath,amssymb,amsfonts,amsthm}
\usepackage{graphicx}
\usepackage{geometry}
\usepackage[utf8]{inputenc}
\usepackage[T2A]{fontenc}
\usepackage{hyperref}
\usepackage[nottoc]{tocbibind}
\renewcommand{\phi}{\varphi}
\renewcommand{\epsilon}{\varepsilon}

\begin{document}
\begin{center}
National Academy of Sciences of Ukraine\\
V.M. Glushkov Institute of Cybernetics
\end{center}

\bigskip
 
\hspace{8cm} Manuscript work, 

\hspace{8cm} State registration number: 0406U001827 
\begin{center}
\bigskip
{\bf Bogdan V. NORKIN}
\end{center} 

\bigskip
\hspace{10cm} UDK 519.812.3; 519.21

\begin{center}
\bigskip

\smallskip
{\bf METHOD OF SUCCESSIVE APPROXIMATIONS FOR SOLVING INTEGRAL EQUATIONS OF ACTUARIAL MATHEMATICS
\\(English translation, originally in Russian)}

\bigskip
01.05.01 – Theoretical Foundations of Informatics and Cybernetics

\end{center}

\bigskip

\begin{center}
Dissertation for the degree of Candidate of Physical and Mathematical Sciences
\end{center}

\vspace{6cm}
\hspace{6cm} Scientific Supervisor:

\hspace{6cm} Doctor of Physical and Mathematical Sciences,  

\hspace{6cm} Professor, Academician

\hspace{6cm} of the National Academy of Sciences of Ukraine

\hspace{6cm} Yuri M. ERMOLIEV

\vspace{5cm}
\begin{center} KYIV - 2006 \end{center}

\setcounter{page}{1} 
\setcounter{chapter}{0}
\setcounter{section}{0}

\tableofcontents

\newpage
\chapter*{Abstract}
\addcontentsline{toc}{chapter}{Abstract}
Norkin B.V. Method of successive approximations for solving integral equations of actuarial mathematics. – Manuscript.

Dissertation for a scientific candidate degree of physical and mathematical sciences by specialty 01.05.01 – theoretical bases of computer science and cybernetics. – V.M. Glushkov Institute of Cybernetics of the National Academy of Sciences of Ukraine, Kiev, 2006.

In the dissertation various generalizations of the classical risk process  describing a stochastic evolution of the capital of an insurance company (Kramer–Lundberg model) are considered. In particular risk processes with variable deterministic premiums, with random premiums, with non Poisson flows of premiums and claims, risk processes in a stochastic markovian  environment are studied. Integral equations for the probability of nonruin as a function of the initial capital of the company for various generalizations of the classical risk process are deduced. For a risk process in a stochastic markovian environment system of integral equations for a set of nonruin probabilities from various initial states of the process are obtained. A general necessary and sufficient, and also concrete sufficient conditions of the existence and uniqueness of solutions of the considered integral equations and systems of equations are established. A successive approximation method for a numerical or analytical solution of the considered integral equations of insurance mathematics is theoretically and practically validated, in particular its uniform convergence and rate of such convergence are established. A technique of estimation of the accuracy of approximate solutions of integral equations of insurance mathematics is developed by construction of approximations to the exact solution from above and from below. The suggested method of consecutive approximations was tested on a number of numerical examples, its comparison with Monte Carlo method, with known approximations of solutions was carried out. 

The developed method of successive approximations for the solution of integral equations of insurance mathematics allows to increase accuracy of actuarial calculations, namely, to calculate within a chosen framework a probability of ruin of an insurance company with any prescribed accuracy, to check applicability and accuracy of various empirical approximate formulas for the probability of ruin, if necessary to improve accuracy of empirical approximations in an iterative way, to estimate accuracy and to correct parameters of Monte Carlo simulation method for calculation of the probability of ruin.

{\bf Keywords:} actuarial mathematics, insurance mathematics, risk process, ruin probability, probability of bankruptcy, integral equations, existence and uniqueness of solutions, successive approximation method.

\chapter*{INTRODUCTION}
\addcontentsline{toc}{chapter}{INTRODUCTION}

\textbf{Relevance of the topic.} Currently, the insurance market is actively developing in Ukraine [23]. Insurance of property, life, health, and civil liability of motorists has become firmly established in our lives. Medical and pension insurance should be mentioned separately – the development of these types of insurance constitutes an important part of the state social program. The business sector requires the development of insurance for investments, bank deposits, loans, cargo transportation, and other business areas. Therefore, improving actuarial calculation methods, insurance company management techniques, substantiation of insurance tariffs, and assessment of the risk of ruin of insurance companies will improve the quality and reduce the cost of insurance services for the population and business structures.

The main feature of insurance activity is that the business is based on randomness, so the corresponding models are stochastic. Insurance activity is aimed at protecting the client from the risk of large random losses through small but deterministic losses, which are called insurance premiums. The ability of an insurance company to conduct such a business is essentially based on the law of large numbers: by collecting a large number of independent insurance contracts, the company ensures an almost predictable dynamics of its capital.

Mathematical models of insurance activity were developed in the works of foreign scientists F. Lundberg, O. Lundberg, H. Cramér, K. Borch, S. Andersen, S. Asmussen, H.U. Gerber, Russian scientists A.N. Shiryaev, V.E. Bening, I.V. Evstigneev, V.V. Kalashnikov, A.V. Melnikov, G.I. Falin, and others, Ukrainian scientists Yu.M. Yermoliev, I.M. Kovalenko, V.S. Korolyuk, B.V. Bondarev, N.S. Bratijchuk, D.V. Gusak, Ya.I. Yeleyko, Yu.S. Mishura, A.N. Nakonechny, and others.

Modern insurance theory is based on the theory of stochastic processes. The oldest classical model of this kind is the F. Lundberg model (1903) (the classical risk process or compound Poisson risk process), which appeared around the same time as Bachelier's stochastic model of stock price changes (Brownian motion model). Unlike the Brownian motion model of stock price changes, the Lundberg model considers a jump stochastic process of capital change: capital monotonically increases between the moments of claim arrivals, and at the moment of an insurance claim, it instantly decreases by a random amount. The time intervals between claim arrivals are random and distributed according to the exponential law, and the claim sizes are independent and identically distributed. This model was then studied and developed in the works of H. Cramér, W. Feller, S. Andersen, K. Borch, H.U. Gerber, and others. Mathematical insurance theory includes such areas as modeling claim distributions, principles for calculating insurance premiums, reinsurance theory, ruin probability estimation, insurance portfolio management, and other sections. One of the important indicators of the stability of the insurance business is the probability of ruin of the company; therefore, the most important component of these studies was the study of the probability of ruin of an insurance company as a function of its initial capital and other parameters. This is a nontrivial problem because one needs to find the probability of the stochastic process trajectory entering the negative region over an infinite time interval. This probability can be calculated by the Monte Carlo method through direct simulation of the stochastic process. The difficulty is that in the most interesting and important case of small ruin probabilities, the Monte Carlo method becomes inefficient, i.e., it does not allow estimating a small probability with the necessary accuracy in a reasonable time. Therefore, various exact and approximate analytical methods, numerical methods, and practical approximate estimates of the ruin probability have been developed. They are based on the fact that for the classical risk process, the ruin probability satisfies a Volterra-type integral equation. It turns out that for more general models, the ruin probability as a function of initial capital also satisfies some more complex integral equations for which neither analytical nor numerical solution methods are known. This dissertation continues and develops this line of research: new more general integral equations and systems of equations for the probability of ruin of an insurance company operating under various conditions are derived, conditions for the existence and uniqueness of their solutions are established, and a general method of successive approximations for finding solutions is substantiated.

\textbf{Connection of the work with scientific programs, plans, and topics.} The work was carried out in accordance with the research plans of Department No. 130 of the V.M. Glushkov Institute of Cybernetics of the NAS of Ukraine within the following scientific topics:
\begin{itemize}
\item "Development of mathematical methods and computational algorithms for risk analysis, optimization, and management." Topic code in the institute plan - VFK 130.09, performed according to the resolution of the Presidium of the NAS of Ukraine, protocol of the Bureau of OI No. 8 of 12.01.2002.
\item "To develop risk assessment methods and their applications in economics, financial and insurance mathematics, and reliability theory." Topic code in the institute plan - VF 130.07, performed according to the resolution of the Presidium of the NAS of Ukraine, protocol of the Bureau of OI No. 6-A of 10.02.2000.
\end{itemize}

\textbf{Goal and objectives of the research.} The main goal of the work is to develop methods for assessing the risk of ruin of insurance companies.

This goal is achieved by constructing mathematical models of insurance activity in the form of stochastic risk processes, deriving integral equations (or systems of integral equations) for the probability of ruin of the risk process as a function of the initial conditions of the process, studying the properties of these equations, and developing methods for their solution.

\textbf{Scientific novelty of the obtained results.} The main results of the work are as follows:
\begin{itemize}
\item generalizations of the classical risk process, which describes the evolution of the capital of an insurance company, are considered, in particular, processes with variable deterministic premiums, with stochastic premiums, with non-Poisson flows of premiums and claims, and risk processes in a Markovian environment;
\item integral equations for the probability of (non-)ruin over an infinite time interval as a function of the initial capital of the company are generalized for the considered generalizations of the classical risk process;
\item systems of equations for the probabilities of non-ruin over an infinite time interval for a risk process in a Markovian environment are derived for the first time;
\item necessary and sufficient conditions for existence and general sufficient conditions for the existence and uniqueness of solutions of integral equations of insurance mathematics for the considered risk processes are established for the first time;
\item the method of successive approximations for the numerical solution of general integral equations of actuarial mathematics, i.e., the method for calculating the probability of ruin of insurance companies, is theoretically and practically substantiated for the first time;
\item conditions for uniform convergence and estimates of the rate of such convergence of the method of successive approximations in solving integral equations of insurance mathematics are obtained for the first time;
\item a methodology for estimating the accuracy of approximate solutions of the considered equations by constructing approximations to the solution from above and below is proposed for the first time.
\end{itemize}
Similar results are known in the literature only for the classical renewal equation; for more general equations and systems of integral equations of actuarial mathematics, these results are obtained for the first time.

\textbf{Practical significance of the obtained results.} The method of successive approximations for solving integral equations of insurance mathematics developed in the dissertation allows increasing the accuracy of actuarial calculations for complex models of insurance activity, in particular, for a chosen model of an insurance company, to estimate the probability of ruin of the company with any accuracy.

\textbf{Personal contribution of the applicant.} All publications on the topic of the dissertation were completed without co-authors.

\textbf{Approbation of the dissertation results.} The results of the research on this dissertation were reported at scientific seminars at the V.M. Glushkov Institute of Cybernetics of the NAS of Ukraine, at the Ivan Franko National University of Lviv, and at the following scientific conferences:
\begin{itemize}
\item International Conference "Problems of decision making under uncertainties (PDMU-2003)", Alushta, September 8-12, 2003. Poster presentation "On calculation of ruin probabilities".
\item International School-Seminar "Issues of optimization of computations", Katsiveli (Crimea), September 12-19, 2003. Report "Method for solving integral equations of insurance mathematics for finding the probability of ruin".
\item International Conference "Prediction and decision making under uncertainties (PDMU-2004)", Ternopil. Report "Calculation of ruin probabilities for a risk process in Markovian environment".
\item Tenth International Conference named after Academician G. Kravchuk, May 13-15, 2004, Kyiv. Report "Integral equations for dividends and the probability of ruin of an insurance company".
\item VI All-Ukrainian Conference of Young Scientists "Information Technologies in Education, Science, and Technology. ITONT-2004", April 28-30, Cherkasy. Report "Mathematical modeling of dividend payments and ruin risk of an insurance company".
\item International Conference "Functional methods in approximation theory, operator theory, stochastic analysis and statistics II" named after A.Ya. Dorogovtsev (1935-2004), October 1-5, 2004, Kyiv. Report "On solution of the basic actuarial integral equation".
\item International Conference "Modern problems and new trends in probability theory", June 19-26, 2005, Chernivtsi. Report "A successive approximation method for solution of actuarial integral equation".
\item International Conference "Problems of decision making under uncertainties (PDMU-2005)", September 12-17, 2005, Berdyansk. Report "On the theory of actuarial integral equations".
\end{itemize}

\textbf{Publications.} The results of the work were published in collections of the V.M. Glushkov Institute of Cybernetics of the NAS of Ukraine, in the scientific journals "Cybernetics and Systems Analysis" and "Problems of Control and Informatics" [40-45], and in conference abstracts [46, 47, 104-108].

\textbf{Main results.}

In the first section, a literature review on the dissertation topic is given.

The second section considers the classical model of insurance activity (F. Lundberg), which consists of modeling the dynamics of the capital of an insurance company as a compound Poisson risk process, and a number of its generalizations. In the classical model, insurance premiums are constant, so the capital of the insurance company linearly monotonically increases over time between claim arrivals, and at the moment of an insurance claim, it instantly decreases by a random amount. The time intervals between claim arrivals are randomly distributed according to the exponential law, and these intervals and claim sizes are independent and identically distributed. Models with a nonlinear premium income rate (depending on current capital), with a non-Poisson claim flow, models with stochastic premiums that arrive at random times, and models of the functioning of an insurance company in a random (Markovian) environment are also considered. For each type of risk process, the corresponding integral equations for the probability of non-ruin $\varphi(u)$ as a monotone function of the initial capital $u$ of the insurance company are given.

The integral equations of insurance mathematics for the probability of non-ruin $\varphi(u)$, $0\le \varphi(u)\le 1$, as a monotone function of the initial capital $u$ have the following general form:
\begin{equation}\label{eqn:0.1}
\varphi(u)=A\varphi(u),
\end{equation}
where $A$ is some linear integral operator, and the solution of the integral equation must satisfy the boundary condition
\begin{equation}\label{eqn:0.2}
\lim_{u\to\infty}\varphi(u)=1.
\end{equation}
This is a natural condition; it means that the company does not go bankrupt if it has an infinitely large initial capital. In the dissertation, it is proposed to solve all problems of actuarial mathematics of the form (\ref{eqn:0.1}), (\ref{eqn:0.2}) by a unified method, namely, the method of successive approximations:
\begin{equation}\label{eqn:0.3}
\varphi^{k+1}(u)=A\varphi^k(u),\qquad 0\le \varphi^0(u)\le 1,\qquad k=0,1,\dots
\end{equation}
It turns out that, generally speaking, the operator $A$ is not a contraction on its domain of definition, and the problem consists in substantiating the method of successive approximations.

In subsequent sections, the properties of the operator $A$ for specific types of risk processes, the conditions for the existence and uniqueness of the solution of problem (\ref{eqn:0.1}), (\ref{eqn:0.2}), the nature of successive approximations, and the conditions for convergence of $\varphi^k(u)$ to the solution of problem (\ref{eqn:0.1}), (\ref{eqn:0.2}) are investigated. The main element of the study of problem (\ref{eqn:0.1}), (\ref{eqn:0.2}) and method (\ref{eqn:0.3}) is the verification of the existence of monotone functions $\varphi_*(u)$, $\varphi^*(u)$ such that
\begin{align}
& 0\le \varphi_*(u)\le \varphi^*(u)\le 1,\nonumber\\
& \lim_{u\to\infty}\varphi_*(u)=1,\label{eqn:0.4}\\
& A\varphi_*(u)\ge \varphi_*(u),\qquad A\varphi^*(u)\le \varphi^*(u)\quad \forall u\ge 0.\nonumber
\end{align}
It turns out that one can always take $\varphi^*(u)\equiv 1$ as $\varphi^*(u)$, and the generalized Cramér-Lundberg bound depending on the form of the operator $A$ as $\varphi_*(u)$. In the dissertation, it is shown that for the considered operators, condition (\ref{eqn:0.4}) is necessary and sufficient for the existence of a solution of problem (\ref{eqn:0.1}), (\ref{eqn:0.2}), and the successive approximations starting from $\varphi^*(u)$ converge monotonically from above, while those starting from $\varphi_*(u)$ converge monotonically from below to the solution of problem (\ref{eqn:0.1}), (\ref{eqn:0.2}). The uniqueness of the solution of problem (\ref{eqn:0.1}), (\ref{eqn:0.2}) is proved for each equation separately using the properties of the corresponding operator $A$.

The third section investigates Picard's method of successive approximations for solving the renewal integral equation (Volterra type) satisfied by the probability of (non-)ruin of the classical risk process.

The classical risk process describing the evolution of the capital of an insurance company is given by the relation
\begin{equation}\label{eqn:0.5}
\xi_t=u+ct-S_t,
\end{equation}
where $t$ is time; $u$ is the initial capital of the insurance company; $c$ is the intensity of premium income; $S_t$ is the aggregate claim payments up to time $t$, $S_t=\sum_{k=1}^{N_t}Y_k$; $Y_k$ are independent identically distributed random variables (claims) with distribution function $F$ and mean $\mu$; $N_t$ is the number of payments up to time $t$ (a Poisson process with intensity $\alpha$). The theory of such processes has been studied in detail. As is known, the probability of non-ruin $\varphi(u)=\Pr\{\xi_t\ge 0\ \forall t\ge 0\}$ with initial capital $u$ satisfies the renewal integral equation (Volterra type):
\begin{equation}\label{eqn:0.6}
\varphi(u)=1-\frac{\alpha\mu}{c}+\frac{\alpha}{c}\int_0^u \varphi(u-z)(1-F(z))\,dz.
\end{equation}

As is known, the Volterra equation has a unique solution, which can be found by Picard's method of successive approximations:
\begin{equation}\label{eqn:0.7}
\varphi^{k+1}(u)=1-\frac{\alpha\mu}{c}+\frac{\alpha}{c}\int_0^u \varphi^k(u-z)[1-F(z)]\,dz,\qquad k=0,1,\dots,
\end{equation}
where $\varphi^0$ is some initial function, and $k$ is the iteration number. The method of successive approximations also allows finding analytical approximate solutions for the probability $\varphi$ as a function of $u$ and other parameters of the process (\ref{eqn:0.5}). As is known, Picard's method for Volterra integral equations converges on each finite interval of values of $u$ by virtue of the contraction mapping principle. However, in the general case, the method is not monotone, and contraction occurs only after a certain number of iterations; therefore, guaranteed convergence of the method for large $u$ may be slow. The method of successive approximations for solving equation (\ref{eqn:0.6}) was proposed in [13], where for the case of a continuous claim distribution function $F(x)$, its convergence and monotonicity were established. In this section, the properties of successive approximations for finding the probability of non-ruin $\varphi(u)$ of the classical risk process are investigated in detail in the general case. It is shown that for this Volterra equation, contraction occurs at each iteration with coefficient $p=\alpha\mu/c<1$ uniformly over all $u\in[0,+\infty)$, which ensures a high convergence rate for all $u$. It is shown that when starting from the initial functions $\varphi^0(u)\equiv 1$ and $\varphi^0(u)\equiv \varphi(0)=1-\alpha\mu/c$, the successive approximations converge monotonically and uniformly to the solution $\varphi(u)$ from above and below, respectively. It is shown that when starting with the initial function $\varphi^0(u)\equiv 1$, all successive approximations $\varphi^k(u)$ satisfy the condition $\varphi^k(u)\ge 1-e^{-Lu}$, where $L>0$ is the Cramér-Lundberg constant, which satisfies the equation
\begin{equation}
(\alpha/c)\int_0^{+\infty} e^{Rz}[1-F(z)]\,dz=1,\nonumber
\end{equation}
and thus $\lim_{k\to\infty}\varphi^k(u)=\varphi(u)\ge 1-e^{-Lu}$ and the condition $\varphi(+\infty)=1$ is satisfied.

In Sections 4 and 5, the results of Section 3 are generalized for the numerical solution of general (already non-Volterra) integral equations of the theory of risk processes. Consider a risk process of the form
\begin{equation}\label{eqn:0.8}
\xi_t=u+\int_0^t c(\xi_s)\,ds - S_t,\qquad t\ge 0,
\end{equation}
where the intensity of premium income $c(\cdot)\ge 0$ depends on the current capital of the company and is a piecewise continuous function; $S_t$ is the aggregate claim payments up to time $t$, $S_t=\sum_{k=1}^{N_t}Y_k$; $Y_k$ are independent identically distributed random variables (claims) with distribution function $F$ and mean $\mu$; $N_t$ is the number of payments up to time $t$. Let $K(t)$ be the distribution function of the time intervals between successive claim arrivals and the time of the first claim arrival. The case of a Poisson claim flow with $K(t)=1-e^{-\alpha t}$ is considered in Section 4, and the general case in Section 5. Define the function $U(u,t)$ as the solution of the Cauchy problem: $dU/dt=c(U)$, $U(0)=u$. For the probability of non-ruin $\varphi(u)$ of an insurance company with initial capital $u\ge 0$ over an infinite time interval, the following integral equation holds (for $U(u,t)=u+ct$ it is given in [35]):
\begin{equation}\label{eqn:0.9}
\varphi(u)=\int_0^{+\infty}\int_0^{U(u,t)}\varphi(U(u,t)-z)\,dF(z)\,dK(t),
\end{equation}
where the integrals are understood in the Lebesgue-Stieltjes sense. In the case of a linear function $U(u,t)$, taking into account the condition $\varphi(+\infty)=1$, this equation reduces to (\ref{eqn:0.6}); however, in the general nonlinear or non-Poisson case, it does not reduce to a Volterra equation.

Consider the following process of successive approximations:
\begin{equation}\label{eqn:0.10}
\varphi^{k+1}(u)=\int_0^{+\infty}\int_0^{U(u,t)}\varphi^k(U(u,t)-z)\,dF(z)\,dK(t),\qquad 0\le \varphi^0(u)\le 1,\qquad k=0,1,\dots
\end{equation}

{\bf Assumption 5.1.}
%\begin{assumption} 
%\label{ass:5.1}

a) The reserve growth function $U(u,t)$ is monotonically nondecreasing (increasing) in $u$ and $t$, $\lim_{t\to+\infty}U(u,t)=+\infty$;

b) there exist constants $u_*\ge 0$, $c_*>0$, and $L>0$ such that $U(u,t)\ge u+c_*t$ for all $u\ge u_*$, and
\begin{equation}
\int_0^{+\infty} e^{Lz}\,dF(z)\int_0^{+\infty} e^{-c_*Lt}\,dK(t)\le 1.\nonumber
\end{equation}

c) $K(t)\cdot F(z)<1\ \forall t,z\ge 0$.
%\end{assumption}

{\bf Theorem 5.1}
%\begin{theorem} \label{thm:5.1}
(On the existence and uniqueness of a solution and pointwise convergence of successive approximations). Suppose that Assumption 5.1 holds. Then
\begin{enumerate}
\item equation (\ref{eqn:0.9}) has a monotone solution $\varphi(u)$ such that
$1\ge \varphi(u)\ge \varphi_*(u)=\max\{0,1-e^{-L(u-u_*)}\}$;
\item $\varphi(u)$ is the unique solution of problem (\ref{eqn:0.9}), (\ref{eqn:0.2});
\item the sequence of approximations (\ref{eqn:0.10}) starting from $\varphi^*(u)\equiv 1$ converges monotonically pointwise from above, and the one starting from $\varphi_*(u)=\max\{0,1-e^{-L(u-u_*)}\}$ converges monotonically pointwise from below to the solution $\varphi(u)$;
\item for any initial approximation $\varphi^0(u)$ such that
$1\ge \varphi^0(u)\ge \varphi_*(u)=\max\{0,1-e^{-L(u-u_*)}\}$,
the corresponding sequence of approximations (\ref{eqn:0.10}) converges pointwise to the solution $\varphi(u)$.
\end{enumerate}
%\end{theorem}

{\bf Theorem 5.2}
%\begin{theorem} \label{thm:5.2}
(On the uniform convergence and the rate of convergence of the method of successive approximations). Under Assumption 5.1, the method of successive approximations (\ref{eqn:0.10}) starting with an initial approximation $\varphi^0(u)$ such that $1\ge \varphi^0(u)\ge \varphi_*(u)$ has the following properties:
\begin{itemize}
\item it converges uniformly to the solution $\varphi(u)$ of problem (\ref{eqn:0.9}), (\ref{eqn:0.2}), i.e., the sequence $\|\varphi^k-\varphi\|=\sup_{u\ge 0}|\varphi^k(u)-\varphi(u)|$ monotonically tends to zero;
\item moreover, it converges to any $\varepsilon$-neighborhood of the solution of problem (\ref{eqn:0.9}), (\ref{eqn:0.2}) at the rate of a geometric progression with a denominator depending on $\varepsilon$:
\begin{equation}
\|\varphi^k-\varphi\|\le (q(\varepsilon))^k\cdot \|\varphi^0-\varphi\|\nonumber
\end{equation}
for all $k$ such that $\|\varphi^k-\varphi\|\ge \varepsilon$.
\end{itemize}
%\end{theorem}

In Subsection 5.3, a numerical comparison of known empirical approximations (Beekman-Bowers, de Vylder, diffusion) for the ruin probability of a non-Poisson risk process (with deterministic time intervals between random claim arrivals) with the exact solution obtained by the method of successive approximations is carried out. An example shows that the approximations can give a significant error for the ruin probability.

In Section 6, a risk process in which both claims and premiums are stochastic is investigated. It describes the stochastic process of customer arrivals to an insurance company more adequately than the classical Cramér-Lundberg model. In [5], for a similar process with independent Poisson flows of stochastic premiums and claims, an integral equation for the probability of non-ruin was derived and some analytical solutions of this equation were found. In this section, more general integral equations for the probability of non-ruin are derived in the case of possibly non-Poisson flows of stochastic premiums and claims and in the presence of deterministic payments, and the method of successive approximations for solving these equations is substantiated. Namely, we establish necessary and sufficient conditions for existence and general sufficient conditions for the uniqueness of solutions of the equations, prove the convergence of the method of successive approximations, and estimate the rate of convergence.

Consider a stochastic process (risk) $\xi_t$ describing the evolution in time $t$ of the capital of an insurance company and satisfying the stochastic integral equation
\begin{equation}\label{eqn:0.11}
\xi_t=u+\int_0^t c(\xi_s)\,ds+P_t-S_t,\qquad t\ge 0,
\end{equation}
where $u$ is the initial capital of the company, $\xi_0=u$; $c(\xi)$ is the intensity of deterministic payments depending on the current capital (a piecewise continuous function); $S_t=\sum_{k=1}^{N_t^c} z_k$ is the aggregate insurance claims; $P_t=\sum_{k=1}^{N_t^p} p_k$ is the aggregate insurance premiums; $z_k$ are independent random variables (claims) with common distribution function $F(z)$; $p_k$ are independent random variables (premiums) with common distribution function $G(z)$; $N_t^c$ is the number of claims up to time $t$, in particular, this may be a Poisson process with intensity $\alpha_c$; similarly, $N_t^p$ is the number of premiums up to time $t$, in particular, this may be a Poisson process with intensity $\alpha_p$. In [5], for the case $c(\cdot)\equiv 0$ and Poisson flows of premiums and claims, the following integral equation for the probability of non-ruin $\varphi(u)$ of process (\ref{eqn:0.11}) was derived:
\begin{equation}\label{eqn:0.12}
\varphi(u)=A_1\varphi(u):=\frac{\alpha_c}{\alpha_c+\alpha_p}\int_0^u \varphi(u-z)\,dF(z)+\frac{\alpha_p}{\alpha_c+\alpha_p}\int_0^{+\infty}\varphi(u+z)\,dG(z).
\end{equation}

Consider a more general process constructed as follows. At the initial time $t=0$ and at each jump time $t'$ of the process (arrival of a stochastic premium or claim), two random variables $\tau_p$ and $\tau_c$ with distribution functions $K_p(t)$ and $K_c(t)$, respectively, are generated. The next jump occurs at time $t''=t'+\min\{\tau_p,\tau_c\}$. If $\tau_p<\tau_c$, the capital increases by a random premium with distribution $G(\cdot)$; otherwise, it decreases by a claim with distribution $F(\cdot)$. Between jumps, the capital changes deterministically according to the equation $d\xi_t/dt=c(\xi_t)$.

The described process includes the case of independent Poisson processes of premium and claim arrivals when $K_p(t)=1-e^{-\alpha_p t}$, $K_c(t)=1-e^{-\alpha_c t}$.

Define the function $U(u,t)$ as the solution of the Cauchy problem: $dU/dt=c(U)$, $U(0)=u$.

{\bf Assumption 6.1.}
%\begin{assumption} \label{ass:6.1}
The functions $K_p(t)$ and $K_c(t)$ have no common discontinuity points, and $K_p(t)<1$, $K_c(t)<1$ for all $t\ge 0$.
%\end{assumption}

{\bf Theorem 6.2.}
%\begin{theorem} \label{thm:6.2}
For the constructed process, the probability of non-ruin $\varphi(u)$ satisfies the integral equation
\begin{align}
\varphi(u)=A_2\varphi(u):=&\int_0^\infty\int_0^{U(u,t)}\varphi(U(u,t)-z)\,dF(z)(1-K_p(t))\,dK_c(t)\nonumber\\
&+\int_0^\infty\int_0^{+\infty}\varphi(U(u,t)+z)\,dG(z)(1-K_c(t))\,dK_p(t),\qquad \varphi(+\infty)=1.\label{eqn:0.13}
\end{align}
%\end{theorem}

{\bf Assumption 6.2.}
%\begin{assumption} \label{ass:6.2}
For the operator $A\in\{A_1,A_2\}$, there exist functions $\varphi_*(u),\varphi^*(u)\in H$ such that conditions (4) hold.
%\end{assumption}

{\bf Theorem 6.3}
%\begin{theorem} \label{thm:6.3}
(Necessary and sufficient conditions for the existence and uniqueness of the solution of problems (\ref{eqn:0.12}), (\ref{eqn:0.2}) and (\ref{eqn:0.13}), (\ref{eqn:0.2})). For the existence and uniqueness of a solution $\varphi\in H$ of problems (\ref{eqn:0.12}), 
(\ref{eqn:0.2}) and (\ref{eqn:0.13}), (\ref{eqn:0.2})), it is necessary and sufficient that Assumption 6.2 holds with the corresponding operator $A\in\{A_1,A_2\}$.
%\end{theorem}

In the section, it is shown that one can take $\varphi^*(u)\equiv 1$ as $\varphi^*(u)$. Another possible function $\varphi^*(u)$ has the form:
\begin{equation}
\varphi^*(u)=\frac{\alpha_p}{\alpha_p+\alpha_c(1-F(u))}.\nonumber
\end{equation}

The following two lemmas give sufficient conditions for the existence of the function $\varphi_*(u)$.

{\bf Lema 6.4.}
%\begin{lemma} \label{lem:6.4}
Let there exist a positive solution $L_1$ of the inequality
\begin{equation}
\frac{\alpha_c}{\alpha_c+\alpha_p}Ee^{L_1 z_i}+\frac{\alpha_p}{\alpha_c+\alpha_p}Ee^{-L_1 p_i}
=\frac{\alpha_c}{\alpha_c+\alpha_p}\int_0^{+\infty} e^{L_1 z}\,dF(z)+\frac{\alpha_p}{\alpha_c+\alpha_p}\int_0^{+\infty} e^{-L_1 z}\,dG(z)\le 1,\nonumber
\end{equation}
and $\lim_{z\to\infty} e^{L_1 z}(1-F(z))=0$. Set $\varphi_*(u)=1-e^{-L_1 u}$. Then $\varphi_*(u)$ is a nondecreasing function such that $0\le \varphi_*(u)\le 1$, $\lim_{u\to+\infty}\varphi_*(u)=1$, and for the operator $A_1$, the relation $A_1\varphi_*(u)\ge \varphi_*(u)$ holds.
%\end{lemma}

{\bf Lema 6.5.}
%\begin{lemma} \label{lem:6.5}
Let there exist constants $u_*\ge 0$, $c_*>0$, $L_2>0$ such that
a) $U(u,t)\ge u+c_*t$ for all $u\ge u_*$, and
b) the inequality
\begin{align}
&\int_0^{+\infty} e^{L_2 z}\,dF(z)\cdot \int_0^\infty e^{-c_*L_2 t}(1-K_p(t))\,dK_c(t)\nonumber\\
&+\int_0^{+\infty} e^{-L_2 z}\,dG(z)\cdot \int_0^\infty e^{-c_*L_2 t}(1-K_c(t))\,dK_p(t)\le 1
\nonumber
\end{align}
holds, and $\lim_{z\to\infty} e^{L_2 z}(1-F(z))=0$. Set $\varphi_*(u)=\max\{0,1-e^{-L_2(u-u_*)}\}$. Then $\varphi_*(u)$ is a nondecreasing function such that $0\le \varphi_*(u)\le 1$, $\lim_{u\to+\infty}\varphi_*(u)=1$, and for the operator $A_2$, the relation $A_2\varphi_*(u)\ge \varphi_*(u)$ holds.
%\end{lemma}

{\bf Corollary 6.1.}
%\begin{corollary} \label{cor:6.1}
When $K_p(t)=1-e^{-\alpha_p t}$, the inequality for the constant $L_2$ takes the form:
\begin{equation}
\frac{\alpha_c}{\alpha_c+\alpha_p+c_*L_2}\int_0^{+\infty} e^{L_2 z}\,dF(z)+\frac{\alpha_p}{\alpha_c+\alpha_p+c_*L_2}\int_0^{+\infty} e^{-L_2 z}\,dG(z)\le 1.\nonumber
\end{equation}
If $F(\cdot)$ and $G(\cdot)$ are such that
$\int_0^{+\infty} e^{L_2 z}\,dF(z)=1+\mu_c L_2+o(L_2)$ and
$\int_0^{+\infty} e^{-L_2 z}\,dG(z)=1-\mu_p L_2+o(L_2)$
(for example, claims and premiums are bounded), where $\mu_c$ and $\mu_p$ are the mean claim and premium sizes, then a positive solution $L_2$ of the inequality exists under the condition $\alpha_c\mu_c<\alpha_p\mu_p+c_*$.
%\end{corollary}

To solve problems (\ref{eqn:0.12}), (\ref{eqn:0.2}) and (\ref{eqn:0.13}), (\ref{eqn:0.2}), we apply the method of successive approximations (\ref{eqn:0.3}) with the corresponding operator $A\in\{A_1,A_2\}$.

{\bf Theorem 6.4}
%\begin{theorem} \label{thm:6.4}
(On the pointwise convergence of approximations from above and below). Under Assumptions 6.1 and 6.2, the sequence of approximations $\{\varphi^k(u)\}$ starting from $\varphi^0(u)\equiv \varphi^*(u)$ decreases monotonically and converges pointwise from above to the solution of problems (\ref{eqn:0.12}), (\ref{eqn:0.2}) and (\ref{eqn:0.13}), (\ref{eqn:0.2}), while the sequence of approximations $\{\varphi^k(u)\}$ starting from $\varphi^0(u)\equiv \varphi_*(u)$ increases monotonically and converges pointwise from below to the solution.
%\end{theorem}

{\bf Corollary 6.2}
%\begin{corollary} \label{cor:6.2}
(On pointwise convergence of approximations). Under Assumptions 6.1 and 6.2, for any initial approximation $\varphi^0\in H^*$, the sequence $\{\varphi^k(u),\ k=0,1,\dots\}$ converges pointwise to the solution of problems (\ref{eqn:0.12}), (\ref{eqn:0.2}) and (\ref{eqn:0.13}), (\ref{eqn:0.2}).
%\end{corollary}

{\bf Theorem 6.5}
%\begin{theorem} \label{thm:6.5}
(On uniform convergence and the rate of convergence of the method of successive approximations). Under Assumptions 6.1 and 6.2 and with the condition $F(\cdot)G(\cdot)<1$, the method of successive approximations (\ref{eqn:0.3}) starting with $\varphi^0(u)$ such that $\varphi_*(u)\le \varphi^0(u)\le 1$ converges uniformly monotonically to the solution of problems (\ref{eqn:0.12}), (\ref{eqn:0.2}) and (\ref{eqn:0.13}), (\ref{eqn:0.2}); moreover, it converges to any $\varepsilon$-neighborhood of the solution at the rate of a geometric progression with a denominator depending on $\varepsilon$.
%\end{theorem}

In the seventh section, a general risk process describing the stochastic evolution of the capital of an insurance company in a Markovian random environment with nonlinear premium income is considered. It is shown that the probabilities of non-ruin of the process (company) as functions of the initial state of the process generally satisfy a system of integral equations with boundary conditions at infinity. Sufficient conditions for the existence of a solution of this system are established, and the method of successive approximations for calculating the ruin probability is substantiated. When starting with the unit initial approximation, the iterations converge to the solution from above, and when starting with the generalized Cramér-Lundberg bound, they converge from below; thus, there is always a possibility to estimate the accuracy of the obtained approximate solution. The workability of the method is verified on a numerical example of a risk process on a Markov chain with two states.

Consider a stochastic risk process $\xi_t$ in a random (Markovian) environment that can be in one of $i=1,\dots,n$ states. Let $\alpha_i$ be the intensity of the Poisson flow of insurance claims; $c_i(\xi_t)$ be the intensity of premium income; $F_i$ be the distribution function of claims in state $i$ of the environment. Let $U_i(u,t)$ be the solution of the Cauchy problem: $dU_i/dt=c_i(U_i)$, $U_i(0)=u$. Let the probability that the chain remains in state $i$ for time $\Delta$ be $e^{-\lambda_i\Delta}$, and the probability of transition of the chain from state $i$ to state $j\ne i$ during time $\Delta$ be $\lambda_i p_{ij}\Delta+o(\Delta)$, where $\sum_{j=1}^n p_{ij}=1$, $p_{ii}=0$, $\lim_{\Delta\to 0} o(\Delta)/\Delta=0$. Then the probability of non-ruin of the process depends on the initial state of the process $(u,i)$, and thus we have a set of probabilities
$\{\varphi_i(u)=P[\inf_{0\le t<\infty}\xi_t>0|\xi_0=u,\text{ init. state of env.}=i]\}$.

{\bf Theorem 7.4.}
%\begin{theorem} \label{thm:7.4}
The vector function $\varphi(u)=\{\varphi_1(u),\varphi_2(u),\dots,\varphi_m(u)\}$ of probabilities of non-ruin $\{\varphi_i(u)\}$ of the process from the initial position $(u,i)$ satisfies the system of integral equations:
\begin{align}
\varphi_i(u)=&\int_0^\infty e^{-t(\alpha_i+\lambda_i)}
\left\{\lambda_i\sum_{j=1}^m p_{ij}\varphi_j(U_i(u,t))\right.\nonumber\\
&\left.+\alpha_i\int_0^{U_i(u,t)}\varphi_i(U_i(u,t)-z)\,dF_i(z)\right\}dt,\qquad i=1,\dots,m.
\label{eqn:0.14}
\end{align}
%\end{theorem}

The system of equations (\ref{eqn:0.14}) is linear. We are interested in a solution satisfying the natural boundary conditions
\begin{equation}
\lim_{u\to+\infty}\varphi_i(u)=1,\qquad i=1,\dots,m,\label{eqn:0.15}
\end{equation}
meaning that with unboundedly large initial capital, the process does not go bankrupt. Such a solution may not exist, for example, if the process goes bankrupt from any initial state, as is the case when the average claims per unit time are greater than the corresponding premiums. Then the solution is trivial, $\varphi_i(u)\equiv 0,\ i=1,\dots,m$.

Consider the following method of successive approximations for finding the functions $\{\varphi_1(u),\dots,\varphi_m(u)\}$:
\begin{align}
\varphi_i^{k+1}(u)=&\int_0^\infty e^{-t(\alpha_i+\lambda_i)}
\left\{\lambda_i\sum_{j=1}^m p_{ij}\varphi_j^k(U_i(u,t))\right.\nonumber\\
&\left.+\alpha_i\int_0^{U_i(u,t)}\varphi_i^k(U_i(u,t)-z)\,dF_i(z)\right\}dt,\qquad i=1,\dots,m;\quad k=0,1,\dots,\label{eqn:0.16}
\end{align}
where $\{\varphi_i^0(u)\}$ are some initial functions such that $0\le \varphi_i^0(u)\le 1$, $\varphi_i^0(+\infty)=1$.

{\bf Assumption 7.1.}
%\begin{assumption} \label{ass:7.1}
Let in each state $i$ of the environment:
\begin{itemize}
\item[a)] the reserve growth function $U_i(u,t)$ is monotonically increasing in $u$ and $t$, $\lim_{t\to+\infty}U_i(u,t)=+\infty$, and $U_i(u,t)\ge u+c_i t$ for $u\ge u_*$ with some constant $c_i\ge 0$;
\item[b)] $\alpha_i\mu_i/c_i<1$, and thus there exists a common constant for all states (an analogue of the Cramér-Lundberg constant) $L>0$ such that $(\alpha_i/c_i)\int_0^{+\infty} e^{Lz}[1-F_i(z)]\,dz\le 1$.
\end{itemize}
%\end{assumption}

{\bf Theorem 7.5}
%\begin{theorem} \label{thm:7.5}
(On the convergence of the method of successive approximations). Suppose that Assumption 7.1 (a), (b) holds. Then
\begin{itemize}
\item system (\ref{eqn:0.14}), (\ref{eqn:0.15}) has a monotone solution with components $\varphi_i(u)$ such that $1\ge \varphi_i(u)\ge \varphi_*(u)=\max\{0,1-e^{-L(u-u_*)}\}$;
\item such a solution is unique;
\item the sequence of approximations (\ref{eqn:0.16}) starting from $\{\varphi_i^0(u)\equiv 1\}$ converges monotonically pointwise from above, and the one starting from $\{\varphi_i^0(u)=\varphi_*(u)=\max\{0,1-e^{-L(u-u_*)}\}\}$ converges monotonically pointwise from below to the solution $\{\varphi_i(u)\}$ for each $i$;
\item for any initial approximation $\{\varphi_i^0(u)\}$ such that $1\ge \varphi_i^0(u)\ge \varphi_*(u)=\max\{0,1-e^{-L(u-u_*)}\}$, the corresponding sequence of approximations (\ref{eqn:0.16}) converges pointwise to the solution $\{\varphi_i(u)\}$.
\end{itemize}
%\end{theorem}

\chapter{LITERATURE REVIEW}
%\addcontentsline{toc}{chapter}{LITERATURE REVIEW}

This section is devoted to a review of mathematical models of insurance activity, mathematical problems, and methods for solving the arising mathematical problems.

\section{Mathematical modeling of insurance activity}

There are several approaches to mathematical modeling of insurance company activities. The main feature of insurance activity is that the business is based on randomness, so the corresponding models are stochastic. Insurance activity is aimed at protecting the client from the risk of large random losses through small but deterministic losses, which are called insurance premiums. The ability of an insurance company to conduct such a business is essentially based on the law of large numbers: by collecting a large number of independent insurance contracts, the company ensures an almost predictable dynamics of its capital. Issues of mathematical modeling of insurance activity are discussed in monographs [6, 35, 72, 87, 92], in professional actuarial journals "Insurance: Mathematics and Economics", "Scandinavian Actuarial Journal", "AUSTIN Bulletin", issues of the journals "Theory of Probability and its Applications", "Review of Applied and Industrial Mathematics" (Russia), "Applied Statistics, Actuarial and Financial Mathematics" (Donetsk), etc.

Insurance is inextricably linked with the concept of risk.

The concept of "risk" came to us from the distant past; in its modern broad sense, it was formed in the Middle Ages [70, p. 147]. The concept of risk is widely used not only in insurance; it has entered everyday life. The word "risk" comes from the Italian word "rischiare", which in turn came from the Greek 'rhiza' (cliff), which since the time of the Greeks meant danger. In Italy, the word risk was used mainly when it came to the loss of ship cargoes coming from the East. The cause was either a storm and shipwreck, or pirates.

An insurance company can be considered as an economic agent maximizing its objective function (or functions) under constraints. The objective function usually used is the expected utility function or the expected profit function. Constraints include technical (regulatory), budget, reliability constraints that regulate the probability of ruin of the company, etc. Probabilistic constraints distinguish this economic agent from other agents (consumers, producers, etc.). Models of this kind belong to stochastic programming models with probabilistic constraints and are considered in [57, 72, 83, 86, 92].

In a broad sense, the model of an economic agent is a system of constraints limiting its freedom of choice. An insurance company can choose the volumes of a particular type of insurance, insurance tariffs, parameters of risk reinsurance, the volume of investment of its capital in other types of financial and economic activity. The parameters of its activity are limited by regulatory requirements for asset liquidity, reliability, etc. For example, since the probability of ruin depends on insurance premium rates (or the amount of insurance load), many works study the question of what level of insurance premiums or load ensures a given level of company reliability [38, 72]. This question is not simple because the probability of ruin is not known explicitly as a function of management parameters.

Any insurance company operates in some economic, in particular, market environment, i.e., among other companies and under conditions of emerging demand and supply of insurance services. On the one hand, to increase profit and improve business reliability, an insurance company wants to receive high insurance premiums, but on the other hand, it must monitor market conditions and reduce tariffs to attract clients. Therefore, there are different models and principles of pricing for insurance products [55, 56, 72, 87].

An insurance company can be considered as an investor directing money to insure various objects and various types of insurance contracts. The assets in this case are insurance contracts, the profitability of which includes both a deterministic component (insurance premium) and a random component (insurance claims). Such activity promises both a constant inflow of capital and is associated with a random flow of insurance claims and the risk of large losses or even ruin. Therefore, for modeling the activity of an insurance company, the methods of portfolio investment theory with necessary modifications are also applicable in principle. Optimization and management in this case is subject to the portfolio of contracts, i.e., the number of contracts of different types or for different insurance objects. Models of this kind are considered in [21, 58, 81, 82]. To build models of this kind, new methods of portfolio theory based on risk measures and the concept of stochastic dominance can be applied [64, 110, 113]. New problems arise in connection with the need to insure catastrophic risks when risks are dependent and the law of large numbers is not applicable [21, 58, 81, 133].

The evolution of the capital of an insurance company can be modeled by methods of the theory of stochastic processes. The oldest classical model of this kind is the F. Lundberg model (1903) [99] (classical risk process or compound Poisson risk process), which appeared around the same time as Bachelier's stochastic model of stock price changes (Brownian motion model). Unlike the Brownian motion model of stock price changes, the Lundberg model considers a jump stochastic process of capital change: capital monotonically increases between claim arrivals, and at the moment of an insurance claim, it instantly decreases by a random amount. The time intervals between claim arrivals are random and distributed according to the exponential law, and the claim sizes are independent and identically distributed. This model was then studied and developed in the works of H. Cramér, W. Feller, S. Andersen, K. Borch, H.U. Gerber, and others [55]. An important component of these studies was the study of the probability of ruin of an insurance company as a function of its initial capital and other parameters. This is a nontrivial problem because one needs to find the probability of the stochastic process trajectory entering the negative region over an infinite time interval. The main method of studying the ruin probability in these works was the method of integral equations. In particular, integral equations for this function were established and studied, particular solutions, asymptotics and estimates of solutions were obtained, and some analytical and numerical methods for their solution were developed. This dissertation continues and develops this line of research: new more general integral equations for the probability of (non-)ruin of an insurance company operating under various conditions are derived, necessary and sufficient conditions for the existence and uniqueness of their solutions are established, and a general method of successive approximations for finding solutions is substantiated.

\section{Classical model of collective risk and methods for finding the probability of ruin}

As early as 1903, F. Lundberg proposed a model of an insurance company based on a stochastic process with independent increments [99], which describes the stochastic evolution of the capital of an insurance company. In this model, on the one hand, capital monotonically and linearly increases over time due to continuously incoming premiums, and on the other hand, at random times (of insurance claim arrivals), it decreases by a random amount (claims) [55]. The company goes bankrupt if its capital becomes less than zero.

The most important characteristic of the company's stability is its probability of non-ruin.

Formally, the classical risk process describing the evolution over time of the capital $\xi_t$ of an insurance company is given by the relation (see [35, 51, 65, 72]):
\begin{equation}
\xi_t=u+ct-S_t,\label{eqn:1.1}
\end{equation}
where $t$ is time; $u$ is the initial capital of the insurance company; $c$ is the intensity of premium income; $S_t$ is the aggregate payments up to time $t$, $S_t=\sum_{k=1}^{N_t} Y_k$; $Y_k$ are independent identically distributed random variables (claims) with distribution function $F$ and mean $\mu$; $N_t$ is the number of payments up to time $t$ (a Poisson process with intensity $\alpha$). The theory of such processes has been studied in detail. As is known, the probability of non-ruin $\varphi(u)=\Pr\{\xi_t\ge 0\ \forall t\ge 0\}$ with initial capital $u$ satisfies the renewal integral equation [35, p. 227]:
\begin{equation}
\varphi(u)=1-\frac{\alpha\mu}{c}+\frac{\alpha}{c}\int_0^u \varphi(u-z)(1-F(z))\,dz.
\label{eqn:1.2}
\end{equation}

This equation can be solved analytically in some cases, and in the general case by numerical methods. Among the analytical methods, the Laplace transform method [35] and the factorization method [11] should be noted; among numerical methods, Picard's method of successive approximations [13, 41-43], the Monte Carlo method [72], and finite difference methods [55, 68, 114]. In practice, a number of approximations are also often used, such as the Cramér-Lundberg, De Vylder, and Beekman-Bowers approximations [35]. The construction of approximations is based on replacing the original process with another risk process for which the exact solution for the ruin probability is known. It should be noted that in special cases, other methods can be used to obtain analytical solutions; for example, [63] found an analytical solution for the ruin probability over a finite time interval in the case of an exponentially distributed claim amount.

The Cramér-Lundberg approximation (from below) for the solution of equation (\ref{eqn:1.2}) has the form $\varphi(u)\ge 1-e^{-Ru},\ u\ge 0$, where the Cramér-Lundberg coefficient $R$ satisfies the equation:
\begin{equation}
\frac{\alpha}{c}\int_0^{+\infty} e^{Rz}[1-F(z)]\,dz=1.\label{eqn:1.3}
\end{equation}
An approximation from above is also known for processes with bounded claims (if the claim size is bounded by some constant $m$, then the ruin probability $\varphi(u)\le 1-e^{-(u+m)}$, see [76]), approximations by V.V. Kalashnikov [24] from above and below for large claims, and a number of others [53]. The main disadvantage of using approximations is their possible large relative error [44] and, moreover, the impossibility of estimating the accuracy in the case of one-sided approximations (in particular, the Cramér-Lundberg approximation).

\section{Generalizations of the classical model}

There are a number of generalizations of the classical risk process [55], such as the risk process with a non-Poisson claim flow, processes with nonlinear and stochastic premiums (in [4], a process with a Poisson flow of premiums of the same size and another independent Poisson flow of independently distributed claims is considered, for which the distribution of the absolute minimum over a finite time interval is obtained), mixed Poisson processes with random claim arrival intensity [88], risk processes on a Markov chain, etc., which more adequately model the operation of an insurance company. For the ruin probability in the case of such more general processes, integro-differential equations can also be derived (and in the case of Markov-modulated processes, systems of equations), described in detail in Chapter 2 of this dissertation. Let us first consider analytical solution methods. The risk process with a dividend barrier (the company's capital is bounded above by some time-dependent function, the so-called dividend barrier, upon reaching which all money above it is paid to shareholders as dividends) was considered in [60, 61, 77, 87]. General risk processes with dependent claims were considered in [21]. In [58], the dependence of the ruin probability on the degree of claim dependence is investigated. In particular, an example is considered in which claims, as well as their arrival time, are dependent within autoregressive models. Calculations were performed by the Monte Carlo method.

\section{Factorization method and other analytical methods}

The idea of using the factorization method developed by M.G. Krein in [34] for Volterra-type integral equations for integro-difference equations belongs to V.S. Korolyuk. He applied this method in [29] and other works for the analysis of the asymptotics of limit problems for random walks. A.A. Borovkov and B.A. Rogozin [8] and other authors developed the factorization approach to limit problems for random walks and processes with independent increments. A fairly wide class of functionals for processes with independent increments was investigated in the works of D.V. Gusak and N.S. Bratijchuk [11] and other students of V.S. Korolyuk using the factorization method.

Among other analytical methods, the potential method developed by V.S. Korolyuk in [30, 33] for random walks and compound Poisson processes, and the combinatorial method most fully presented in the monograph by L. Takács [50], should be noted. In [96, 97], J. Kemperman developed the Wiener-Hopf method and the projection method in the study of Markovian random walks and homogeneous Markov processes with constraints. The potential method for semi-continuous processes with independent increments in the case of unbounded variance found further development in the works of V.S. Korolyuk, V.M. Shurenkov, V.M. Suprun [31, 32], and N.S. Bratijchuk [9, 10]. The potential method allows refining the asymptotic behavior of the distribution of many limit functionals.

Compound Poisson processes with jumps of one sign have the property of semi-continuity. Based on this property, in [11, 18], the following important result is established: the factorization components in the main factorization identity are determined by the distribution of the positive or negative part of the values of the initial process. Therefore, the distributions of some limit functionals are expressed through the distributions of the positive or negative values of the process. In addition, in the case of semi-continuity, the property of exponentiality of the distribution of the corresponding functional, which determines one of the factorization components, holds.

Using the combinatorial method, L. Takács [50, 115, 116] obtained interesting results for the distribution of extrema of compound Poisson processes and recurrent processes, which generalize the class of Poisson processes (by removing the constraint on the exponentiality of the duration distribution between adjacent jumps). In addition to limit functionals for processes and random walks, functionals of additive type attract attention. The study of these functionals and other aspects of the study of stochastic processes are devoted to the monographs by I.I. Gikhman and A.V. Skorokhod [19, 49], J. Cohen and O. Boxma [24], J. Kemperman [96], and others.

\section{Risk processes in a Markovian environment}

The approach for obtaining the distribution of the absolute minimum of a process, as well as the corresponding probabilities of (non-)ruin, for general processes with independent increments on Markov chains is based on matrix integro-differential equations with partial derivatives for the matrix distribution function [18, p. 88]
$\{\varphi_{ij}(t,x)=P[\min_{0\le \tau\le t}\xi(\tau)>x | \xi(0)=0, \text{init. state}=i, \text{fin. state}=j]\}$. However, these matrix equations are very complex. They contain partial derivatives with respect to time $t$, in contrast to similar equations for the classical risk process. To solve these equations, analytical methods using matrix characteristic functions are developed in [18]. Therefore, it makes sense to derive similar simpler equations directly for the probabilities
$\{\varphi_i(u)=P[\min_{0\le \tau<\infty}\xi(\tau)>0|\xi(0)=u, \text{init. state of chain}=i]\}$
of non-ruin of the process over an infinite time interval with initial capital $u$. These equations no longer contain partial derivatives with respect to time $t$, unlike the equations given in [18, (1.8)], and direct methods, for example, the known Laplace transform, can be applied to solve them. The effectiveness of this approach is illustrated by an example in Section 7 for a chain with two states.

In the monograph by D.V. Gusak [18], a significant part of the results for the distribution of functionals from ordinary processes with independent increments and from processes on a Markov chain is obtained by using the factorization method and expressing the desired distributions in terms of characteristic functions of extrema and their complements, i.e., in terms of components of the main factorization identities. Identities of this type with the corresponding generalization are established for processes on a Markov chain and for inhomogeneous processes of semi-Markov type defined on a superposition of two renewal processes. In some cases, the relations for the desired distributions are simplified and make it possible to find limit distributions of functionals (in particular, distributions of functionals related to crossing a level $x$ when $x\to 0$ or $x\to \pm\infty$, distributions of absolute extrema or their complements). For sums of random variables on a Markov chain, factorization identities were first considered in the works of foreign authors, in particular, G. Baxter [71] and H. Miller [102]. In the works of E.L. Presman [47], a probabilistic interpretation of the factorization components is given, with the help of which the asymptotic behavior of the distribution of functionals of sums was investigated.

One of the first objects of study for homogeneous processes with independent increments on a Markov chain was the central limit theorem (G. Alyashkyavichyus [1], J. Keilson, D. Wishart [95], G. Fukuschima, M. Hitsuda [85]). A complete description of homogeneous processes with independent increments on a Markov chain $Z(t)=\{\xi(t),x(t)\}$ is given in the work of I.I. Ezhov and A.V. Skorokhod [20], where these processes were called Markov processes homogeneous in the corresponding component.

Among the works devoted to walks and processes on a Markov chain, one should also mention some works of Novosibirsk mathematicians: A.A. Mogulsky [37], where the construction of factorization identities is first carried out for sums on a Markov chain, and then the corresponding identities for processes on a Markov chain are established by passing to the limit; K. Arndt [3], where the distribution of extrema of a random walk on a Markov chain was investigated and the asymptotics of these distributions were analyzed; V.S. Lugavov [36], where the distribution of the sojourn time in a strip of processes with independent increments on a Markov chain was studied.

The multidimensional process of inventory evolution in a Markovian environment (which is very close to risk, since the probability of ruin is analogous to the probability of inventory depletion) was studied in [14].

In the works of Jansen [89], Jansen and Reinhard [90], Reinhard [117], Asmussen [66], Schmidli [118], and others, risk processes in a Markovian environment were considered, which are defined on a finite Markov chain $x(t), t\ge 0$ with phase space $X=\{1,2,\dots,m\}$ and a set of lower semicontinuous processes $\xi_k(t)=S_k(t)-t$. The jumps of the processes occur with intensity $\lambda_k\ (k=1,m)$ and have a continuous distribution $F_k(x),\ F_k(0)=0$.

The two main results are as follows. This is the analogue of the Cramér-Lundberg approximation for the ruin probability obtained by Asmussen [66, 67]:
\begin{equation}
\psi_i(u)\approx C_i e^{-\gamma u},\nonumber
\end{equation}
where $(u,i)$ is the initial state of the process, and Asmussen requires that the tails of the distributions $F_i(x)$ decay exponentially; the second result is the algorithm of Asmussen and Rolski [68] for computing $\psi_i(u)$ in the case of phase-type distributions $F_i(x)$. However, these distributions do not contain heavy-tailed distributions such as log-normal, Pareto, or log-gamma distributions. Denote $\overline{F}(u)=\int_u^{+\infty} F(du)$ the tail of the measure, $\mu_F$ the mean of the claim distribution $F$, and $F_0(u)=\overline{F(u)}/\mu_F$. The classical result for the Poisson model of the risk process with heavy-tailed claims (without Markov modulation) is that the ruin probability $\psi(u)\approx C\overline{F_0(u)}$, where $C$ is some constant (see von Bahr [119], Embrechts and Veraverbeke [79], Embrechts and Vilaseñor [80], Klüppelberg [98], Zinchenko [120]). This formula holds for subexponential distributions $F$, i.e., $\overline{F^{*2}}(u)/\overline{F}(u)\to 2$, $u\to\infty$. In the work of Asmussen, Henriksen, and Klüppelberg [69], this classical result is extended to Markov-modulated risk processes. Namely, suppose that there exists a subset of states $E_F$ of the chain for which the claim distributions $F_i$, $i\in E_F$ satisfy $(1-F_i(x))/(1-F(x))\to b_i\in(0,+\infty)$, $x\to\infty$, for some subexponential function $F$, and $(1-F_i(x))/(1-F(x))\to 0$, $x\to\infty$, for $i\notin E_F$. Then it is shown that the ruin probability $\psi_i(u)\approx C_i\int_u^\infty (1-F(x))\,dx$, where $C_i$ are some constants.

In the work of Schmidli [118], an overview of results on risk processes, in particular, on risk processes in a Markovian environment perturbed by a Wiener process, is given. In [122], an integral equation is derived for a risk process with a stochastic return on investment rate.

\section{Numerical methods, the method of successive approximations}

A number of methods for estimating the ruin probability are based on the Pollaczek-Khinchine formula [65] for solving equation (23). In particular, in [38], an effective variant of the Monte Carlo method is proposed; in [24], using the same approach, good upper and lower bounds are obtained in the presence of large claims; and in [111], an iterative method of approximations from above and below using a discrete approximation of the claim distribution function is proposed.

Note that (\ref{eqn:1.2}) is a Volterra integral equation. Indeed, making the change of variables $x=u-z$ in the integral, instead of (\ref{eqn:1.2}), we obtain:
\begin{equation}
\varphi(u)=1-\frac{\alpha\mu}{c}+\frac{\alpha}{c}\int_0^u \varphi(x)K(x,u)\,dx,\label{eqn:1.4}
\end{equation}
where the kernel $K(x,u)=1-F(u-x)$ is a measurable bounded function. As is known [27], the right-hand side of the Volterra equation (\ref{eqn:1.2}) is a contraction operator; therefore, it has a unique solution, which can be found by the method of successive approximations. However, equations (\ref{eqn:1.2}), (\ref{eqn:1.4}) are Volterra equations of a special kind; therefore, their solutions and successive approximations to the solution have specific properties.

In [13, 41-44, 77], the method of successive approximations was used to find the probability of ruin. In [13, 41], this method was used to obtain estimates of the probability of non-ruin and to solve integral equations of actuarial mathematics of Volterra type. In [77], numerical examples of solving actuarial integral equations by the Runge-Kutta method and the method of successive approximations (without justification) are given. In [41-44], conditions for the existence and uniqueness of solutions of general actuarial integral equations were obtained, the convergence of the method of successive approximations was proved, and the results of numerical experiments were given.

The XTREMALS 3.01 package, developed by German scientists [112], is a system for processing and visualizing statistical data. In relation to risk analysis, the package allows finding the ruin probability over a finite interval for the classical risk process, provided that the claim size has either a Poisson distribution $F(x)=1-e^{-x},\ x\ge 0$, or a Pareto distribution $F(x)=1-x^{-\alpha},\ x\ge 1,\ a>0$, or a beta distribution $F(x)=1-(-x)^{-\alpha},\ -1<x\le 0,\ a<0$. Calculations are performed by the Monte Carlo method. The output is the ruin probability (numerically and graphically) for initial capitals from zero to a given level with step 1.

\section{Conclusions to the section}

As can be seen, exact theoretical methods for solving the problem of finding the ruin probability are diverse and sufficiently developed. However, their practical application to general problems is extremely difficult. Even the calculation of the ruin probability for the special case of a process on a Markov chain with two states [40] is very laborious and raises doubts about the possibility of practical application of exact methods in the general situation. There are also a number of exact numerical methods for finding the ruin probability of the classical risk process, based on the fact that this probability satisfies some Volterra integral equation. Because of this, their generalization to more general processes is difficult.

On the other hand, empirical analytical approximations are widely used in practice, but they are not always good and there are no guaranteed estimates of their accuracy.

A general numerical method is the Monte Carlo simulation method, but it has poor accuracy in the most interesting case of small ruin probabilities.

An alternative to these methods can be the method of successive approximations, which guarantees any given accuracy of the solution and is very practical, as it is applicable to a wide class of processes and has good convergence rate characteristics.

This dissertation is devoted to the study of various generalized risk processes, the study of integral equations for the ruin probability as a function of initial capital, the substantiation of a general method of successive approximations for finding the probability of non-ruin, and the study of the convergence rate of the method.

\chapter{MODELS OF INSURANCE ACTIVITY AND INTEGRAL EQUATIONS FOR THE PROBABILITY OF RUIN}

In this section, various mathematical models of insurance activity are considered that generalize the classical collective risk model of F. Lundberg [99]. For each model, the corresponding integral equations for the probability of ruin of an insurance company as a function of initial capital are written out. The models and equations of this section are published in [40-47], and the general methodology for their study in [106-108].

\section{Integral equations of the theory of the classical risk process}

The classical risk process describing the evolution over time of the capital $\xi_t$ of an insurance company is given by the relation (see [35, 51, 65, 72]):
\begin{equation}
    \xi _t = u + ct - S_t,\label{(2.1)}
\end{equation}
where $t$ is time; $u$ is the initial reserve of the insurance company; $c$ is the (constant) premium income rate; $S_t$ is the aggregated claim payments by time $t$, $S_t=\sum_{k=1}^{N_t} z_k$; $z_k$ are independent identically distributed random variables (losses) with distribution function $F(z)$ and mean $\mu$, $F(z)=0$ for $z\le 0$; $N_t$ is the number of payments by time $t$ (a Poisson process with intensity $\alpha$). The theory of such processes is studied in detail and constitutes an important part of actuarial mathematics [35, 51, 65, 72]. As is known, the probability of non-ruin on an infinite time interval $\varphi (u)=\Pr\{\xi _t\ge 0\ \forall t\ge 0\}$ with initial capital $u$ satisfies the integral equation [35, p. 227]:
\begin{equation}
    \varphi (u)=\int\limits_{0}^{+\infty }{\alpha e^{-\alpha t}\int\limits_{0}^{u+ct}{\varphi (u+ct-z)dF(z)dt}}\label{(2.2)}
\end{equation}
and the boundary conditions: $\varphi (u)=0$ for $u<0$ and $\varphi (+\infty)=1$.

Equation (\ref{(2.2)}) with boundary condition $\phi(+\infty)=1$ can be rewritten in the form [35, (3.13)]:
\begin{equation}
    \varphi (u)=1-\frac{\alpha \mu }{c}+\frac{\alpha }{c}\int\limits_{0}^{u}{\varphi (u-z)\left( 1-F(z) \right)dz}.\label{(2.3)}
\end{equation}

Note that (\ref{(2.3)}) is a Volterra integral equation. Indeed, making the change of variables $x=u-z$ in the integral, instead of (\ref{(2.3)}) we get:
\begin{equation}
    \varphi (u)=1-\frac{\alpha \mu }{c}+\frac{\alpha }{c}\int_{\,0}^{\,u}{\varphi (x)K(x,u)dx},
		\label{(2.4)}
\end{equation}
where the kernel $K(x,u)=1-F(u-x)$ is a measurable bounded function. A number of approaches to its solution are known (see, for example, [11, 12, 30, 35, 38, 65]).

In the cases of exponential and fixed claims, exact solutions are obtained using the Laplace transform. The Cramér-Lundberg estimate is also known: $\varphi (u)\ge 1-e^{-Lu},\ u\ge 0$, where the Lundberg coefficient $L$ satisfies the equation:
\begin{equation}
    \frac{\alpha}{c}\int_{\,0}^{\,+\infty }{e^{Lz}\left[ 1-F(z) \right]dz}=1.\label{(2.5)}
\end{equation}
A number of approximations of $\varphi(u)$ for large $u$ are also known (see, for example, [35, 72]).

However, the problem of analytical and numerical solution of equation (\ref{(2.3)}) continues to attract the attention of researchers.

As is known [27], the Volterra equation has a unique solution, which can be found by the Picard method of successive approximations. But (\ref{(2.3)}) is a special Volterra equation, so the corresponding method of successive approximations has features. These issues are addressed in Chapter 3.

\section{Nonlinear risk process with a Poisson stream of independent claims}

Consider the risk process with nonlinear capital growth
\begin{equation*}
    \xi _t=u+\int_{\,0}^{\,t}{c(\xi _s)ds}\ -S_t,\quad t\ge 0,
\end{equation*}
where the premium arrival rate $c(\cdot )\ge 0$ depends on the current capital of the company and is a piecewise continuous function. This is a piecewise-deterministic Markov process [55]. We will assume that the reserves (capital) $U$ of the insurance company in the absence of insurance claims change over time $t$ according to the differential equation $dU/dt=c(U)$, $U(0)=u$. The solution of this Cauchy problem depends on $u$ and $t$; denote it by $U(u,t)$, $U(u,0)=u$.

For example, in the classical case $U(u,t)=u+ct$. If, in addition, insurance premiums are withdrawn from reserves starting from some reserve level $v\ge u$, then $U(u,t)=\min\{v,u+ct\}$.
If reserves are kept in a bank account with a continuous interest rate $\delta$, and insurance premiums $c=c(U)$ depend on the current level of reserves, then $U(t,u)$ satisfies the differential equation almost everywhere:
\begin{equation*}
    \frac{dU}{dt}=\delta U+c(U),\quad U(0)=u.
\end{equation*}
In particular, if $c=const$, then
\begin{equation*}
    U(u,t)=\left( u+\frac{c}{\delta} \right)e^{\delta t}-\frac{c}{\delta}=ue^{\delta t}+\frac{c}{\delta}(e^{\delta t}-1).
\end{equation*}

Insurance contracts can be of two main types: proportional and with limits on payment amounts. Let $x$ be a random loss in an individual insurance event. In general, payments (claim function) $Z(x)\le x$ under an insurance contract can be represented as
\begin{equation*}
    Z(x)=\min\{ b,\gamma \max\{0,x-a\} \},
\end{equation*}
where $a$ is the amount of loss from which the insurance contract comes into effect; $\gamma$ is the insured share of the loss exceeding the initial level $a$; $b$ is the maximum level of insurance payments. Let us introduce the distribution function of insurance claims $F(z)=P\{Z(x)<z\}$.

It is often assumed that insurance events occur independently and are distributed in time according to some probability distribution, usually the Poisson law with parameter $\alpha$. General risk processes with non-Poisson streams of dependent insurance claims are considered in [21].

We will also assume that the company maintains a constant composition of the insurance portfolio, i.e., instead of departing clients, it immediately recruits new subscribers with similar parameters of insurance contracts.

The most important parameter of the activity of an insurance company is the probability of insolvency $\psi(y)$ as a function of the initial reserve $u$ and other parameters of the insurance policy, i.e., $y=\{u,\alpha,c,F,U(\cdot,\cdot),Z(\cdot)\}$ (often the dependence of $\psi(\cdot)$ only on the initial reserve $u$ is indicated).

The function $\psi(u)$ is used to select the insurance reserve $u$ that guarantees a certain level of company stability, i.e., a certain probability of non-ruin $\varphi(u)=1-\psi(u)$ under a given insurance policy [72]. Conversely, the function $\varphi(u)$ can be used to select the amount of the insurance premium $c$ for given $u,\alpha,F$, which guarantees a certain level of stability.

The functions $\varphi(\cdot),\psi(\cdot)$ can be considered as risk measures for the solution (insurance policy) $y=\{u,\alpha,c,F,U(\cdot,\cdot),Z(\cdot)\}$, which can be chosen from some finite or infinite set of alternatives $Y$. Another important indicator of the activity of an insurance company is some measure of profitability $\Phi(y)$ of policy $y$, which can be the amount of the insurance premium $c$, the insurance loading $\rho=c/(\alpha\mu)-1$ (showing the excess of premiums over the average costs of covering insurance claims), or the mathematical expectation of insurance reserves reduced to the initial moment, dividends, etc. Thus, each possible solution $y$ is described by profitability indicators $\Phi(y)$ and risk $\varphi(y)$, and the entire set of possible alternatives $Y$ can be represented as a certain set in the "profitability-risk" space, as is usually done in portfolio investment theory. The decision-maker in the insurance company only has to choose a suitable point in this set.

It is easy to see that for the probability of non-ruin $\varphi(u)$ of an insurance company with initial reserve $u$ on an infinite time interval, the following integral equation holds (compare with [35, (3.12)]):
\begin{equation}
    \varphi (u)=\int\limits_{0}^{+\infty }{\alpha e^{-\alpha t}\int\limits_{0}^{U(u,t)}{\varphi (U(u,t)-z)dF(z)dt}}\label{(2.6)}
\end{equation}
with boundary conditions $\varphi(u)=0$ for $u<0$ and $\varphi(+\infty)=1$, and the probability of non-ruin $\varphi(u,T)=\Pr\{\xi_t\ge 0\ \forall t\in[0,T]\}$ on the finite time interval $[0,T]$ is determined by the equation
\begin{equation}
    \varphi (u,T)=\int\limits_{0}^{T}{\alpha e^{-\alpha t}\int\limits_{0}^{U(u,t)}{\varphi (U(u,t)-z,T-t)dF(z)dt}}\label{(2.7)}
\end{equation}
with boundary conditions $\varphi(u)=0$ for $u<0$ and $\varphi(+\infty,T)=1\ \forall T\ge 0$. Equations (\ref{(2.6)}) and (\ref{(2.7)}) in the general case do not reduce to the Volterra integral equation (\ref{(2.3)}). Equation (\ref{(2.6)}) may have no solutions satisfying the boundary condition $\varphi(+\infty)=1$. For example, its simplified version (\ref{(2.2)}) has such a solution only under the condition $\alpha\mu/c<1$. Otherwise, equation (\ref{(2.2)}) has only the solution identically equal to zero. This model is considered in Section 5.

\section{Risk process in the presence of random premiums}

Consider a risk process $\xi_t$ of the form
\begin{equation*}
    \xi _t = u + ct + P_t - S_t,\quad t\ge 0,
\end{equation*}
$S_t=\sum_{k=1}^{N_t^c} z_k$ are the aggregated insurance claims, $P_t=\sum_{k=1}^{N_t^p} p_k$ are the aggregated premiums, $z_k$ are independent random variables (claims) with common distribution function $F(z)$, $p_k$ are independent random variables (premiums) with common distribution function $G(z)$, $N_t^c$ is the number of claims received by time $t$ – in this model, it is a Poisson process with intensity $\alpha_c$. Similarly, $N_t^p$ is the number of premiums received by time $t$, a Poisson process with intensity $\alpha_p$. This model describes the process of premium arrival more adequately.

{\bf Theorem 2.1.}
The probability of non-ruin $\varphi(u)$ satisfies the integro-differential equation
\begin{equation}
    (\alpha_c+\alpha_p)\varphi(u)-c\varphi'(u)=\alpha_c\int_{0}^{u}{\varphi(u-z)dF(z)}+\alpha_p\int_{0}^{+\infty}{\varphi(u+z)dG(z)}.\label{(2.8)}
\end{equation}

This model is studied in detail in Section 6.

\section{Reinsurance model}

Consider a mathematical model of a reinsurance company with initial reserve $u_0$, reserve growth function $U_0(u_0,t)$. Assume that companies $i=1,\ldots,n$ with initial reserve $u_i$, reserve growth function $U_i(u_i,t)$, Poisson claim stream with intensity parameter $\alpha_i$, and claim distribution function $F_i(z)$ are subject to reinsurance. We will assume that all streams of insurance events and claim sizes are independent. Then the probability of non-ruin $\varphi(u_0)$ of the reinsurance company on an infinite time interval as a function of the variable $u_0$ (and parameters $u_1,\ldots,u_n$) similarly (\ref{(2.6)}) satisfies the equation:
\begin{equation}
    \varphi (u_0)=\int\limits_{0}^{+\infty }{\sum\limits_{i=1}^{n}{\alpha_i e^{-\alpha_i t}\int\limits_{0}^{U_0(u_0,t)}{\varphi (U_0(u_0,t)-z)dF_i(z)dt}}},\label{(2.9)}
\end{equation}
with boundary conditions $\varphi(u_0)=0$ for $u_0<0$ and $\varphi(+\infty)=1$, and the probability of ruin $\varphi(u_0,T)$ on the finite time interval $[0,T]$ satisfies the equation:
\begin{equation}
    \varphi (u_0,T)=\int\limits_{0}^{T}{\sum\limits_{i=1}^{n}{\alpha_i e^{-\alpha_i t}\int\limits_{0}^{U_0(u_0,t)}{\varphi (U_0(u_0,t)-z,T-t)dF_i(z)dt}}}\label{(2.10)}
\end{equation}
with boundary conditions $\varphi(u_0)=0$ for $u_0<0$ and $\varphi(+\infty,T)=1\ \forall T\ge 0$.

Reinsurance models are studied, for example, in [15, 75].

\section{Risk process with a non-Poisson flow of independent claims}

Let $\tau_1,\tau_2,\tau_3,\ldots$ be successive claim arrival times, where $\tau_1$ has distribution function $K_1(t)$, and the intervals $\tau_i-\tau_{i-1},\ i>1$, have the same distribution function $K(t)$ (renewal process). Denote by $\varphi(u)=\Pr\{\xi_t\ge 0\ \forall t\ge 0\}$ the probability of non-ruin of the process $\{\xi_t\}_{t\ge 0}$ on an infinite time interval with initial capital $u\ge 0$. This function satisfies the following integral equations (see [35, Section III]).

\begin{itemize}
    \item{a)} For the ordinary renewal risk process, we have $K_1(t)=K(t)$ and (see [35, (3.74)]):
    \begin{equation}
        \varphi (u)=\int\limits_{0}^{+\infty }{dK(t)\int\limits_{0}^{u+ct}{\varphi (u+ct-z)dF(z)}},\quad \varphi(+\infty)=1;\label{(2.11)}
    \end{equation}
    \item{b)} For the general renewal risk process (delayed renewal process), $K_1(t)\ne K(t)$, and its probability of non-ruin $\varphi_1(u)$ is expressed in terms of the probability of non-ruin $\varphi(\cdot)$ of the ordinary renewal risk process [35, (3.72)]:
    \begin{equation}
        \varphi_1(u)=\int\limits_{0}^{+\infty }{dK_1(t)\int\limits_{0}^{u+ct}{\varphi (u+ct-z)dF(z)dt}};\label{(2.12)}
    \end{equation}
    \item{c)} For the special case of a stationary renewal risk process, the distribution functions $K_1(t)$, $K(t)$ are related by $K_1(t)=\alpha\int_{0}^{t}[1-K(\tau)]d\tau$, where $\alpha=1/\int_{0}^{+\infty}(1-K(\tau))d\tau$, and its probability of non-ruin $\varphi_1(u)$ is expressed in terms of the probability of non-ruin $\varphi(\cdot)$ of the ordinary renewal risk process via (\ref{(2.12)}) or [35, (3.75)]:
    \begin{equation*}
        \varphi_1(u)=1-\frac{\alpha \mu}{c}+\frac{\alpha}{c}\int_{\,0}^{\,u}{\varphi (u-z)\left( 1-F(z) \right)dz}.
    \end{equation*}
    \item{d)} For the special case of the classical risk process with a Poisson stream of claim occurrence times $\tau_0,\tau_1,\tau_2,\ldots$, we have $K_1(t)=K(t)=1-e^{-\alpha t}$, and (\ref{(2.11)}) holds, and the probability of non-ruin $\varphi(u)$ satisfies the equation [35, (3.13)]:
    \begin{equation}
        \varphi (u)=1-\frac{\alpha \mu}{c}+\frac{\alpha}{c}\int_{\,0}^{\,u}{\varphi (u-z)\left( 1-F(z) \right)dz},\quad \varphi(+\infty)=1.\label{(2.13)}
    \end{equation}
\end{itemize}

The non-Poisson renewal risk process is studied in Chapter 4.

\section{Mixed risk process}

Let there be a Poisson stream (intensity $\alpha$) of insurance events with loss distribution $F(z)$. Let there also be another independent stream of special events with random intervals $\tau$ between events, $K(t)=P\{\tau<t\}$, and with loss size distribution $H(z)$. Denote by $\varphi_t(u)$ the probability of ruin of such a mixed risk process if the initial capital is $u$ and time $t$ has passed since the last special event. We also introduce the conditional distribution function
\begin{equation*}
    K_t(s)=\frac{P\{\tau>t,\tau>t+s\}}{P\{\tau>t\}}=\frac{P\{\tau>t+s\}}{P\{\tau>t\}}=\frac{1-K(t+s)}{1-K(t)}
\end{equation*}
as the conditional probability that a special claim will not arrive for another time $s$ if it has not occurred during time $t$. Define the conditional density
\begin{equation*}
    k_t(s)=- [K_t(s)]'_s.
\end{equation*}
Then the integral equation for $\varphi_t(u)$ has the form (compare with [35, (3.72)]):
\begin{align}
    \varphi_t(u) & = \int\limits_{0}^{\infty }{ \left\{ k_t(s) e^{-\alpha s} \int\limits_{0}^{U(u,s)} \varphi_0(U(u,s)-z) dH(z) \right. } \nonumber \\
    & \quad + \left. K_t(s) e^{-\alpha s} \alpha \int\limits_{0}^{U(u,s)} \varphi_{t+s}(U(u,s)-z) dF(z) \right\} ds, \label{(2.14)} \\
    & \varphi_t(u)=0,\quad \varphi_t(+\infty)=1,\quad \forall u<0,\ t\ge 0.\nonumber
\end{align}

\section{Risk process in a Markovian environment}

Consider a risk process $\xi_t$ in a Markovian environment (on a Markov chain) with $n$ states $\{1,2,\ldots,n\}$ [18, 40, 88, 91]:
\begin{align*}
    \xi_t & = u + \left( \sum_{h=0}^{M_t-1} c(h)(\tau_{h+1}-\tau_h) + c(M_t)(t-\tau_{M_t}) \right) \nonumber \\
    & \quad - \left( \sum_{h=0}^{M_t-1} \sum_{k=1}^{N(\tau_{h+1}-\tau_h,\alpha(h))} Y(k,h) + \sum_{k=1}^{N(t-\tau_{M_t},\alpha(M_t))} Y(k,M_t) \right),
\end{align*}
where $u$ is the initial capital of the company; $h$ is the ordinal number of the environment state switch; $M_t$ is the number of environment state switches by time $t$; $N(\Delta,\alpha)$ is the number of insurance events during time $\Delta$ with intensity $\alpha$; $\alpha_i$, $\alpha(h)$ are the intensities of insurance events in state $i$ and after the $h$-th switch, respectively; $c_i$ and $c(h)$ are the premium income rates in state $i$ and after the $h$-th switch, respectively; $\tau_h$ is the $h$-th moment of change of the chain state, $\tau_0=0$; $Y(k,h)$ is the $k$-th claim after the $h$-th switch.
The probability that the chain remains in state $i$ for time $\Delta$ is $e^{-\lambda_i\Delta}$, the probability of the chain transitioning from state $i$ to state $j\ne i$ during time $\Delta$ is $\lambda_i p_{ij}\Delta+o(\Delta)$, where $\sum_{j=1}^{n} p_{ij}=1$, $p_{ii}=0$, $\lim_{\Delta\to 0} o(\Delta)/\Delta=0$.

\begin{figure}[h]
    \centering
    \includegraphics[width=0.8\textwidth]{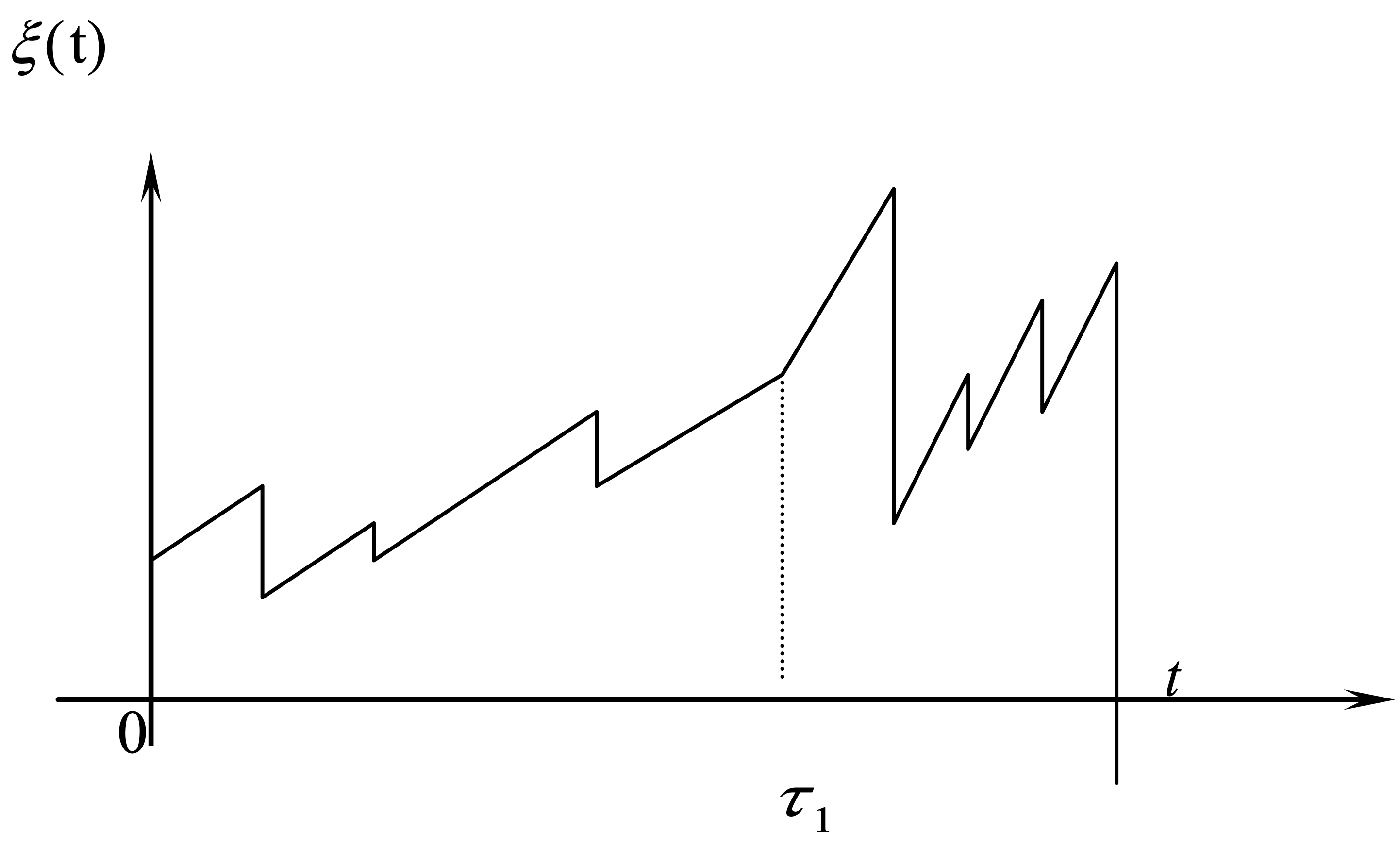}
    \caption{Trajectory of the risk process}
    \label{fig:3.1}
\end{figure}

Consider the probability $\psi_i(u)$ of ruin from the $i$-th initial state of the chain as a function of initial capital $u$ and the probability that ruin does not occur when starting from the $i$-th state, $\varphi_i(u)=1-\psi_i(u)$.

{\bf Theorem 2.2} [40]. The vector-function $\varphi(u)=\{\varphi_1(u),\varphi_2(u),\ldots,\varphi_n(u)\}$ satisfies the system of integral equations ($i=1,\ldots,n$):
    \begin{equation}
        \varphi_i(u)=\int_{\,0}^{\,\infty }{e^{-t(\alpha_i+\lambda_i)}
        \left( \lambda_i \sum_{j=1}^{n} p_{ij} \varphi_j(U_i(u,t))
        + \alpha_i \int\limits_{0}^{U_i(u,t)} \varphi_i(U_i(u,t)-z) dF_i(z) \right) dt},
				\label{(2.15)}
    \end{equation}
    where $U_i(u,t)=u+c_i t$. Obviously, (\ref{(2.15)}) also holds for general reserve growth functions $U_i(u,t)$ in environment state $i$.

This model is studied in detail in Chapter 7.

\section{Risk process in a semi-Markovian environment}

Let the environment be in $i=1,2,\ldots,n$ states, and let $\tau_i$ be the random time spent by the environment in the $i$-th state. Denote by $K_i(t)=P\{\tau_i<t\}$ the distribution function of $\tau_i$. Introduce the conditional distribution function
\begin{equation*}
    K_{it}(s)=\frac{P\{\tau_i>t,\tau_i>t+s\}}{P\{\tau_i>t\}}=\frac{P\{\tau_i>t+s\}}{P\{\tau_i>t\}}=\frac{1-K_i(t+s)}{1-K_i(t)}
\end{equation*}
as the conditional probability that the environment will remain in state $i$ for at least another time $s$ if it has already been in this state for time $t$; introduce also the corresponding density
\begin{equation*}
    k_{it}(s)=- [K_{it}(s)]'_s.
\end{equation*}
Denote by $\varphi_{it}(u)$ the probability that the process will not be ruined from the initial state $(i,t,u)$, i.e., under the condition that the environment has been in state $i$ for time $t$, and the initial capital is $u$. Then, similarly to [40], the functions $\varphi_{it}(u)$ satisfy the system of integral equations:
\begin{align}
    \varphi_{it}(u) & = \int\limits_{0}^{\infty }{ \left\{ \sum_{j\ne i} k_{it}(s) e^{-\alpha_i s} \pi_{ij} \varphi_{j0}(U_i(u,s)) \right. } \nonumber \\
    & \quad + \left. K_{it}(s) e^{-\alpha_i s} \alpha_i \int\limits_{0}^{U_i(u,s)} \varphi_{i(t+s)}(U_i(u,s)-z) dF_i(z) \right\} ds,\label{(2.16)}  \\
    & \varphi_{it}(u)=0,\quad \varphi_{it}(+\infty)=1,\quad \forall i,u<0,\ t\ge 0,\quad i=1,\ldots,n,\nonumber
\end{align}
where $U_i(u,t)$ are the reserve growth functions in state $i$.

\section{Conclusions to the section}

The integral equations of insurance mathematics for the probability of non-ruin $\varphi(u)$, $0\le \varphi(u)\le 1$, as a monotone function of the initial capital $u$ have the following general form:
\begin{equation}
    \varphi (u)=A\varphi (u),\label{(2.17)}
\end{equation}
where $A$ is a certain linear integral operator, and the solution of the integral equation must satisfy the boundary condition
\begin{equation}
    \lim_{u\to\infty}\varphi(u)=1.\label{(2.18)}
\end{equation}
The dissertation proposes to solve all problems of actuarial mathematics of the form (\ref{(2.17)}), (\ref{(2.18)}) by a unified method, namely, the method of successive approximations:
\begin{equation}
    \varphi ^{k+1}(u)=A\varphi ^{k}(u),\quad 0\le \varphi ^{0}(u)\le 1,\quad k=0,1,\ldots
		\label{(2.19)}
\end{equation}
It turns out that, generally speaking, the operator $A$ is not a contraction on its domain of definition, and the problem consists in justifying the method of successive approximations.

In subsequent sections, the properties of the operator $A$ for specific types of risk processes, conditions for the existence and uniqueness of the solution of problems (\ref{(2.17)}), (\ref{(2.18)}), the nature of successive approximations, and conditions for the convergence of $\varphi^k(u)$ to the solution of problems (\ref{(2.17)}), (\ref{(2.18)}) are investigated. The main point of the study of problems (\ref{(2.17)}), (\ref{(2.18)}) and method (\ref{(2.19)}) is to find monotone functions $\varphi_*(u)$, $\varphi^*(u)$ such that $0\le \varphi_*(u)\le \varphi^*(u)\le 1$ and
\begin{equation*}
    A\varphi_*(u)\ge \varphi_*(u),\quad A\varphi^*(u)\le \varphi^*(u)\ \forall u\ge 0.
\end{equation*}

It turns out that one can always take $\varphi^*(u)\equiv 1$, and as $\varphi_*(u)$ the generalized Cramér-Lundberg bound. Moreover, the existence of functions $\varphi_*(u)$, $\varphi^*(u)$ is necessary and sufficient for the existence of a solution to problems (\ref{(2.17)}), (\ref{(2.18)}), and the successive approximations starting with $\varphi^0(u)=\varphi^*(u)$ converge monotonically from above, and those starting with $\varphi^0(u)=\varphi_*(u)$ converge monotonically from below to the solution of problems (\ref{(2.17)}), (\ref{(2.18)}). The uniqueness of the solution of problems (\ref{(2.17)}), (\ref{(2.18)}) is proved for each equation separately using the properties of the corresponding operator $A$.

\chapter{ON THE METHOD OF SUCCESSIVE APPROXIMATIONS FOR CALCULATING THE PROBABILITY OF RUIN OF THE CLASSICAL RISK PROCESS}

In this section, the Picard method of successive approximations for solving the renewal integral equation (of Volterra type), which is satisfied by the probability of (non-)ruin of the classical risk process, is investigated. Convergence, monotonicity, and the rate of convergence of the approximations are established in the entire domain of possible initial values of the risk process. The results are illustrated by numerical calculations. The results of the section are published in [41, 42, 44, 104].

\section{Introduction}

Recall that the classical risk process describing the evolution of the capital of an insurance company is given by the relation
\begin{equation}
    \xi _t = u + ct - S_t,\label{(3.1)}
\end{equation}
where $t$ is time; $u$ is the initial capital of the insurance company; $c$ is the premium income rate; $S_t$ is the aggregated payments by time $t$, $S_t=\sum_{k=1}^{N_t} Y_k$; $Y_k$ are independent identically distributed random variables (claims) with distribution function $F(y)$ and mean $\mu$; $N_t$ is the number of payments by time $t$ (a Poisson process with intensity $\alpha$). As is known, the probability of non-ruin $\varphi(u)=\Pr\{\xi_t\ge 0\ \forall t\ge 0\}$ with initial capital $u$ satisfies the renewal integral equation
\begin{equation}
    \varphi (u)=1-\frac{\alpha \mu}{c}+\frac{\alpha}{c}\int_{\,0}^{\,u}{\varphi (u-z)\left( 1-F(z) \right)dz}.\label{(3.2)}
\end{equation}

However, equations (\ref{(2.3)}), (\ref{(2.4)}) and (\ref{(3.2)}) are Volterra equations of a special type, so their solutions and successive approximations to the solution have specific properties. The purpose of this section is to study these properties.

As is known [27], the right-hand side of the Volterra equation (\ref{(3.2)}) is a contraction operator; therefore, it has a unique solution, which can be found by the Picard method of successive approximations.
\begin{equation}
    \varphi ^{k+1}(u)=1-\frac{\alpha \mu}{c}+\frac{\alpha}{c}\int\limits_{0}^{u}{\varphi ^{k}\left( u-z \right)\left[ 1-F\left( z \right) \right]dz},\quad k=0,1,\ldots,\label{(3.3)}
\end{equation}
where $\varphi^0(u)$ is some initial function. The method of successive approximations also allows finding analytical approximate solutions for the probability $\varphi$ as a function of $u$ and other parameters of the process (\ref{(3.3)}). As is known, Picard's method for Volterra integral equations converges on every finite interval of values of $u$ by virtue of the contraction mapping principle. However, in general, the method is not monotone, and contraction occurs only after a certain number of iterations; therefore, the guaranteed convergence of the method for large $u$ may be slow. The method of successive approximations for solving equation (\ref{(3.2)}) was proposed to be applied in [13], see also [6, Chap. 5, §4], where for the case of a continuous claim distribution function $F(x)$, its convergence and monotonicity were established. In this section, in the general case, the properties of successive approximations for finding the probability of non-ruin $\varphi(u)$ of the classical risk process are investigated in detail. It is shown that for the Volterra equation (\ref{(3.2)}) under consideration, contraction occurs at each iteration with coefficient $p=\alpha\mu/c<1$ uniformly over all $u\in[0,+\infty)$, which ensures a high convergence rate for all $u$. It is shown that when starting from the initial functions $\varphi^0(u)\equiv 1$ and $\varphi^0(u)\equiv \varphi(0)=1-\alpha\mu/c$, the successive approximations converge monotonically to the solution $\varphi(u)$. It is shown that when starting from the initial function $\varphi^0(u)\equiv 1$, all successive approximations $\varphi^k(u)$ satisfy the condition $\varphi^k(u)\ge 1-e^{-Lu}$, where $L>0$ is the Lundberg constant, and thus the boundary conditions $\varphi^k(+\infty)=1$ are satisfied at each iteration and in the limit, $\varphi(+\infty)=1$.

\section{The method of successive approximations}

The solution of equation (\ref{(3.2)}) can be found by the following method of successive approximations:
\begin{equation}
    \varphi ^{k+1}(u)=1-\frac{\alpha \mu}{c}+\frac{\alpha}{c}\int_{\,0}^{\,u}{\varphi ^{k}\left( u-z \right)\left[ 1-F\left( z \right) \right]dz},\quad k=0,1,\ldots,\label{(3.4)}
\end{equation}
where $\varphi^0(u)$ is some initial function. All functions $\varphi^k(u)$ are considered to be defined on some interval $[0,v]$, $v<+\infty$. Note that to find $\varphi^{k+1}(u)$, it is sufficient to know $\varphi^k(u)$ on the interval $[0,v]$. Let $L_\infty([0,v])$ be the space of continuous functions $f(u)$ defined on the interval $[0,v]$ with norm $\|f\|=\max_{u\in[0,v]}|f(u)|$. As is known, $L_\infty([0,v])$ is a complete normed (Banach) space for any $v<+\infty$.

{\bf Theorem 3.1.}
    If $\frac{\alpha \mu}{c}<1$, $\int_{\,0}^{\,+\infty }(1-F(z))dz<\infty$, the function $\varphi^0(u),\ u\in[0,v]$, is monotonically non-decreasing and satisfies $0\le \varphi^0(u)\le 1$, then
    \begin{enumerate}
        \item[(i)] all $\varphi^k(u)$ are continuous and monotonically non-decreasing on the interval $[0,v]$, $0\le \varphi^k(u)\le 1$;
        \item[(ii)] the sequence $\{\varphi^k(u)\}$ is fundamental in $L_\infty([0,v])$, and it converges to the solution of the original equation, and $\|\varphi-\varphi^k\|\le \frac{p^k}{1-p}$, $p=\frac{\alpha \mu}{c}<1$;
        \item[(iii)] the limit function $\varphi(u)=\lim_{k\to\infty}\varphi^k(u)$ is continuous and monotone (non-decreasing), and $0\le \varphi(u)\le 1$.
    \end{enumerate}
%\end{theorem}

\begin{proof}
    (i). The function $1-F(z)$ is monotone in $z$, therefore Riemann integrable (see [39]), hence $\varphi(u)$ is continuous and monotonically non-decreasing. By induction, all $\varphi^k(u),\ k\ge 1$, are continuous and monotonically non-decreasing. Note that $\mu=\int_{0}^{+\infty}(1-F(z))dz$. If $0\le \varphi^k(u)\le 1$, then the inequalities hold:
    \begin{equation*}
        1-\frac{\alpha \mu}{c}\le \varphi^{k+1}(u)\le 1-\frac{\alpha \mu}{c}+\frac{\alpha}{c}\int_{0}^{u}{\varphi^k(u-z)[1-F(z)]dz}\le 1.
    \end{equation*}
    Thus, by induction, it follows that $1-\alpha\mu/c\le \varphi^{k+1}(u)\le 1,\ k\ge 1$.

    (ii). Let us estimate the convergence rate of the process of successive approximations. We have
    \begin{equation*}
        \varphi^{k+1}(u)=1-\frac{\alpha \mu}{c}+\frac{\alpha}{c}\int_{0}^{u}{\varphi^k(u-z)[1-F(z)]dz},
    \end{equation*}
    \begin{equation*}
        \varphi^{k}(u)=1-\frac{\alpha \mu}{c}+\frac{\alpha}{c}\int_{0}^{u}{\varphi^{k-1}(u-z)[1-F(z)]dz},
    \end{equation*}
    whence
    \begin{equation*}
        \varphi^{k+1}(u)-\varphi^{k}(u)=\frac{\alpha}{c}\int_{0}^{u}{(\varphi^k(u-z)-\varphi^{k-1}(u-z))[1-F(z)]dz},
    \end{equation*}
    \begin{align*}
        \|\varphi^{k+1}-\varphi^{k}\| & \le \frac{\alpha}{c}\int_{0}^{u}{\|\varphi^k-\varphi^{k-1}\|[1-F(z)]dz} \nonumber \\
        & \le \frac{\alpha}{c}\int_{0}^{+\infty}{[1-F(z)]dz}\cdot \|\varphi^k-\varphi^{k-1}\| \nonumber \\
        & \le \frac{\alpha \mu}{c}\cdot \|\varphi^k-\varphi^{k-1}\| = p\|\varphi^k-\varphi^{k-1}\|,\quad k=1,2,\ldots
    \end{align*}
    Then for $n\le m$
    \begin{align*}
        \|\varphi^n-\varphi^m\| & \le \sum_{k=m+1}^{n}\|\varphi^k-\varphi^{k-1}\| \le \|\varphi^{m+1}-\varphi^{m}\|\sum_{k=0}^{n-m-1}p^k \nonumber \\
        & \le \frac{\|\varphi^{m+1}-\varphi^{m}\|}{1-p} \le \frac{\|\varphi^1-\varphi^0\|p^m}{1-p} \le \frac{p^m}{1-p}.
    \end{align*}
    Consequently, the sequence $\{\varphi^k(u)\}$ is fundamental and has a unique continuous limit $\varphi(u)$ in $L_\infty([0,v])$ for each finite $v<+\infty$. Passing to the limit in (\ref{(3.4)}) as $k\to\infty$ for each fixed $u$, we get that $\varphi(u)$ satisfies the original equation (\ref{(3.2)}). Passing to the limit in the inequality $\|\varphi^n-\varphi^m\|=\max_{u\in[0,v]}|\varphi^n(u)-\varphi^m(u)|\le p^m/(1-p)$ as $n\to\infty$, we obtain $\|\varphi-\varphi^m\|\le p^m/(1-p)$.

    (iii). The limit function is continuous and monotone, as the limit of continuous and monotone functions such that $0\le \varphi^k(u)\le 1$.
\end{proof}

{\bf Lemma 3.1.}
(i). If ${{\varphi }^{0}}(u)\equiv \varphi (0)=1-{\alpha \mu }/{c}\;$, then the iterative sequence ${{\varphi }^{k}}(u)$ is monotonically non-decreasing (increasing).
(ii) If ${{\varphi }^{0}}(u)\equiv 1$, then $\left\{ {{\varphi }^{k}}\left( u \right) \right\}$ is monotonically non-increasing (decreasing).

\begin{proof}
  (i). For ${{\varphi }^{0}}(u)=1-{\alpha \mu }/{c}\;$ we have:
	\[{{\varphi }^{1}}(u)=1-\frac{\alpha \mu }{c}+\frac{\alpha }{c}\int\limits_{0}^{u}{{{\varphi }^{0}}\left( u-z \right)\left[ 1-F\left( z \right) \right]dz}=1-\frac{\alpha \mu }{c}+\frac{\alpha }{c}\int\limits_{0}^{u}{\left( 1-\frac{\alpha \mu }{c} \right)\left[ 1-F\left( z \right) \right]dz}=\]		\[=\left( 1-\frac{\alpha \mu }{c} \right)\cdot \left( 1+\frac{\alpha }{c}\int\limits_{0}^{u}{\left[ 1-F\left( z \right) \right]dz} \right)\ge 1-\frac{\alpha \mu }{c}={{\varphi }^{0}}(u)\ .\]
It holds:
\begin{equation}
{{\varphi }^{k+1}}(u)-{{\varphi }^{k}}(u)=\frac{\alpha }{c}\int\limits_{0}^{u}{\left( {{\varphi }^{k}}(u-z)-{{\varphi }^{k-1}}(u-z) \right)}\left[ 1-F(z) \right]dz.\label{(3.5)}
\end{equation}
Thus, obviously, if ${{\varphi }^{k}}(u)\ge {{\varphi }^{k-1}}(u)\ \ \forall u\in [0,v]$, then ${{\varphi }^{k+1}}(u)\ge {{\varphi }^{k}}(u)\ \ \forall u\in [0,v]$. The first statement has been proven.

(ii). Let ${{\varphi }^{0}}(u)\equiv 1$ , then 
\begin{align} 
& {{\varphi }^{1}}(u)=1-\frac{\alpha \mu }{c}+\frac{\alpha }{c}\int\limits_{0}^{u}{{{\varphi }^{0}}\left( u-z \right)\left[ 1-F\left( z \right) \right]dz}=1-\frac{\alpha \mu }{c}+\frac{\alpha }{c}\int\limits_{0}^{u}{\left[ 1-F\left( z \right) \right]dz}\le \nonumber\\ 
& \ \ \ \ \ \ \ \ \le 1-\frac{\alpha \mu }{c}+\frac{\alpha }{c}\int\limits_{0}^{+\infty }{\left[ 1-F\left( z \right) \right]dz}=1={{\varphi }^{0}}(u).\nonumber 
\end{align}
By induction, inequality (\ref{(3.5)}) implies that ${{\varphi }^{k}}(u)\ge {{\varphi }^{k-1}}(u)\ \ \forall u\in [0,v]$ , $k>0$. Thus, starting from the initial functions ${{\varphi }^{0}}(u)\equiv 1$ or ${{\varphi }^{0}}(u)=1-{\alpha \mu }/{c}$, the iterative process monotonically converges to the solution of the original integral equation.
\end{proof}

{\bf Lemma 3.2.}
If the initial function ${{\varphi }^{0}}(u)\equiv 1$ , then for all k ${{\lim }_{k\to \infty }}{{\varphi }^{k}}(u)=1$.

\begin{proof}
Let $u={{u}_{1}}+{{u}_{2}},\ \ {{u}_{1}}\ge 0,\ \ {{u}_{2}}\ge 0$. By Lemma 3.1, all ${{\varphi }^{k}}(u)$ are monotonically non-decreasing, therefore
\[\int_{0}^{u}{{{\varphi }^{k}}(u-z)\left[ 1-F(z) \right]dz}\ge \int_{0}^{{{u}_{1}}}{{{\varphi }^{k}}({{u}_{1}}+{{u}_{2}}-z)\left[ 1-F(z) \right]dz}\ge \]
\[\ge \int\limits_{0}^{{{u}_{1}}}{{{\varphi }^{k}}({{u}_{2}})\left[ 1-F(z) \right]dz}={{\varphi }^{k}}({{u}_{2}})\int\limits_{0}^{{{u}_{1}}}{\left[ 1-F(z) \right]dz}\ .\]
Suppose that ${{\lim }_{k\to \infty }}{{\varphi }^{k}}(u)=1$ . Let $u={{u}_{1}}+{{u}_{2}}\to \infty $ and ${{u}_{1}}\to \infty ,\ \ {{u}_{2}}\to \infty $. Then 
\begin{align} 
& 1\ge {{\lim }_{u\to \infty }}{{\varphi }^{k+1}}(u)={{\lim }_{u\to \infty }}\left( 1-\frac{\alpha \mu }{c}+\frac{\alpha }{c}\int\limits_{0}^{u}{{{\varphi }^{k}}(u-z)\left[ 1-F(z) \right]dz} \right)\ge \nonumber\\ 
& \ge 1-\frac{\alpha \mu }{c}+{{\lim }_{u\to \infty }}{{\varphi }^{k}}(u)\cdot {{\lim }_{u\to \infty }}\frac{\alpha }{c}\int\limits_{0}^{u}{\left[ 1-F(z) \right]dz}=1-\frac{\alpha \mu }{c}+\frac{\alpha \mu }{c}=1, \nonumber
\end{align}
from this it follows by induction that ${{\lim }_{k\to \infty }}{{\varphi }^{k+1}}(u)=1$.
\end{proof}

{\bf Corollary 3.1.}
If ${{\varphi }^{0}}(u)\equiv 1$ , then for all k the approximations ${{\varphi }^{k}}(u)$ are monotone (non-decreasing) continuous functions such that $1-{\alpha \mu }/{c}\;=\varphi (0)\le {{\varphi }^{k}}(u)\le 1$ and ${{\lim }_{k\to \infty }}{{\varphi }^{k}}(u)=1$.

{\bf Lemma 3.3.}
If ${{\varphi }^{k}}(u)\ge 1-{{e}^{-Lu}}$, then ${{\varphi }^{k+1}}(u)\ge 1-{{e}^{-Lu}}$, where $L$ is the Lundberg constant, which is the root of equation (\ref{(2.5)}).
\begin{proof}
Let ${{\varphi }^{k}}(u)\ge 1-{{e}^{-Lu}}$, then
	\[{{\varphi }^{k+1}}(u)\ge 1-\frac{\alpha \mu }{c}+\frac{\alpha }{c}\int\limits_{0}^{u}{\left( 1-{{e}^{-L(u-z)}} \right)\left( 1-F(z) \right)dz}=\]
	\[=1-\frac{\alpha \mu }{c}+\frac{\alpha }{c}\int\limits_{0}^{u}{\left( 1-F(z) \right)dz}-{{e}^{-Lu}}\frac{\alpha }{c}\int\limits_{0}^{u}{{{e}^{Lz}}\left( 1-F(z) \right)dz}=\]
	\[=1-\frac{\alpha \mu }{c}+\frac{\alpha }{c}\int\limits_{0}^{u}{\left( 1-F(z) \right)dz}+\frac{\alpha }{c}\int\limits_{u}^{+\infty }{{{e}^{L(z-u)}}\left( 1-F(z) \right)dz}-{{e}^{-Lu}}\frac{\alpha }{c}\int\limits_{0}^{+\infty }{{{e}^{Lz}}\left( 1-F(z) \right)dz}\ge \]
	\[\ge 1-\frac{\alpha \mu }{c}+\frac{\alpha }{c}\int\limits_{0}^{u}{\left( 1-F(z) \right)dz}+\frac{\alpha }{c}\int\limits_{u}^{+\infty }{\left( 1-F(z) \right)dz}-{{e}^{-Lu}}=\]
	\[=1-{{e}^{-Lu}}-\frac{\alpha \mu }{c}+\frac{\alpha }{c}\int\limits_{0}^{+\infty }{\left( 1-F(z) \right)dz}=1-{{e}^{-Lu}}\ .\ \]
\end{proof}

{\bf Corollary 3.2.}
If ${{\varphi }^{0}}(u)\equiv 1$, then ${{\varphi }^{k}}(u)\ge 1-{{e}^{-Lu}}\ \ \forall \ k\ge 0$.

{\bf Corollary 3.3.}
If ${{\varphi }^{0}}(u)=1-{{e}^{-Lu}}$, then the sequence of functions $\left\{ {{\varphi }^{k}}(u) \right\}$ monotonically increases and converges to the solution $\varphi (u)$ of equation 
(\ref{(3.2)}).

\section{The method of expansion with respect to a small parameter}

For the case of small insurance loadings, when the method of successive approximations may converge slowly, another method for finding $\varphi(u)$ is justified in [41], namely, the method of expansion of $\varphi$ with respect to a small parameter, where the expansion coefficients are found by solving a sequence of some Volterra integral equations that no longer contain a small parameter.

Denote $\varepsilon=1-\alpha\mu/c$, then $0<\varepsilon<1$, and equation (\ref{(3.2)}) takes the form:
\begin{equation}
    \varphi(u)=\varepsilon+\frac{\alpha}{c}\int_0^u{\varphi(u-z)[1-F(z)]dz}.\label{(3.6)}
\end{equation}
As the estimate of Theorem 3.1 shows, the method of successive approximations for solving the integral equation (75) works the better the closer $\varepsilon$ is to 1, and the convergence rate estimate becomes poor when $\varepsilon$ is close to 0. In the latter case, it is advisable to apply the method of expansion of the solution $\varphi$ with respect to the small parameter $\varepsilon$. Represent $\varphi(u)=\sum_{k=1}^\infty \varepsilon^k \varphi_k(u)$. Since $1=1/\varepsilon$, we have $\varphi(u)=\sum_{k=0}^\infty \varepsilon^k \varphi_k(u)$ with $\varphi_0(u)=1$. Substituting this expansion into (\ref{(3.6)}) and equating coefficients of equal powers, we get:
\[{{\varphi }_{1}}(u)=1+\left( {1}/{\mu }\; \right)\int_{\,0}^{\,u}{{{\varphi }_{1}}(u-z)\left( 1-F(z) \right)dz},\]
	\[{{\varphi }_{k}}(u)=\frac{1}{\mu }\int_{\,0}^{\,u}{{{\varphi }_{k}}(u-z)\left( 1-F(z) \right)dz}-\frac{1}{\mu }\int_{\,0}^{\,u}{{{\varphi }_{k-1}}(u-z)\left( 1-F(z) \right)dz},\ \ \ \ \ k=2,3,...\]
Thus, the coefficients of the expansion $\varphi_k(u)$ satisfy Volterra integral equations, which can be solved sequentially by numerical methods.

Let us denote ${{f}_{0}}(u)\equiv 1$, \[{{f}_{k}}(u)=-\left( {1}/{\mu }\; \right)\int_{\,0}^{\,u}{{{\varphi }_{k}}(u-z)\left( 1-F(z) \right)dz}.\] Obviously, 
\[{{\varphi }_{k}}(u)=\left( {1}/{\mu }\; \right)\int_{\,0}^{\,u}{{{\varphi }_{k}}(u-z)\left( 1-F(z) \right)dz}+{{f}_{k-1}}(u),\] 
$\left\| {{f}_{0}} \right\|=1$, $\left\| {{f}_{k-1}} \right\|\le {{\max }_{u\in [0,v]}}{{\varphi }_{k-1}}(u)\cdot \left( {1}/{\mu }\; \right)\int_{\,0}^{\,+\infty }{\left( 1-F(z) \right)dz}=\left\| {{\varphi }_{k-1}} \right\|,\ \ \ k=1,2,...$,
For the function ${{\varphi }_{k}}(u)$, taking into account the fact that $0\le 1-F(u)\le 1$, the following estimate is valid (see [27, Chapter 9, §3]): 
\[\left\| {{\varphi }_{k}} \right\|\le \sum\nolimits_{i=0}^{+\infty }{\frac{{{v}^{i}}}{{{\mu }^{i}}i!}\left\| {{f}_{k-1}} \right\|}={{e}^{{v}/{\mu }\;}}\left\| {{f}_{k-1}} \right\|\le {{e}^{k{v}/{\mu }\;}}.\]
The series \[\varphi (u,q)=\sum\nolimits_{k=1}^{+\infty }{{{\varphi }_{k}}(u){{q}^{k}}}\] is known to be absolutely convergent on the interval $[0,v]$ for $q{{e}^{{v}/{\mu }\;}}<1$ and is a solution of equation (\ref{(3.6)}). Thus, for a given $q$, successive approximations \[{{\varphi }^{k}}(u,q)=\sum\nolimits_{i=1}^{k}{{{\varphi }_{i}}(u){{q}^{i}}}\] certainly converge to the solution of equation (\ref{(3.6)}) on the interval $\left[ 0,\mu \cdot \ln \left( {1}/{q}\; \right) \right]$.

\section{Example}

Let
\begin{equation*}
    F(z)=\varepsilon(1-e^{-\varepsilon z})+(1-\varepsilon)(1-e^{-z}),\quad 0<\varepsilon<1;
		\;\;\;\mu =\int_{0}^{+\infty }{\left( 1-F(z) \right)dz}=2-\varepsilon. 
\end{equation*}
Denote $q=1-\alpha\mu/c$, then $\alpha/c=(1-q)/\mu$, and the probability of non-ruin $\varphi(u)$ depends on the initial capital $u$ and two parameters $q,\mu$. The parameter $q$ is expressed in terms of the insurance loading $\rho=c/(\alpha\mu)-1$ as follows: $q=\rho/(1+\rho)$. Since the real values of the insurance loading may be close to zero, the values of the parameter $q$ may also be close to zero. Let $\varepsilon=q=0.1$.

The iterative process for solving equation (\ref{(3.6)}) can be written in the form
\begin{equation*}
    \varphi^{k+1}(u)=\frac{1-q}{\mu}\int_{u}^{+\infty}{e^{-\frac{1-q}{\mu}(y-u)}\int_0^y{\varphi^k(y-z)dF(z)dy}},\quad k=0,1,\ldots
\end{equation*}

If $\varphi^0(u)\equiv 1$, then according to Lemma 3.1, the sequence of approximations decreases monotonically and converges from above to the solution, and according to Lemma 3.2, $\varphi^k(u)\ge 1-e^{-Lu}$. The course of iterations is shown in Fig. 3.1.

\begin{figure}[h]
    \centering
    \includegraphics[width=0.8\textwidth]{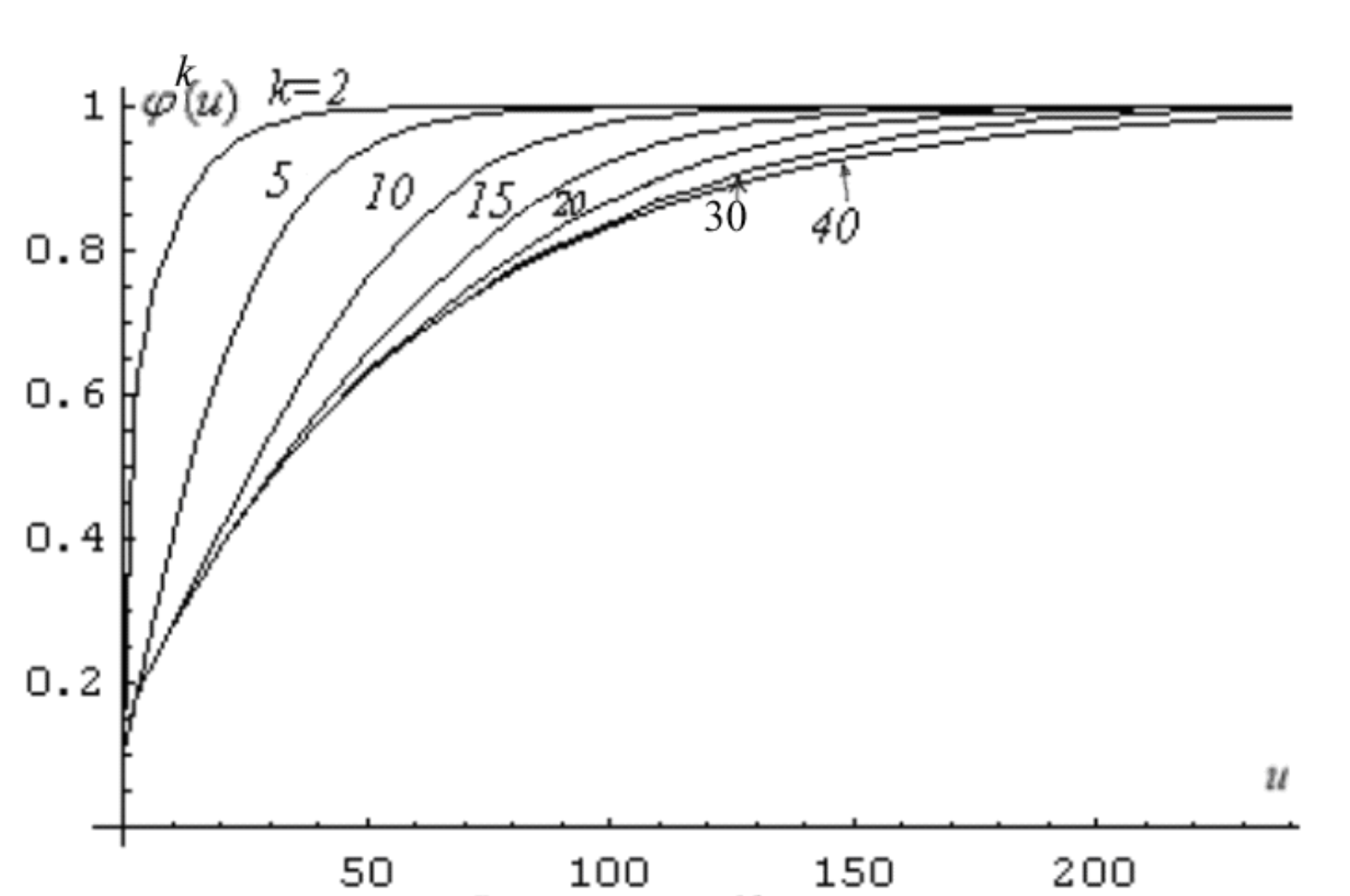}
    \caption{Iterations $\varphi^{0}(u), \ldots, \varphi^{40}(u)$}
    \label{fig:3.1}
\end{figure}

If $\varphi^0(u)=\varphi(0)=1-\alpha\mu/c$, then according to Lemma 3.1, the sequence of approximations increases monotonically; the course of the corresponding iterations is shown in Fig. 3.2.

\begin{figure}[h]
    \centering
    \includegraphics[width=0.8\textwidth]{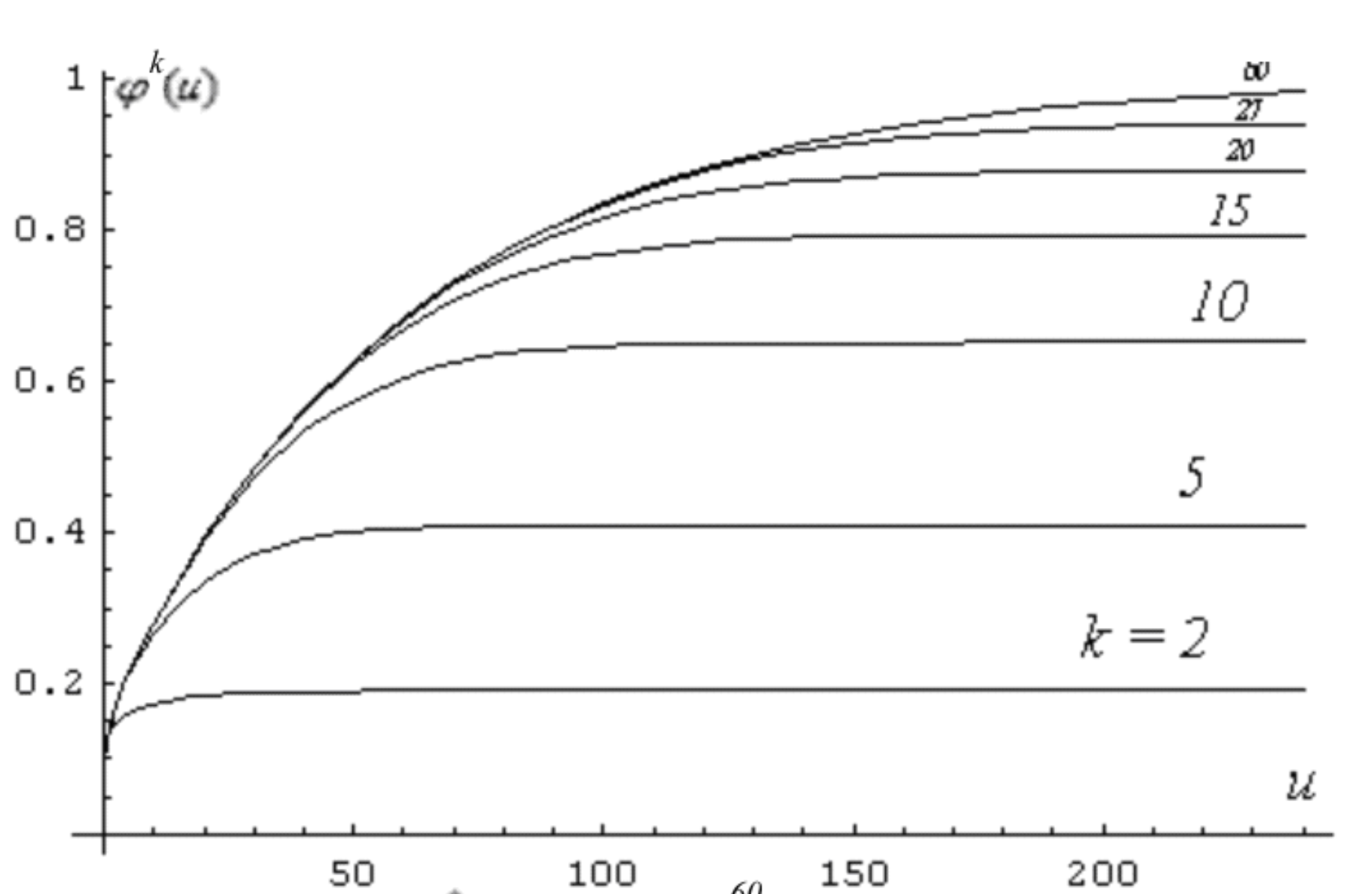}
    \caption{Iterations $\varphi^{0}(u), \ldots, \varphi^{60}(u)$}
    \label{fig:3.2}
\end{figure}

According to Theorem 3.1, the a priori estimate of the accuracy of the approximations has the form:
\begin{equation*}
    \|\varphi-\varphi^k\|\le \frac{p^k}{1-p}.
\end{equation*}
Calculations show that the actual accuracy of the approximations is significantly higher than the indicated a priori accuracy, for example, $\|\varphi-\varphi^{5}\|\le 0.001$.

\section{Conclusions to the section}

It is known that the probabilities of ruin and non-ruin of the classical risk process, as functions of the initial capital $u$, satisfy renewal integral equations of Volterra type (see [35, (3.24), (3.13)]) and boundary conditions at infinity, $\varphi(+\infty)=1$, $\psi(+\infty)=0$. Therefore, in principle, they can be found by the Picard method of successive approximations; it is also necessary to guarantee the fulfillment of the boundary conditions. The method of successive approximations also allows finding analytical approximate solutions for the probability $\varphi$ as a function of $u$ and other parameters of the process. As is known, Picard's method for Volterra integral equations converges on every finite interval of values of $u$ by virtue of the contraction mapping principle. However, in general, contraction occurs only after several iterations; therefore, the guaranteed convergence of the method for large $u$ may be very slow. In this section, the properties of successive approximations for finding the probability of non-ruin $\varphi(u)$ of the classical risk process are investigated in detail. It is shown that for the Volterra equation under consideration, which is satisfied by $\varphi(u)$, contraction occurs at each iteration with coefficient $p=\alpha\mu/c$ uniformly over all $u$, which ensures a high convergence rate for all $u$. It is shown that when starting from the initial functions $\varphi^0(u)\equiv 1$ and $\varphi^0(u)=\varphi(0)=1-\alpha\mu/c$, the successive approximations converge monotonically to the solution $\varphi(u)$, and when starting from the initial function $\varphi^0(u)\equiv 1$, all successive approximations $\varphi^k(u)$ satisfy the condition $\varphi^k(u)\ge 1-e^{-Lu}$, where $L$ is the Lundberg constant, and thus the boundary conditions $\varphi^k(+\infty)=1$ are satisfied at each iteration and in the limit, $\varphi(+\infty)=1$. For the case of small insurance loadings, when the method of successive approximations may converge slowly, another method for finding $\varphi(u)$ is justified, namely, the method of expansion of $\varphi$ with respect to a small parameter, where the expansion coefficients are found by solving a sequence of some Volterra integral equations that no longer contain a small parameter. The theoretical results are illustrated by a numerical example.

\chapter{APPLICATION OF THE METHOD OF SUCCESSIVE APPROXIMATIONS FOR FINDING THE PROBABILITY OF NON-RUIN OF A RISK PROCESS WITH A NONLINEAR PREMIUM GROWTH FUNCTION}

The purpose of this section is to generalize the results presented in Chapter 3 for the numerical solution of general (already non-Volterra) integral equations of the theory of risk processes. The results of this chapter are published in [42].

\section{Risk process with a nonlinear premium growth function}

Consider a risk process of the form
\begin{equation*}
    \xi _t = u + \int_{\,0}^{\,t}{c(\xi_\tau)d\tau} - S_t,\quad t\ge 0,
\end{equation*}
where $t$ is time; $\xi_t$ is the current capital of the company; $u$ is the initial capital of the insurance company; $c(\cdot)$ is a non-negative piecewise continuous function expressing the premium income rate depending on the current capital; $S_t$ is the aggregated payments by time $t$, $S_t=\sum_{k=1}^{N_t} z_k$; $z_k$ are independent identically distributed random variables (claims) with distribution function $F(z)$ and mean $\mu$, $F(z)=0$ for $z\le 0$; $N_t$ is the number of payments by time $t$ (a Poisson process with intensity $\alpha$).

Introduce the function $U(u,t)$ as the solution of the Cauchy problem: $dU/dt=c(U),\ U(0)=u$.

It is easy to see that for the probability of non-ruin $\varphi(u)$ of an insurance company with initial reserve $u$ on an infinite time interval, the following integral equation holds (compare with [35, (3.12)]):
\begin{equation}
    \varphi (u)=A\varphi (u):=\int\limits_{0}^{+\infty }{\alpha e^{-\alpha t}\int\limits_{0}^{U(u,t)}{\varphi (U(u,t)-z)dF(z)dt}}\label{(4.1)}
\end{equation}
with boundary conditions
\begin{equation}
    \varphi (u)=0\ \text{for}\ u<0,\quad \varphi (+\infty)=1,\label{(4.2)}
\end{equation}
and the probability of non-ruin $\varphi(u,T)=\Pr\{\xi_t\ge 0\ \forall t\in[0,T]\}$ on the finite time interval $[0,T]$ is determined by the equation
\begin{equation}
    \varphi (u,T)=\int\limits_{0}^{T}{\alpha e^{-\alpha t}\int\limits_{0}^{U(u,t)}{\varphi (U(u,t)-z,T-t)dF(z)dt}}\label{(4.3)}
\end{equation}
with boundary conditions $\varphi(u)=0$ for $u<0$ and $\varphi(+\infty,T)=1\ \forall T\ge 0$. Equation (\ref{(4.1)}) in the general case does not reduce to the Volterra integral equation (\ref{(2.3)}). Equation (\ref{(4.1)}) may have no solutions satisfying the boundary condition $\varphi(+\infty)=1$. For example, its simplified version (\ref{(2.2)}) with $U(u,t)=u+ct$ has such a solution only under the condition $\alpha\mu/c<1$. Otherwise, equation (\ref{(2.2)}) has only the solution identically equal to zero.

\section{The method of successive approximations}

The considered integral equations (\ref{(2.6)}) - (\ref{(2.16)}), (\ref{(4.1)}) have the general form $\varphi=A\varphi$, where $A$ is the corresponding integral operator, $\varphi$ is a scalar or vector function of one or more arguments, satisfying boundary conditions. The method of successive approximations has the form $\varphi^{k+1}=A\varphi^k$, $\varphi^0(u)$ is some initial approximation, for example, $\varphi^0(u)\equiv 1$. Let us investigate this method using equation ((\ref{(4.1)})) as an example.

Consider the following process of successive approximations:
\begin{equation}\label{(4.4)}
    \begin{aligned}
        \varphi^{k+1}(u) & = \int\limits_{0}^{+\infty }{\alpha e^{-\alpha t}\int\limits_{0}^{U(u,t)}{\varphi^k(U(u,t)-z)dF(z)dt}}; \\
        0 & \le \varphi^0(u)\le 1,\quad u\ge 0;\quad k=0,1,\ldots
    \end{aligned}
\end{equation}

%\begin{assumption}
{\bf Assumption 4.1.}
    Let
    \begin{enumerate}
        \item[(a)] The reserve growth function $U(u,t)$ be monotonically non-decreasing (increasing) in $u$ and $t$, $\lim_{t\to +\infty}U(u,t)=+\infty$.
        \item[(b)] There exist constants $u_*\ge 0$, $c_*>0$ and $L>0$ such that $U(u,t)\ge u+c_* t$ for all $u\ge u_*$ and
        \begin{equation*}
            \left( \alpha/c_* \right)\int_{\,0}^{\,+\infty }{e^{Lz}[1-F(z)]dz}\le 1
        \end{equation*}
        (then $\alpha\mu/c_*<1$).
    \end{enumerate}
    
		Denote by $\{\bar{\varphi}^k(u)\}$ the sequence of approximations (\ref{(4.4)}) starting from the initial approximation $\varphi^0(u)\equiv 1$, and by $\{\underline{\varphi}^k(u)\}$ the sequence of approximations (\ref{(4.4)}) starting from the initial approximation $\varphi^0(u)=\varphi_*(u):=\max\{0,1-e^{-L(u-u_*)}\}$.
%\end{assumption}

%\begin{lemma}
{\bf Lemma 4.1.}
    Let Assumptions 4.1(a) be satisfied. Then
    \begin{enumerate}
        \item[(i)] each function $\bar{\varphi}^k(u)$ is monotonically non-decreasing (increasing) in $u\ge 0$ and $0\le \bar{\varphi}^k(u)\le 1$, $k\ge 0$;
        \item[(ii)] the sequence of values $\{\bar{\varphi}^k(u)\ge 0,\ k=0,1,\ldots\}$ decreases monotonically with increasing $k$;
        \item[(iii)] the limit function $\bar{\varphi}(u)=\lim_{k\to\infty}\bar{\varphi}^k(u)$ is a solution of equation (85), $\bar{\varphi}(u)$ is monotone (non-decreasing) in $u$ and $0\le \bar{\varphi}(u)\le 1$;
        \item[(iv)] under the additional Assumption 4.1(b), for all $k$, $\bar{\varphi}^k(u)\ge \varphi_*(u)=\max\{0,1-e^{-L(u-u_*)}\}$, where $L$ is the Lundberg constant, and thus $1\ge \bar{\varphi}(u)\ge \max\{0,1-e^{-L(u-u_*)}\}$.
    \end{enumerate}
%\end{lemma}

\begin{proof}
(i) The statement is proved by induction. For $k=0$ it is true. If $0\le {{\bar{\varphi }}^{k}}(u)\le 1$, then 
\[0\le {{\bar{\varphi }}^{k+1}}(u)\le \int\limits_{0}^{+\infty }{\alpha {{e}^{-\alpha t}}\int\limits_{0}^{U(u,t)}{dF(z)dt}}=\int\limits_{0}^{+\infty }{\alpha {{e}^{-\alpha t}}\int\limits_{0}^{+\infty }{dF(z)}dt\le }\int\limits_{0}^{+\infty }{\alpha {{e}^{-\alpha t}}dt=1}.\]
If ${{\bar{\varphi }}^{k}}(u)$ is monotone in $u$, then due to the monotonicity of $U(u,t)$ the function ${{\varphi }^{k+1}}(u)$ is also monotone in $u$.

(ii) The statement follows by induction from the relations:
	\[{{\bar{\varphi }}^{1}}(u)=\int\limits_{0}^{+\infty }{\alpha {{e}^{-\alpha t}}\int\limits_{0}^{U(u,t)}{dF(z)}dt\le }\int\limits_{0}^{+\infty }{\alpha {{e}^{-\alpha t}}dt=1}={{\bar{\varphi }}^{0}}(u),\]
	\[{{\bar{\varphi }}^{k+1}}(u)-{{\bar{\varphi }}^{k}}(u)=\int\limits_{0}^{+\infty }{\alpha {{e}^{-\alpha t}}\int\limits_{0}^{U(u,t)}{\left[ {{{\bar{\varphi }}}^{k}}\left( U(u,t)-z \right)-{{{\bar{\varphi }}}^{k-1}}\left( U(u,t)-z \right) \right]dF(z)}dt}.\]
	(iii) For any $u$, the sequence $\left\{ {{{\bar{\varphi }}}^{k}}(u) \right\}$ is monotone decreasing and bounded below, so there exists a limit $\bar{\varphi }(u)={{\lim }_{k\to \infty }}{{\bar{\varphi }}^{k}}(u)$. The limit function $\bar{\varphi }(u)$, like all ${{\bar{\varphi }}^{k}}(u)$, is monotone in $u$ and satisfies the relation $0\le \bar{\varphi }(u)\le 1$. Passing to the limit in equality (\ref{(4.4)}) with respect to $k\to \infty $, by virtue of Lebesgue's theorem on the limit transition under the integral sign, we obtain that the limit function $\bar{\varphi }(u)$ satisfies equation (\ref{(4.1)}).
	
	(iv) Obviously, ${{\bar{\varphi }}^{0}}(u)\equiv 1\ge \max \left\{ 0.1-{{e}^{-L\left( u-{{u}_{*}} \right)}} \right\}$. Let us show that if ${{\bar{\varphi }}^{k}}(u)\ge \max \left\{ 0.1-{{e}^{-L\left( u-{{u}_{*}} \right)}} \right\}$, then ${{\bar{\varphi }}^{k+1}}(u)\ge \max \left\{ 0.1-{{e}^{-L\left( u-{{u}_{*}} \right)}} \right\}$. Using the monotonicity of ${{\bar{\varphi }}^{k}}$ and the fact that ${{\bar{\varphi }}^{k}}(U(u,t)-z)\ge \max \left\{ 0,1-{{e}^{-L\left( U(u,t)-z-{{u}_{*}} \right)}} \right\}$, we obtain
	\begin{align}
  & {{{\bar{\varphi }}}^{k+1}}(u)=\int\limits_{0}^{\infty }{\alpha {{e}^{-\alpha t}}}\int\limits_{0}^{U(u,t)}{{{{\bar{\varphi }}}^{k}}(U(u,t)-z)}dF(z)dy\ge \nonumber \\ 
 & \,\,\,\,\,\,\,\,\,\,\,\,\,\,\,\,\,\ge \int\limits_{0}^{\infty }{\alpha {{e}^{-\alpha t}}}\int\limits_{0}^{U(u,t)}{\max \left\{ 0,1-{{e}^{-L\left( U(u,t)-z-{{u}_{*}} \right)}} \right\}}dF(z)dy= \nonumber\\ 
 & \,\,\,\,\,\,\,\,\,\,\,\,\,\,\,\,\,=\int\limits_{0}^{\infty }{\alpha {{e}^{-\alpha t}}}\int\limits_{0}^{\max \left\{ 0,U(u,t)-{{u}_{*}} \right\}}{\left( 1-{{e}^{-L\left( U(u,t)-z-{{u}_{*}} \right)}} \right)}dF(z)dy. \nonumber 
\end{align}
Let $u\ge {{u}_{*}}$. Integrating the inner integral by parts and using $F(0)=0$, we obtain
	\[{{\bar{\varphi }}^{k+1}}(u)\ge \int\limits_{0}^{\infty }{\alpha {{e}^{-\alpha t}}}{{e}^{-L(U(u,t)-{{u}_{*}})}}L\int\limits_{0}^{U(u,t)-{{u}_{*}}}{{{e}^{Lz}}}F(z)dzdt.\]
Let's transform	
	\begin{align}
  & L\int\limits_{0}^{U(u,t)-{{u}_{*}}}{{{e}^{Lz}}}F(z)dz=L\int\limits_{0}^{U(u,t)-{{u}_{*}}}{{{e}^{Lz}}}\left( 1-\left( 1-F(z) \right) \right)dz=\nonumber\\
	&={{e}^{L(U(u,t)-{{u}_{*}})}}-1-L\int\limits_{0}^{U(u,t)-{{u}_{*}}}{{{e}^{Lz}}}\left( 1-F(z) \right)dz\ge  \nonumber\\ 
 & \ge {{e}^{L(U(u,t)-{{u}_{*}})}}-1-L\int\limits_{0}^{+\infty }{{{e}^{Lz}}}\left( 1-F(z) \right)dz.\nonumber 
\end{align}
Thus, for $u\ge {{u}_{*}}$ it will be
	\begin{align}
  & {{{\bar{\varphi }}}^{k+1}}(u)\ge \int\limits_{0}^{\infty }{\alpha {{e}^{-\alpha t}}}{{e}^{-L(U(u,t)-{{u}_{*}})}}\left( {{e}^{L(U(u,t)-{{u}_{*}})}}-1-L\int\limits_{0}^{+\infty }{{{e}^{Lz}}}\left( 1-F(z) \right)dz \right)dt= \nonumber\\ 
 & \ \ \ \ \ \ \ \ \ \ \ =1-\alpha \left( 1+L\int\limits_{0}^{+\infty }{{{e}^{Lz}}}\left( 1-F(z) \right)dz \right)\cdot \int\limits_{0}^{+\infty }{{{e}^{-\alpha t-L(U(u,t)-{{u}_{*}})}}dt} \nonumber 
\end{align}
and taking into account that \[\int_{0}^{+\infty }{{{e}^{Lz}}\left( 1-F(z) \right)dz}\le {{{c}_{*}}}/{\alpha }\;\] and $U(u,t)\ge u+{{c}_{*}}t$, we obtain ${{\bar{\varphi }}^{k+1}}(u)\ge 1-{{e}^{-L(u-{{u}_{*}})}}$.
The lemma is proved.	
\end{proof}

%\begin{lemma}
{\bf Lemma 4.2.}
    Let Assumptions 4.1(a) and (b) be satisfied. Then
    \begin{enumerate}
        \item[(i)] each function $\underline{\varphi}^k(u)$ is monotonically non-decreasing (increasing) in $u\ge 0$ and $0\le \underline{\varphi}^k(u)\le 1$, $k\ge 0$;
        \item[(ii)] the sequence of values $\{\underline{\varphi}^k(u)\ge 0,\ k=0,1,\ldots\}$ increases monotonically with increasing $k$, and thus for all $k$, $\underline{\varphi}^k(u)\ge \max\{0,1-e^{-L(u-u_*)}\}$, where $L$ is the Lundberg constant;
        \item[(iii)] the limit function $\underline{\varphi}(u)=\lim_{k\to\infty}\underline{\varphi}^k(u)$ is a solution of equation (\ref{(4.1)}), $\underline{\varphi}(u)$ is monotone (non-decreasing) in $u$ and $1\ge \underline{\varphi}(u)\ge \varphi_*(u)=\max\{0,1-e^{-L(u-u_*)}\}$.
    \end{enumerate}
%\end{lemma}

\begin{proof}
The proof of item (i) is the same as in Lemma 4.1. From the proof of item (iv) of Lemma 4.1 it follows that ${{\underset{\scriptscriptstyle-}{\varphi }}^{1}}(u)\ge \max \left\{ 0,1-{{e}^{-L\left( u-{{u}_{*}} \right)}} \right\}={{\underset{\scriptscriptstyle-}{\varphi }}^{0}}(u)$. From this, similarly to the proof of point (ii) of Lemma 4.1, it follows that ${{\underset{\scriptscriptstyle-}{\varphi }}^{k+1}}(u)\ge {{\underset{\scriptscriptstyle-}{\varphi }}^{k}}(u)\ge {{\underset{\scriptscriptstyle-}{\varphi }}^{0}}(u)$ for any $k\ge 0$ and $u\ge 0$. The proof of point (iii) of Lemma 4.2 is similar to the proof of point (iv) of Lemma 4.1.
\end{proof}

%\begin{lemma}
{\bf Corollary 4.1.}
    If equation (\ref{(4.1)}) has a unique monotone solution $\varphi(u)$, then under Assumption 4.1(a), $\bar{\varphi}(u)=\varphi(u)$, and under Assumptions 4.1(a), (b), $\bar{\varphi}(u)=\underline{\varphi}(u)=\varphi(u)$ and the boundary condition $\varphi(+\infty)=1$ is satisfied.
%\end{lemma}

The following lemmas establish the contraction properties of the operator $A$ from (\ref{(4.1)}). For a bounded function $\varphi(u),\ u\ge 0$, define the norm as $\|\varphi\|=\sup_{u\ge 0}|\varphi(u)|$.

%\begin{lemma}
{\bf Lemma 4.3}
    (On the weak contraction property of the operator $A$ when $F(\cdot)<1$). For any monotone non-decreasing bounded functions $\varphi_1(u),\varphi_2(u)$ and for any $v>0$, the following holds:
    \begin{equation*}
        \sup_{u\in[0,v]}|A\varphi_1(u)-A\varphi_2(u)|\le q(v)\cdot \|\varphi_1-\varphi_2\|,
    \end{equation*}
    where $\|\varphi_1-\varphi_2\|=\sup_{u\in[0,+\infty)}|\varphi_1(u)-\varphi_2(u)|$, and
    \begin{equation*}
        q(v)=\int_{0}^{+\infty}{\alpha e^{-\alpha t}F(U(v,t))dt}<1.
    \end{equation*}
%\end{lemma}

\begin{proof}
Let ${{\varphi }_{1}}(u),\,\,{{\varphi }_{2}}(u)$ be monotone (non-decreasing) bounded functions. For $u\in \left[ 0,v \right]$, due to the monotonicity of $U(\cdot ,t)$, the following estimates hold:
\begin{align}
& \left| A{{\varphi }_{1}}(u)-A{{\varphi }_{2}}(u) \right|\le \int\limits_{0}^{+\infty }{\alpha {{e}^{-\alpha t}}\int\limits_{0}^{U(u,t)}{\left| {{\varphi }_{1}}(U(u,t)-z)-{{\varphi }_{2}}(U(u,t)-z) \right|dF(z)}dt}\le \nonumber\\ 
& \ \ \ \ \ \ \ \ \ \ \ \ \ \ \ \ \ \ \le \left\| {{\varphi }_{1}}-{{\varphi }_{2}} \right\|\int\limits_{0}^{+\infty }{\alpha {{e}^{-\alpha t}}\int\limits_{0}^{U(u,t)}{dF(z)}dt}=\left\| {{\varphi }_{1}}-{{\varphi }_{2}} \right\|\int\limits_{0}^{+\infty }{\alpha {{e}^{-\alpha t}}F\left( U(u,t) \right)dt}\le \nonumber\\ 
& \,\,\,\,\,\,\,\,\,\,\,\,\,\,\,\,\,\,\,\,\,\,\,\,\,\,\,\,\,\,\le \left\| {{\varphi }_{1}}-{{\varphi }_{2}} \right\|\int\limits_{0}^{+\infty }{\alpha {{e}^{-\alpha t}}F\left( U(v,t) \right)dt}=q(v)\cdot \left\| {{\varphi }_{1}}-{{\varphi }_{2}} \right\|. \nonumber
\end{align}
Since $F(\cdot )$ and $U(\cdot ,t)$ are monotonic (not decreasing) and $F(\cdot )<1$, then $F\left( U(v,t) \right)<1$ and, therefore, $q(v)<1$. The lemma is proven.
\end{proof}

%\begin{lemma}
{\bf Lemma 4.4}
    (On the contraction property of the operator $A$ when $F(\cdot)<1$). For any monotone non-decreasing bounded functions $\varphi_1(u),\varphi_2(u)$ such that $\|A\varphi_1-A\varphi_2\|\ge \varepsilon>0$ and $0\le \varphi_*(u)\le A\varphi_1(u)\le 1$, $0\le \varphi_*(u)\le A\varphi_2(u)\le 1$, with a monotonically tending to one function $\varphi_*(u)$, the following holds:
    \begin{equation*}
        \|A\varphi_1-A\varphi_2\|\le q(\varepsilon)\cdot \|\varphi_1-\varphi_2\|,\quad q(\varepsilon)<1.
    \end{equation*}
%\end{lemma}

\begin{proof}
Let $\left\| A{{\varphi }_{1}}-A{{\varphi }_{2}} \right\|={{\lim }_{s\to \infty }}\left| A{{\varphi }_{1}}({{u}^{s}})-A{{\varphi }_{2}}({{u}^{s}}) \right|\ge \varepsilon >0$ with some sequence $\left\{ {{u}^{s}} \right\}$. Obviously $\left| A{{\varphi }_{1}}(u)-A{{\varphi }_{2}}(u) \right|\le 1-{{\varphi }_{*}}(u)$. We define $${{u}^{*}}(\varepsilon )={{\sup }_{v\ge 0}}\left\{ v:\,\,1-{{\varphi }_{*}}(v)\ge {\varepsilon }/{2}\; \right\}<+\infty. $$ Then for sufficiently large $s$ we have $1-{{\varphi }_{*}}({{u}^{s}})\ge \left| A{{\varphi }_{1}}({{u}^{s}})-A{{\varphi }_{2}}({{u}^{s}}) \right|\ge {\varepsilon }/{2}\;$ and, consequently, ${{u}^{s}}\le {{u}^{*}}(\varepsilon )$. The assessments hold
\begin{align} 
& \left| A{{\varphi }_{1}}({{u}^{s}})-A{{\varphi }_{2}}({{u}^{s}}) \right|\le \int\limits_{0}^{+\infty }{\alpha {{e}^{-\alpha t}}\int\limits_{0}^{U({{u}^{s}},t)}{\left| {{\varphi }_{1}}(U({{u}^{s}},t)-z)-{{\varphi }_{2}}(U({{u}^{s}},t)-z) \right|dF(z)}dt}\le \nonumber\\ 
& \ \ \ \ \ \ \ \ \ \ \ \ \ \ \ \ \ \ \le \left\| {{\varphi }_{1}}-{{\varphi }_{2}} \right\|\int\limits_{0}^{+\infty }{\alpha {{e}^{-\alpha t}}\int\limits_{0}^{U({{u}^{s}},t)}{dF(z)}dt}=\left\| {{\varphi }_{1}}-{{\varphi }_{2}} \right\|\int\limits_{0}^{+\infty }{\alpha {{e}^{-\alpha t}}F\left( U({{u}^{s}},t) \right)dt}\le \nonumber\\ 
& \,\,\,\,\,\,\,\,\,\,\,\,\,\,\,\,\,\,\,\,\,\,\,\,\,\,\,\,\,\,\le \left\| {{\varphi }_{1}}-{{\varphi }_{2}} \right\|\int\limits_{0}^{+\infty }{\alpha {{e}^{-\alpha t}}F\left( U({{u}^{*}}(\varepsilon ),t) \right)dt}=q(\varepsilon )\cdot \left\| {{\varphi }_{1}}-{{\varphi }_{2}} \right\|,\nonumber
\end{align}
where \[q(\varepsilon )=\int\limits_{0}^{+\infty }{\alpha {{e}^{-\alpha t}}F\left( U({{u}^{*}}(\varepsilon ),t) \right)dt}<1.\] The lemma is proven.
\end{proof}

%\begin{lemma}
{\bf Corollary 4.2.}
    Let Assumption 4.1 be satisfied, set $\varphi_*(u)=\max\{0,1-e^{-L(u-u_*)}\}$. Then for any functions $\varphi_1(u)\ne \varphi_2(u)$ such that $\varphi_*(u)\le \varphi_1(u)\le 1$, $\varphi_*(u)\le \varphi_2(u)\le 1$, the following holds:
    \begin{equation*}
        \|A\varphi_1-A\varphi_2\|<\|\varphi_1-\varphi_2\|.
    \end{equation*}
%\end{lemma}

\begin{proof}
If $\left\| A{{\varphi }_{1}}-A{{\varphi }_{2}} \right\|=0$, the statement of the corollary is obvious. Let $\left\| A{{\varphi }_{1}}-A{{\varphi }_{2}} \right\|>0$. From the proof of Lemma 4.1 it follows that
${{\varphi }_{*}}(u)\le A{{\varphi }_{*}}(u)\le A{{\varphi }_{1}}(u)\le 1$, ${{\varphi }_{*}}(u)\le A{{\varphi }_{*}}(u)\le A{{\varphi }_{2}}(u)\le 1$.
Then, when considering the case of the statement of the consequences, the statement follows from Lemma 4.4.
\end{proof}

%\begin{lemma}
{\bf Lemma 4.5}
    (On the contraction property of the operator $A$ when $F(\cdot)\le 1$). Let $F(z)=1$ for $z\ge \bar{z}$. For any monotone non-decreasing bounded functions $\varphi_1(u),\varphi_2(u)$ such that $\|\varphi_1-\varphi_2\|\ge \varepsilon>0$ and $0\le \varphi_*(u)\le \varphi_1(u)\le 1,\ 0\le \varphi_*(u)\le \varphi_2(u)\le 1$ with a monotonically tending to one function $\varphi_*(u)$, the following holds:
    \begin{equation*}
        \|A\varphi_1-A\varphi_2\|\le q(\varepsilon)\cdot \|\varphi_1-\varphi_2\|,
    \end{equation*}
    where $q(\varepsilon)<1$, $\|\varphi_1-\varphi_2\|=\sup_{u\in[0,+\infty)}|\varphi_1(u)-\varphi_2(u)|$, $\|A\varphi_1-A\varphi_2\|=\sup_{u\in[0,+\infty)}|A\varphi_1(u)-A\varphi_2(u)|$.
%\end{lemma}

\begin{proof}
Obviously, $\left| {{\varphi }_{1}}(u)-{{\varphi }_{2}}(u) \right|\le 1-{{\varphi }_{*}}(u)$. Then for any $u\ge 0$ the following estimates are valid:
\begin{align}
  & \left| A{{\varphi }_{1}}(u)-A{{\varphi }_{2}}(u) \right|\le \int\limits_{0}^{+\infty }{\alpha {{e}^{-\alpha t}}\int\limits_{0}^{U(u,t)}{\left| {{\varphi }_{1}}(U(u,t)-z)-{{\varphi }_{2}}(U(u,t)-z) \right|dF(z)}dt}\le \nonumber \\ 
 & \,\,\,\,\,\,\,\,\,\,\,\,\,\,\,\,\,\,\,\le \int\limits_{0}^{+\infty }{\alpha {{e}^{-\alpha t}}\int\limits_{0}^{U(u,t)}{\min \left\{ \left\| {{\varphi }_{1}}-{{\varphi }_{2}} \right\|,1-{{\varphi }_{*}}(U(u,t)-z) \right\}dF(z)}dt}\le \nonumber \\ 
 & \,\,\,\,\,\,\,\,\,\,\,\,\,\,\,\,\,\,\,\le \int\limits_{0}^{+\infty }{\alpha {{e}^{-\alpha t}}\int\limits_{0}^{\min \left\{ U(u,t),\bar{z} \right\}}{\min \left\{ \left\| {{\varphi }_{1}}-{{\varphi }_{2}} \right\|,1-{{\varphi }_{*}}(U(u,t)-z) \right\}dF(z)}dt}\le \nonumber \\ 
 & \,\,\,\,\,\,\,\,\,\,\,\,\,\,\,\,\,\,\,\le \int\limits_{0}^{+\infty }{\alpha {{e}^{-\alpha t}}\int\limits_{0}^{\min \left\{ U(u,t),\bar{z} \right\}}{\min \left\{ \left\| {{\varphi }_{1}}-{{\varphi }_{2}} \right\|,1-{{\varphi }_{*}}(\max \left\{ 0,U(u,t)-\bar{z} \right\}) \right\}dF(z)}dt}\le \nonumber 
\end{align}
\begin{align}
  & \,\,\,\,\,\,\,\,\,\,\,\,\,\,\,\,\,\,\,\le \int\limits_{0}^{+\infty }{\alpha {{e}^{-\alpha t}}\min \left\{ \left\| {{\varphi }_{1}}-{{\varphi }_{2}} \right\|,1-{{\varphi }_{*}}(\max \left\{ 0,U(u,t)-\bar{z} \right\}) \right\}F\left( \min \left\{ U(u,t),\bar{z} \right\} \right)dt}\le \nonumber \\ 
 & \,\,\,\,\,\,\,\,\,\,\,\,\,\,\,\,\,\,\,\,\le \int\limits_{0}^{+\infty }{\alpha {{e}^{-\alpha t}}\min \left\{ \left\| {{\varphi }_{1}}-{{\varphi }_{2}} \right\|,1-{{\varphi }_{*}}(\max \left\{ 0,U(0,t)-\bar{z} \right\}) \right\}dt}\le \nonumber 
\end{align}
\begin{align}
  & \,\,\,\,\,\,\,\,\,\,\,\,\,\,\,\,\,\,\le \left\| {{\varphi }_{1}}-{{\varphi }_{2}} \right\|.\int\limits_{0}^{+\infty }{\alpha {{e}^{-\alpha t}}\min \left\{ 1,{\left( 1-{{\varphi }_{*}}(\max \left\{ 0,U(0,t)-\bar{z} \right\}) \right)}/{\left\| {{\varphi }_{1}}-{{\varphi }_{2}} \right\|}\; \right\}dt}\le \nonumber  \\ 
 & \,\,\,\,\,\,\,\,\,\,\,\,\,\,\,\,\,\,\le \left\| {{\varphi }_{1}}-{{\varphi }_{2}} \right\|.\int\limits_{0}^{+\infty }{\alpha {{e}^{-\alpha t}}\min \left\{ 1,{\left( 1-{{\varphi }_{*}}(\max \left\{ 0,U(0,t)-\bar{z} \right\}) \right)}/{\varepsilon }\; \right\}dt}. \nonumber 
\end{align}
But since $U(0,t)\to +\infty $, then $\left( 1-{{\varphi }_{*}}(\max \left\{ 0,U(0,t)-\bar{z} \right\}) \right)\to 0$ as $t\to +\infty $ and $\min \left\{ 1,{\left( 1-{{\varphi }_{*}}(\max \left\{ 0,U(0,t)-\bar{z} \right\}) \right)}/{\varepsilon }\; \right\}<1$ for all sufficiently large $t$. It follows that
\[q(\varepsilon )=\int\limits_{0}^{+\infty }{\alpha {{e}^{-\alpha t}}\min \left\{ 1,{\left( 1-{{\varphi }_{*}}(\max \left\{ 0,U(0,t)-\bar{z} \right\}) \right)}/{\varepsilon }\; \right\}dt}<1.\]
The lemma is proved.
\end{proof}

%\begin{corollary}
{\bf Coroolary 4.3.}
    From Lemmas 4.3 and 4.5, it follows that the operator $A$ defined in (\ref{(4.1)}) is a contraction on the set of functions $\varphi(u)$ such that $\varphi_*(u)=\max\{0,1-e^{-L(u-u_*)}\}\le \varphi(u)\le 1$, i.e., for any functions $\varphi_1(u)\ne \varphi_2(u)$ such that $\varphi_*(u)\le \varphi_1(u)\le 1,\ \varphi_*(u)\le \varphi_2(u)\le 1$, the following holds:
    \begin{equation*}
        \|A\varphi_1-A\varphi_2\|<\|\varphi_1-\varphi_2\|.
    \end{equation*}
%\end{corollary}

%\begin{theorem}
{\bf Theorem 4.1}
    (On the convergence of the sequence of approximations). Let Assumptions 4.1 be satisfied. Then:
    \begin{enumerate}
        \item[1)] equation (\ref{(4.1)}) has a monotone solution $\varphi(u)$ such that $1\ge \varphi(u)\ge \varphi_*(u)=\max\{0,1-e^{-L(u-u_*)}\}$;
        \item[2)] $\varphi(u)$ is the unique solution of problems (\ref{(4.1)}), (\ref{(4.2)});
        \item[3)] for any initial approximation $\varphi^0(u)$ such that $1\ge \varphi^0(u)\ge \varphi_*(u)=\max\{0,1-e^{-L(u-u_*)}\}$, the corresponding sequence of approximations (\ref{(4.4)}) converges pointwise to the solution $\varphi(u)$.
    \end{enumerate}
%\end{theorem}

\begin{proof}
Statement 1) follows from Lemmas 4.1, 4.2.

Let us prove assertion 2) of Theorem 4.1 on the uniqueness of the solution to problem (\ref{(4.1)}), (\ref{(4.2)}). Assume the contrary: let there be two solutions ${{\varphi }_{1}}\ne {{\varphi }_{2}}$ (i.e. ${{\varphi }_{1}}(u)\ne {{\varphi }_{2}}(u)$ for some $u\ge 0)$ of problem (\ref{(4.1)}), (\ref{(4.2)}). Put ${{\varphi }^{0}}(u)=\max \left\{ {{\varphi }_{1}}(u),{{\varphi }_{2}}(u) \right\}$. Obviously, ${{\varphi }^{0}}\ne {{\varphi }_{1}}$, ${{\varphi }^{0}}\ne {{\varphi }_{2}}$. Consider the sequence $\left\{ {{\varphi }^{k}}(u):=A{{\varphi }^{k}}(u),\,\,k=0,1,... \right\}$. Obviously, $1\ge {{\varphi }^{0}}(u)\ge {{\varphi }_{1}}(u)$ and $1\ge {{\varphi }^{0}}(u)\ge {{\varphi }_{2}}(u)$. Since the operator $A$ maps a large function to a large one (is monotone), then, taking into account Lemma 4.1, $1\ge {{\varphi }^{1}}(u)=A{{\varphi }^{0}}(u)\ge A{{\varphi }_{1}}(u)={{\varphi }_{1}}(u)$, $1\ge {{\varphi }^{1}}(u)=A{{\varphi }^{0}}(u)\ge A{{\varphi }_{2}}(u)={{\varphi }_{2}}(u)$ and, thus, ${{\varphi }^{1}}(u)\ge \max \left\{ {{\varphi }_{1}}(u),{{\varphi }_{2}}(u) \right\}={{\varphi }^{0}}(u)$. From this it follows by induction that ${{\varphi }^{k+1}}(u)\ge {{\varphi }^{k}}(u)\ge {{\varphi }^{0}}(u)$ for all $k\ge 0$. The sequence $\left\{ {{\varphi }^{k}}(\cdot ) \right\}$ is monotonically increasing and is bounded above by unity, so it has a pointwise limit $\varphi (u)$, $1\ge \varphi (u)\ge {{\varphi }^{0}}(u)$, which is a solution of equation (\ref{(4.1)}) by virtue of Lebesgue's theorem on the limit passage under the integral sign. Obviously, the obtained function $\varphi (u)$ is a solution of problem (\ref{(4.1)}), (\ref{(4.2)}), different from ${{\varphi }_{1}}$ and ${{\varphi }_{2}}$. Thus, two different functions ${{\varphi }_{1}}(u)$ and $\varphi (u)\ge {{\varphi }_{1}}(u)$ are solutions of equation (\ref{(4.1)}) and satisfy the conditions of Lemmas 4.4, 4.5 with ${{\varphi }_{*}}\left( u \right)={{\varphi }_{1}}(u)$. But this contradicts the contraction property of the operator $A$ from Lemma 4.4 or 4.5.

Statement 3 of Theorem 4.1 is proved similarly. Let $1\ge {{\varphi }^{0}}(u)\ge \max \left\{ 0,1-{{e}^{-L(u-{{u}_{*}})}} \right\}$, i.e. ${{\bar{\varphi }}^{0}}(u)\ge {{\varphi }^{0}}(u)\ge {{\underset{\scriptscriptstyle-}{\varphi }}^{0}}(u)$. Then ${{\bar{\varphi }}^{1}}(u)\ge {{\varphi }^{1}}(u)\ge {{\underset{\scriptscriptstyle-}{\varphi }}^{1}}(u)$ and, moreover, for any $k$ ${{\bar{\varphi }}^{k}}(u)\ge {{\varphi }^{k}}(u)\ge {{\underset{\scriptscriptstyle-}{\varphi }}^{k}}(u)$. It follows that there exists a limit $\varphi (u)={{\lim }_{k\to \infty }}{{\varphi }^{k}}(u)\equiv \bar{\varphi }(u)\equiv \underset{\scriptscriptstyle-}{\varphi }(u)$. The theorem is proved.
\end{proof}

%\begin{theorem}
{\bf Theorem 4.2}
    (On uniform convergence and the rate of convergence of the method of successive approximations). Under Assumptions 4.1, the method of successive approximations (\ref{(4.4)}), starting from an initial approximation $\varphi^0(u)$ such that $1\ge \varphi^0(u)\ge \max\{0,1-e^{-L(u-u_*)}\}$, converges uniformly monotonically to the solution $\varphi(u)$ of problems (\ref{(4.4)}), i.e., $\|\varphi^k-\varphi\|=\sup_{u\ge 0}|\varphi^k(u)-\varphi(u)|$ converges monotonically to zero, and moreover, the method (\ref{(4.4)}) converges uniformly to any $\varepsilon$-neighborhood of the solution of problem (\ref{(4.1)}) at the rate of a geometric progression with denominator $q(\varepsilon)$ depending on $\varepsilon$, i.e.,
    \begin{equation*}
        \|\varphi^k-\varphi\|\le (q(\varepsilon))^k\cdot \|\varphi^0-\varphi\|,\quad 0\le q(\varepsilon)<1,
    \end{equation*}
    for all $k$ such that $\|\varphi^k-\varphi\|\ge \varepsilon$.
%\end{theorem}

\begin{proof}
By Theorem 4.1, the solution $\varphi (u)$ of problem (\ref{(4.1)}) exists and is unique. It follows from Lemmas 4.4 and 4.5 that the sequence $\left\| {{\varphi }^{k}}-\varphi \right\|$ is monotonically decreasing and for any $\varepsilon >0$ there exists $q(\varepsilon )$, $0\le q(\varepsilon )<1$, such that for all $k$ such that $\left\| {{\varphi }^{k}}-\varphi \right\|\ge \varepsilon $, $\left\| {{\varphi }^{k}}-\varphi \right\|\le q(\varepsilon )\left\| {{\varphi }^{k-1}}-\varphi \right\|\le {{\left( q(\varepsilon ) \right)}^{k}}\cdot \left\| {{\varphi }^{0}}-\varphi \right\|$. It remains to show that ${{\lim }_{k\to \infty }}\left\| {{\varphi }^{k}}-\varphi \right\|=0$. Suppose on the contrary, that ${{\lim }_{k\to \infty }}\left\| {{\varphi }^{k}}-\varphi \right\|=\varepsilon >0$. Then, by the above, for all $k$ we have
$\varepsilon \le \left\| {{\varphi }^{k}}-\varphi \right\|=\left\| A{{\varphi }^{k-1}}-A\varphi \right\|\le q(\varepsilon )\left\| {{\varphi }^{k-1}}-\varphi \right\|\le {{\left( q(\varepsilon ) \right)}^{k}}\cdot \left\| {{\varphi }^{0}}-\varphi \right\|$,
which is impossible for sufficiently large $k$. The theorem is proved.
\end{proof}

%\begin{remark}
{\bf Remark 4.1.}
    The method (\ref{(4.4)}) for solving equation (\ref{(4.1)}) and the results of Lemmas 4.1, 4.2 and Theorem 4.1 on its convergence fully extend to equation (\ref{(4.3)}) in view of the obvious relation between their solutions: $\varphi(u)\le \varphi(u,T)$ (the probability of non-ruin on a finite time interval is greater than the probability of non-ruin on an infinite time interval).
%\end{remark}

\section{Example}

Let $U(u,t)=u+ct$, $x(z)=z$, and the claim distribution $F(z)$ be a mixture of two exponential distributions:
\begin{equation*}
    F(z)=\varepsilon(1-\exp(-\varepsilon z))+(1-\varepsilon)(1-\exp(-z)),\quad \mu=\int_{\,0}^{\,+\infty}{(1-F(z))dz}=2-\varepsilon.
\end{equation*}
Denote $q=1-\alpha\mu/c$, then $\alpha/c=(1-q)/\mu$, and thus the probability of non-ruin $\varphi(u)$ depends on the initial capital $u$ and two parameters $q,\mu$. The parameter $q$ is expressed in terms of the insurance loading $\rho=c/(\alpha\mu)-1$ as follows: $q=\rho/(1+\rho)$. Since the real values of the insurance loading may be close to zero, the values of the parameter $q$ may also be close to zero. Take $\varepsilon=q=0.1$.

The iterative process (\ref{(4.4)}) for solving equation (\ref{(4.1)}) can be rewritten in the form:
\begin{equation}
    \varphi^{k+1}(u)=\frac{1-q}{\mu}\int_{u}^{+\infty}{e^{-\frac{1-q}{\mu}(y-u)}\int_0^y{\varphi^k(y-z)dF(z)dy}},\quad k=0,1,\ldots\label{(4.5)}
\end{equation}

The functional nature of the iterations, the presence of double integrals and infinite integration limits, and the potentially large number of iterations require an economical numerical implementation of process (\ref{(4.5)}) and consideration of roundoff errors. The integrals in (\ref{(4.5)}) can be calculated using parallel processors.

If $\varphi^0(u)\equiv 1$, then according to Lemma 4.1(ii), the sequence of approximations decreases monotonically and converges from above to the solution of equation (\ref{(4.1)}), and according to Lemma 4.1(iv), $\lim_{u\to+\infty}\varphi^k(u)=1$. The course of iterations is shown in Fig. 4.1.
If $\varphi^0(u)=1-e^{-Lu}$, where the Lundberg constant $L\approx 0.0164$, then according to Lemma 4.2, the sequence of approximations increases monotonically and converges from below to the solution of equation (\ref{(4.1)}); the course of the corresponding iterations is shown in Fig. 4.2. From the graphs, it follows that the difference between the obtained approximations from above and below does not exceed 0.035.

\begin{figure}[h]
    \centering
    \includegraphics[width=0.8\textwidth]{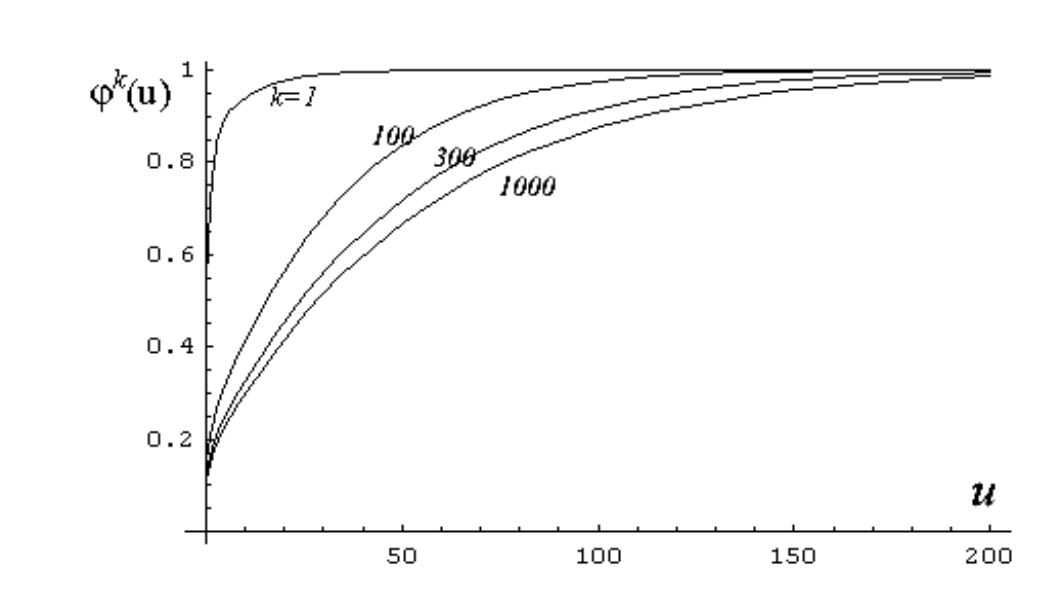}
    \caption{}
    \label{fig:4.1}
\end{figure}

\begin{figure}[h]
    \centering
    \includegraphics[width=0.8\textwidth]{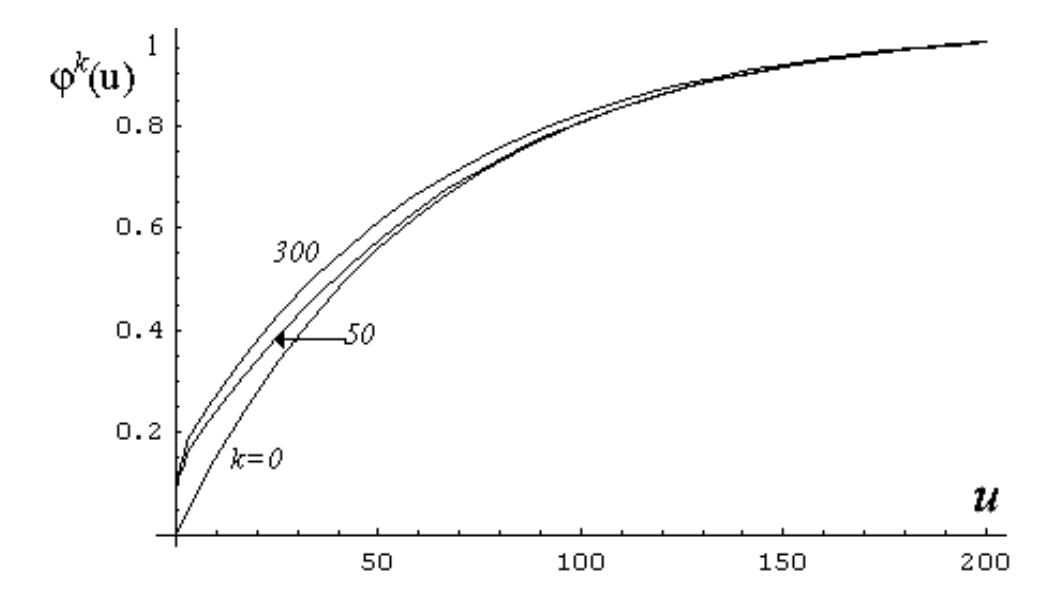}
    \caption{}
    \label{fig:4.2}
\end{figure}

\chapter{APPLICATION OF THE METHOD OF SUCCESSIVE APPROXIMATIONS FOR CALCULATING THE PROBABILITY OF RUIN OF A RISK PROCESS WITH A NON-POISSON CLAIM FLOW}

In this section, the method of successive approximations for finding the probability of ruin of a renewal risk process is studied. This process and the integral equations for the probability of non-ruin were described in Section 2.5. The results of this chapter are published in [44, 106, 107, 108].

\section{The method of successive approximations for finding the probability of ruin}

Consider a renewal risk process with nonlinear capital growth:
\begin{equation}
    \xi _t = u + \int_{\,0}^{\,t}{c(\xi_s)ds} - S_t,\quad t\ge 0,\label{(5.1)}
\end{equation}
where the premium arrival rate $c(\cdot)\ge 0$ depends on the current capital of the company and is a piecewise continuous function; $S_t$ is the aggregated payments by time $t$, $S_t=\sum_{k=1}^{N_t} z_k$; $z_k$ are independent identically distributed random variables (claims) with distribution function $F(z)$ and mean $\mu$, $F(z)=0$ for $z\le 0$; $N_t$ is the number of payments by time $t$. We assume that the time intervals between successive claims are random, independent, and have the same distribution function $K(\cdot)$.

We assume that the reserves (capital) of the insurance company in the absence of claims increase over time $t$ from the initial level $u$ to the value $U(u,t)$, $U(u,0)=u$, according to the equation $dU/dt=c(U)$. For example, in the case $c(\cdot)=const$, we have $U(u,t)=u+ct$.

Then equation (\ref{(2.1)}) takes the form
\begin{equation}
    \varphi (u)=A\varphi (u):=\int_{\,0}^{\,+\infty }{\int_{\,0}^{\,U(u,t)}{\varphi (U(u,t)-z)dF(z)}\,dK(t)}.\label{(5.2)}
\end{equation}
We seek a solution satisfying the boundary condition
\begin{equation}
    \varphi (+\infty)=1.\label{(5.3)}
\end{equation}

%\begin{assumption}
{\bf Assumption 5.1.}
    Let
    \begin{enumerate}
        \item[(a)] the reserve growth function $U(u,t)$ be monotonically non-decreasing (increasing) in $u$ and $t$, $\lim_{t\to +\infty}U(u,t)=+\infty$;
        \item[(b)] there exist constants $u_*\ge 0$, $c_*>0$ and $L>0$ such that $U(u,t)\ge u+c_* t$ for all $u\ge u_*$, and
        \begin{equation}
            \int_{\,0}^{\,+\infty }{e^{Lz}dF}(z)\,\int_{\,0}^{\,+\infty }{e^{-c_*Lt}dK(t)\le 1};
        \label{(5.4)}
				\end{equation}
        \item[(c)] $K(t)\cdot F(z)<1\ \forall t,z\ge 0$.
    \end{enumerate}
%\end{assumption}

Consider the following process of successive approximations for solving equation (\ref{(5.2)}):
\begin{equation}
    \begin{aligned}
        \varphi^{k+1}(u) & = \int_{\,0}^{\,+\infty }{\int_{\,0}^{\,U(u,t)}{\varphi^k(U(u,t)-z)dF(z)}\,dK(t)}, \\
        0 & \le \varphi^0(u)\le 1,\quad u\ge 0,\quad k=0,1,\ldots
    \end{aligned}
		\label{(5.5)}
\end{equation}
Denote by $\{\bar{\varphi}^k(u)\}$ the sequence of approximations (\ref{(5.5)}) starting from the initial approximation $\varphi^0(u)=\varphi^*(u)\equiv 1$, and by $\{\underline{\varphi}^k(u)\}$ the sequence of approximations (\ref{(5.5)}) starting from the initial approximation $\varphi^0(u)=\varphi_*(u)=\max\{0,1-e^{-L(u-u_*)}\}$.

%\begin{lemma}
{\bf Lemma 5.1.}
    Let Assumption 5.1(a) be satisfied. Then
    \begin{enumerate}
        \item[(i)] each function $\bar{\varphi}^k(u)$ is monotonically non-decreasing (increasing) in $u\ge 0$ and $0\le \bar{\varphi}^k(u)\le 1,\quad k=0,1,\ldots$;
        \item[(ii)] the sequence $\{\bar{\varphi}^k(u)\ge 0,\ k=0,1,\ldots\}$ decreases monotonically with increasing $k$;
        \item[(iii)] the limit function $\bar{\varphi}(u)=\lim_{k\to\infty}\bar{\varphi}^k(u)$ is a solution of equation (98), $\bar{\varphi}(u)$ is monotone (non-decreasing) in $u$ and $0\le \bar{\varphi}(u)\le 1$;
        \item[(iv)] under the additional Assumption 5.1(b), for all $k$, $\bar{\varphi}^k(u)\ge \max\{0,1-e^{-L(u-u_*)}\}$, and thus $1\ge \bar{\varphi}(u)\ge \max\{0,1-e^{-L(u-u_*)}\}$, where $\varphi_*(u)=\max\{0,1-e^{-L(u-u_*)}\}$ is the generalized Cramér-Lundberg bound for the probability of non-ruin.
    \end{enumerate}
%\end{lemma}

\begin{proof}
(i) The assertion is proved by induction. For $k=0$ it is obviously true. If $0\le {{\bar{\varphi }}^{k}}(u)\le 1$, then
\[0\le {{\bar{\varphi }}^{\kappa +1}}(u)\le \int_{0}^{+\infty }{dK(t)}\int_{0}^{U(u,t)}{dF(z)}=\int_{0}^{+\infty }{dK(t)}\int_{0}^{+\infty }{dF(z)}\le \int_{0}^{+\infty }{dK(t)}=1.\] 
If ${{\bar{\varphi }}^{k}}(u)$ is monotone in u, then due to the monotonicity of $U(u,t)$, the function ${{\bar{\varphi }}^{k+1}}(u)$ is also monotone in $u$.

(ii)The second statement of the lemma follows by induction from the relations:
	\[{{\bar{\varphi }}^{1}}(u)=\int_{0}^{+\infty }{dK(t)}\int_{0}^{U(u,t)}{dF(z)}\le \int_{0}^{+\infty }{dK(t)}=1={{\bar{\varphi }}^{0}}(u),\]
i.e., ${{\bar{\varphi }}^{1}}(u)\le {{\bar{\varphi }}^{0}}(u)$ for all $u\ge 0$. Since	
\[{{\bar{\varphi }}^{k+1}}(u)-{{\bar{\varphi }}^{k}}(u)=\int_{0}^{+\infty }{dK(t)}\int_{0}^{U(u,t)}{[{{{\bar{\varphi }}}^{k}}(U(u,t)-z)-{{{\bar{\varphi }}}^{k-1}}(U(u,t)-z)]dF(z)},\] 
then the relation ${{\bar{\varphi }}^{k+1}}(u)-{{\bar{\varphi }}^{k}}(u)\le 0$, $u\ge 0$, is satisfied for all $k\ge 0$, i.e., the sequence $\{{{\bar{\varphi }}^{k}}(u)\}$ decreases monotonically with increasing $k$.

(iii) For any u, the sequence $\{{{\bar{\varphi }}^{k}}(u)\}$ is monotone decreasing and bounded below, so there exists a limit $\bar{\varphi }(u)={{\lim }_{k\to \infty }}{{\bar{\varphi }}^{k}}(u)$. The limit function $\bar{\varphi }(u)$, like all ${{\bar{\varphi }}^{k}}(u)$, is monotone in $u$ and satisfies the relation $0\le \bar{\varphi }(u)\le 1$. Passing to the limit in equality (\ref{(5.5)}) with respect to $k\to \infty $ by virtue of Lebesgue's theorem on the limit transition under the integral sign, we obtain that the limit function $\bar{\varphi }(u)$ satisfies equation (\ref{(5.2)}).

(iv) Obviously, ${{\bar{\varphi }}^{0}}(u)\equiv 1\ge \max \left\{ 0.1-{{e}^{-L(u-{{u}_{*}})}} \right\}$. Let us show that if
${{\bar{\varphi }}^{k}}(u)\ge \max \left\{ 0.1-{{e}^{-L(u-{{u}_{*}})}} \right\}$, then ${{\bar{\varphi }}^{k+1}}(u)\ge \max \left\{ 0.1-{{e}^{-L(u-{{u}_{*}})}} \right\}$. Using the monotonicity of ${{\bar{\varphi }}^{k}}(u)$ and the fact that
\[\,\,\,{{\bar{\varphi }}^{k}}(U(u,t)-z)\ge \max \left\{ 0.1-{{e}^{-L(U(u,t)-{{u}_{*}}-z)}} \right\},\]
we obtain:
\begin{align}
  & {{{\bar{\varphi }}}^{k+1}}(u)=\int_{0}^{+\infty }{dK(t)}\int_{0}^{U(u,t)}{{{{\bar{\varphi }}}^{k}}(U(u,t)-z)dF(z)}\ge  \nonumber\\ 
 & \,\,\,\,\,\,\,\,\,\,\,\,\,\,\,\,\,\ge \int_{0}^{+\infty }{dK(t)}\int_{0}^{U(u,t)}{\max \left\{ 0,1-{{e}^{-L(U(u,t)-{{u}_{*}}-z)}} \right\}dF(z)}=\, \nonumber\\ 
 & \,\,\,\,\,\,\,\,\,\,\,\,\,\,\,\,\,=\int_{0}^{+\infty }{dK(t)}\int_{0}^{\max \left\{ 0,U(u,t)-{{u}_{*}} \right\}}{\left( 1-{{e}^{-L(U(u,t)-{{u}_{*}}-z)}} \right)dF(z)}. \nonumber 
\end{align}	
Let $u\ge {{u}_{*}}$. Integrating the inner integral by parts and using $F(0)=0$, we obtain:
\[{{\bar{\varphi }}^{k+1}}(u)\ge \int_{0}^{+\infty }{dK(t)}{{e}^{-L(U(u,t)-{{u}_{*}})}}L\int_{0}^{U(u,t)-{{u}_{*}}}{{{e}^{Lz}}F(z)}dz.\]
We transform
\begin{align}
  & L\int_{0}^{U(u,t)}{{{e}^{Lz}}F(z)}dz=L\int_{0}^{U(u,t)-{{u}_{*}}}{{{e}^{Lz}}(1-(1-F(z))}dz=
	\nonumber\\
	&\,\,\,\,\,\,\,\,\,\,\,\,\,\,\,\,\,\,\,\,\,\,\,\,\,\,
	={{e}^{L(U(u,t)-{{u}_{*}})}}-1-L\int_{0}^{U(u,t)-{{u}_{*}}}{{{e}^{Lz}}(1-F(z)})dz\ge 
	\nonumber \\ 
 & \,\,\,\,\,\,\,\,\,\,\,\,\,\,\,\,\,\,\,\,\,\,\,\,\,\,\ge {{e}^{L(U(u,t)-{{u}_{*}})}}-1-L\int_{0}^{+\infty }{{{e}^{Lz}}(1-F(z))}dz={{e}^{L(U(u,t)-{{u}_{*}})}}-\int_{0}^{+\infty }{{{e}^{Lz}}dF(z)}.  \nonumber
\end{align}
Thus, for $u\ge {{u}_{*}}$, there will be 
\begin{align} 
& {{{\bar{\varphi }}}^{k+1}}(u)\ge \int_{0}^{+\infty }{dK(t)}{{e}^{-L(U(u,t)-{{u}_{*}})}}({{e}^{L(U(u,t)-{{u}_{*}})}}-\int_{0}^{+\infty }{{{e}^{Lz}}dF(z)})=\nonumber\\
&\,\,\,\,\,\,\,\,\,\,\,\,\,\,=1-\int_{0}^{+\infty }{{{e}^{Lz}}dF(z)}\int_{0}^{+\infty }{{{e}^{-L(U(u,t)-{{u}_{*}})}}dK(t)}\ge \nonumber\\ 
& \,\,\,\,\,\,\,\,\,\,\,\,\,\,\ge 1-{{e}^{-L(u-{{u}_{*}})}}\int_{0}^{+\infty }{{{e}^{Lz}}dF(z)}\int_{0}^{+\infty }{{{e}^{-{{c}_{*}}Lt}}dK(t)}\ge 1-{{e}^{-L(u-{{u}_{*}})}}\,. \nonumber
\end{align}
The lemma is proven.
\end{proof}

%\begin{lemma}
{\bf Lemma 5.2.}
    Let Assumptions 5.1(a), (b) be satisfied. Then
    \begin{enumerate}
        \item[(i)] each function $\underline{\varphi}^k(u)$ is monotonically non-decreasing (increasing) in $u\ge 0$ and $0\le \underline{\varphi}^k(u)\le 1,\quad k\ge 0$;
        \item[(ii)] the sequence of values $\{\underline{\varphi}^k(u)\ge 0,\ k=0,1,\ldots\}$ increases monotonically with increasing $k$, and thus for all $k$, $\underline{\varphi}^k(u)\ge \varphi_*(u)=\max\{0,1-e^{-L(u-u_*)}\}$;
        \item[(iii)] the limit function $\underline{\varphi}(u)=\lim_{k\to\infty}\underline{\varphi}^k(u)$ is a solution of equation (98), $\underline{\varphi}(u)$ is monotone (non-decreasing) in $u$ and $1\ge \underline{\varphi}(u)\ge \varphi_*(u)=\max\{0,1-e^{-L(u-u_*)}\}$.
    \end{enumerate}
%\end{lemma}

\begin{proof}
The proof of item (i) is the same as in Lemma 5.1. From the proof of item (iv) of Lemma 5 it follows that ${{\underline{\varphi }}^{1}}(u)\ge 1-{{e}^{-Lu}}={{\underline{\varphi }}^{0}}(u)$. This, similar to the proof of item (ii) of Lemma 5.1, implies that ${{\underline{\varphi }}^{k+1}}(u)\ge {{\underline{\varphi }}^{k}}(u)\ge {{\underline{\varphi }}^{0}}(u)$ for any $k\ge 0$ and $u\ge 0$. The proof of (iii) is similar to the proof of Lemma 5.1(iv).
\end{proof}

%\begin{lemma}
{\bf Corollary 5.1.}
    If equation (\ref{(5.2)}) has a unique monotone solution $\varphi(u)$, then under conditions (a), $\bar{\varphi}(u)=\varphi(u)$, and under conditions (a), (b), $\bar{\varphi}(u)=\underline{\varphi}(u)=\varphi(u)$ and the boundary condition $\varphi(+\infty)=1$ is satisfied.
%\end{lemma}

The following two lemmas establish the contraction properties of the operator $A$ from (98). For a bounded function $\varphi(u),\ u\ge 0$, define the norm as $\|\varphi\|=\sup_{u\ge 0}|\varphi(u)|$.

%\begin{lemma}
{\bf Lemma 5.3.}
    (On the weak contraction property of the operator $A$ when $F(\cdot)<1$). For any monotone non-decreasing bounded functions $\varphi_1(u),\varphi_2(u)$ and for any $v>0$, the following holds:
    \begin{equation*}
        \sup_{u\in[0,v]}|A\varphi_1(u)-A\varphi_2(u)|\le q(v)\cdot \|\varphi_1-\varphi_2\|,
    \end{equation*}
    where $\|\varphi_1-\varphi_2\|=\sup_{u\in[0,+\infty)}|\varphi_1(u)-\varphi_2(u)|$, and
    \begin{equation*}
        q(v)=\int_{\,0}^{\,+\infty}{F(U(v,t))dK(t)}<1.
    \end{equation*}
%\end{lemma}

\begin{proof}
Let ${{\varphi }_{1}}(u),\,\,{{\varphi }_{2}}(u)$ be monotone (non-decreasing) bounded functions. For $u\in \left[ 0,v \right]$, due to the monotonicity of $U(\cdot ,t)$, the estimates hold
\begin{align}
& \left| A{{\varphi }_{1}}(u)-A{{\varphi }_{2}}(u) \right|\le \int\limits_{0}^{+\infty }{\int\limits_{0}^{U(u,t)}{\left| {{\varphi }_{1}}(U(u,t)-z)-{{\varphi }_{2}}(U(u,t)-z) \right|dF(z)}dK(t)}\le \nonumber\\ 
&  \ \ \ \ \ \ \ \ \ \ \ \ \ \ \ \le \left\| {{\varphi }_{1}}-{{\varphi }_{2}} \right\|\int\limits_{0}^{+\infty }{\int\limits_{0}^{U(u,t)}{dF(z)}dt}=\left\| {{\varphi }_{1}}-{{\varphi }_{2}} \right\|\int\limits_{0}^{+\infty }{\alpha {{e}^{-\alpha t}}F\left( U(u,t) \right)dK(t)}\le \nonumber\\ 
& \,\,\,\,\,\,\,\,\,\,\,\,\,\,\,\,\,\,\,\,\,\,\,\,\,\,\,\,\,\,\le \left\| {{\varphi }_{1}}-{{\varphi }_{2}} \right\|\int\limits_{0}^{+\infty }{F\left( U(v,t) \right)dK(t)}=q(v)\cdot \left\| {{\varphi }_{1}}-{{\varphi }_{2}} \right\|. \nonumber
\end{align}
Since $F(\cdot )$ and $U(\cdot ,t)$ are monotone (nondecreasing) and $F(\cdot )<1$, then $F\left( U(v,t) \right)<1$ and, consequently, $q(v)=\int\limits_{0}^{+\infty }{F\left( U(v,t) \right)dK(t)}<1$. The lemma is proved.
\end{proof}

%\begin{lemma}
{\bf Lemma 5.4.}
    (On the contraction property of the operator $A$ when $F(\cdot)<1$). For any monotone non-decreasing bounded functions $\varphi_1(u),\varphi_2(u)$ such that $\|A\varphi_1-A\varphi_2\|\ge \varepsilon>0$ and $0\le \varphi_*(u)\le A\varphi_1(u)\le 1$, $0\le \varphi_*(u)\le A\varphi_2(u)\le 1$, with a monotonically tending to one function $\varphi_*(u)$, the following holds:
    \begin{equation*}
        \|A\varphi_1-A\varphi_2\|\le q(\varepsilon)\cdot \|\varphi_1-\varphi_2\|,\quad q(\varepsilon)<1.
    \end{equation*}
%\end{lemma}

\begin{proof}
Let $\left\| A{{\varphi }_{1}}-A{{\varphi }_{2}} \right\|={{\lim }_{s\to \infty }}\left| A{{\varphi }_{1}}({{u}^{s}})-A{{\varphi }_{2}}({{u}^{s}}) \right|\ge \varepsilon >0$ with some sequence $\left\{ {{u}^{s}} \right\}$. Obviously $\left| A{{\varphi }_{1}}(u)-A{{\varphi }_{2}}(u) \right|\le 1-{{\varphi }_{*}}(u)$. We define $${{u}^{*}}(\varepsilon )={{\sup }_{v\ge 0}}\left\{ v:\,\,1-{{\varphi }_{*}}(v)\ge {\varepsilon }/{2}\; \right\}<+\infty. $$ Then for sufficiently large $s$ we have $1-{{\varphi }_{*}}({{u}^{s}})\ge \left| A{{\varphi }_{1}}({{u}^{s}})-A{{\varphi }_{2}}({{u}^{s}}) \right|\ge {\varepsilon }/{2}\;$ and, consequently, ${{u}^{s}}\le {{u}^{*}}(\varepsilon )$. Fair assessments
\begin{align} 
& \left| A{{\varphi }_{1}}({{u}^{s}})-A{{\varphi }_{2}}({{u}^{s}}) \right|\le \int\limits_{0}^{+\infty }{\int\limits_{0}^{U({{u}^{s}},t)}{\left| {{\varphi }_{1}}(U({{u}^{s}},t)-z)-{{\varphi }_{2}}(U({{u}^{s}},t)-z) \right|dF(z)}dK(t)}\le \nonumber\\ 
& \ \ \ \ \ \ \ \ \ \ \ \ \ \ \ \ \ \ \le \left\| {{\varphi }_{1}}-{{\varphi }_{2}} \right\|\int\limits_{0}^{+\infty }{\int\limits_{0}^{U({{u}^{s}},t)}{dF(z)}dK(t)}=\left\| {{\varphi }_{1}}-{{\varphi }_{2}} \right\|\int\limits_{0}^{+\infty }{F\left( U({{u}^{s}},t) \right)dK(t)}\le \nonumber\\ 
& \,\,\,\,\,\,\,\,\,\,\,\,\,\,\,\,\,\,\,\,\,\,\,\,\,\,\,\,\,\,\le \left\| {{\varphi }_{1}}-{{\varphi }_{2}} \right\|\int\limits_{0}^{+\infty }{F\left( U({{u}^{*}}(\varepsilon ),t) \right)dK(t)}=q(\varepsilon )\cdot \left\| {{\varphi }_{1}}-{{\varphi }_{2}} \right\|, \nonumber
\end{align}
where $q(\varepsilon )=\int\limits_{0}^{+\infty }{F\left( U({{u}^{*}}(\varepsilon ),t) \right)dK(t)}<1$. The lemma is proven.
\end{proof}

%\begin{corollary}
{\bf Corollary 5.2.}
    Let Assumptions 5.1 be satisfied. If in Lemma 5.4 the functions $\varphi_1(u)\ne \varphi_2(u)$ satisfy $\varphi_*(u)\le \varphi_1(u)\le 1$, $\varphi_*(u)\le \varphi_2(u)\le 1$, where $\varphi_*(u)=\max\{0,1-e^{-L(u-u_*)}\}$, then
    \begin{equation*}
        \|A\varphi_1-A\varphi_2\|<\|\varphi_1-\varphi_2\|.
    \end{equation*}
%\end{corollary}

\begin{proof}
If $\left\| A{{\varphi }_{1}}-A{{\varphi }_{2}} \right\|=0$, then the statement of the corollary is obvious. Let $\left\| A{{\varphi }_{1}}-A{{\varphi }_{2}} \right\|={{\lim }_{s\to \infty }}\left| A{{\varphi }_{1}}({{u}^{s}})-A{{\varphi }_{2}}({{u}^{s}}) \right|>0$ with some sequence $\left\{ {{u}^{s}} \right\}$. From the proof of Lemma 5.1 (iv) it follows that
${{\varphi }_{*}}(u)\le A{{\varphi }_{*}}(u)\le A{{\varphi }_{1}}(u)\le 1$, ${{\varphi }_{*}}(u)\le A{{\varphi }_{*}}(u)\le A{{\varphi }_{2}}(u)\le 1$.
Since ${{\lim }_{u\to +\infty }}{{\varphi }_{*}}(u)=1$, it follows that the sequence $\left\{ {{u}^{s}} \right\}$ is bounded. Then the assertion of the corollary follows from Lemma 5.4.
\end{proof}

%\begin{lemma}
{\bf Lemma 5.5.}
    (On the contraction property of the operator $A$ when $F(\cdot)\le 1$). Let $K(\cdot)<1$ and $F(z)=1$ for $z\ge \bar{z}$. Then for any monotone non-decreasing bounded functions $\varphi_1(u),\varphi_2(u)$ such that $\|\varphi_1-\varphi_2\|\ge \varepsilon>0$ and $0\le \varphi_*(u)\le \varphi_1(u)\le 1,\ 0\le \varphi_*(u)\le \varphi_2(u)\le 1$ with a monotonically tending to one function $\varphi_*(u)$, the following holds:
    \begin{equation*}
        \|A\varphi_1-A\varphi_2\|\le q(\varepsilon)\cdot \|\varphi_1-\varphi_2\|,
    \end{equation*}
    where $q(\varepsilon)<1$, $\|\varphi_1-\varphi_2\|=\sup_{u\in[0,+\infty)}|\varphi_1(u)-\varphi_2(u)|$, $\|A\varphi_1-A\varphi_2\|=\sup_{u\in[0,+\infty)}|A\varphi_1(u)-A\varphi_2(u)|$.
%\end{lemma}

\begin{proof}
Obviously, $\left| {{\varphi }_{1}}(u)-{{\varphi }_{2}}(u) \right|\le 1-{{\varphi }_{*}}(u)$. Then for any $u\ge 0$ the following estimates are valid:
\begin{align}
  & \left| A{{\varphi }_{1}}(u)-A{{\varphi }_{2}}(u) \right|\le \int\limits_{0}^{+\infty }{\int\limits_{0}^{U(u,t)}{\left| {{\varphi }_{1}}(U(u,t)-z)-{{\varphi }_{2}}(U(u,t)-z) \right|dF(z)}dK(t)}\le  \nonumber\\ 
 & \,\,\,\,\,\,\,\,\,\,\,\,\,\,\,\,\,\,\,\,\,\le \int\limits_{0}^{+\infty }{\int\limits_{0}^{U(u,t)}{\min \left\{ \left\| {{\varphi }_{1}}-{{\varphi }_{2}} \right\|,1-{{\varphi }_{*}}(U(u,t)-z) \right\}dF(z)}dK(t)}\le \nonumber \\ 
 & \,\,\,\,\,\,\,\,\,\,\,\,\,\,\,\,\,\,\,\,\,\le \int\limits_{0}^{+\infty }{\int\limits_{0}^{\min \left\{ U(u,t),\bar{z} \right\}}{\min \left\{ \left\| {{\varphi }_{1}}-{{\varphi }_{2}} \right\|,1-{{\varphi }_{*}}(U(u,t)-z) \right\}dF(z)}dK(t)}\le 
\nonumber 
\end{align}
\begin{align}
  & \,\,\,\,\,\,\,\,\,\,\,\,\,\,\,\,\,\,\,\,\,\le \int\limits_{0}^{+\infty }{\int\limits_{0}^{\min \left\{ U(u,t),\bar{z} \right\}}{\min \left\{ \left\| {{\varphi }_{1}}-{{\varphi }_{2}} \right\|,1-{{\varphi }_{*}}(\max \left\{ 0,U(u,t)-\bar{z} \right\}) \right\}dF(z)}dK(t)}\le \nonumber \\ 
 & \,\,\,\,\,\,\,\,\,\,\,\,\,\,\,\,\,\,\,\,\,\,\le \int\limits_{0}^{+\infty }{\min \left\{ \left\| {{\varphi }_{1}}-{{\varphi }_{2}} \right\|,1-{{\varphi }_{*}}(\max \left\{ 0,U(u,t)-\bar{z} \right\}) \right\}F\left( \min \left\{ U(u,t),\bar{z} \right\} \right)dK(t)}\le \nonumber \\ 
 & \,\,\,\,\,\,\,\,\,\,\,\,\,\,\,\,\,\,\,\,\,\,\le \int\limits_{0}^{+\infty }{\min \left\{ \left\| {{\varphi }_{1}}-{{\varphi }_{2}} \right\|,1-{{\varphi }_{*}}(\max \left\{ 0,U(0,t)-\bar{z} \right\}) \right\}dK(t)}\le \nonumber \\ 
 & \,\,\,\,\,\,\,\,\,\,\,\,\,\,\,\,\,\,\,\,\,\,\le \left\| {{\varphi }_{1}}-{{\varphi }_{2}} \right\|\cdot \int\limits_{0}^{+\infty }{\min \left\{ 1,{\left( 1-{{\varphi }_{*}}(\max \left\{ 0,U(0,t)-\bar{z} \right\}) \right)}/{\varepsilon }\; \right\}dK(t)}. \nonumber 
\end{align}
But since $U(0,t)\to +\infty $, then $\left( 1-{{\varphi }_{*}}(\max \left\{ 0,U(0,t)-\bar{z} \right\}) \right)\to 0$ as $t\to +\infty $ and, thus, $\min \left\{ 1,{\left( 1-{{\varphi }_{*}}(\max \left\{ 0,U(0,t)-\bar{z} \right\}) \right)}/{\varepsilon }\; \right\}<1$ for all sufficiently large $t$. It follows that
\[q(\varepsilon )=\int\limits_{0}^{+\infty }{\min \left\{ 1,{\left( 1-{{\varphi }_{*}}(\max \left\{ 0,U(0,t)-\bar{z} \right\}) \right)}/{\varepsilon }\; \right\}dK(t)}<1.\]
The lemma is proved.
\end{proof}

%\begin{corollary}
{\bf Corollary 5.3.}
    Let Assumptions 5.1 be satisfied. From Lemmas 5.4 and 5.5, it follows that the operator $A$ defined in (98) is a contraction on the set of functions $\varphi(u)$ such that $\varphi_*(u)=\max\{0,1-e^{-L(u-u_*)}\}\le \varphi(u)\le 1$, i.e., for any functions $\varphi_1(u)\ne \varphi_2(u)$ such that $\varphi_*(u)\le \varphi_1(u)\le 1,\ \varphi_*(u)\le \varphi_2(u)\le 1$, the following holds:
    \begin{equation*}
        \|A\varphi_1-A\varphi_2\|<\|\varphi_1-\varphi_2\|.
    \end{equation*}
%\end{corollary}

%\begin{theorem}
{\bf Theorem 5.1.}
    Assume that Assumptions 5.1 are satisfied. Then
    \begin{enumerate}
        \item[1)] equation (\ref{(5.2)}) has a monotone solution $\varphi(u)$ such that
        $1\ge \varphi(u)\ge \varphi_*(u)=\max\{0,1-e^{-L(u-u_*)}\}$;
        \item[2)] $\varphi(u)$ is the unique monotone solution of problems (\ref{(5.2)}), (\ref{(5.3)});
        \item[3)] the sequence of approximations (101) starting with $\varphi^*(u)\equiv 1$ converges monotonically from above, and starting with $\varphi_*(u)=\max\{0,1-e^{-L(u-u_*)}\}$, converges monotonically from below to the solution $\varphi(u)$;
        \item[4)] for any initial approximation $\varphi^0(u)$ such that
        $1\ge \varphi^0(u)\ge \varphi_*(u)=\max\{0,1-e^{-L(u-u_*)}\}$,
        the corresponding sequence of approximations (\ref{(5.5)}) converges pointwise to the solution $\varphi(u)$.
    \end{enumerate}
%\end{theorem}

\begin{proof}
Statements 1), 3) follow from Lemmas 5.1, 5.2.

Let us prove assertion 2) of Theorem 5.1 on the uniqueness of the solution to problem (5.2), (5.3). Assume the contrary: let there be two solutions ${{\varphi }_{1}}\ne {{\varphi }_{2}}$ (i.e. ${{\varphi }_{1}}(u)\ne {{\varphi }_{2}}(u)$ for some $u\ge 0$) of problem (5.2), (5.3). Put ${{\varphi }^{0}}(u)=\max \left\{ {{\varphi }_{1}}(u),{{\varphi }_{2}}(u) \right\}$. Obviously, ${{\varphi }^{0}}\ne {{\varphi }_{1}}$, ${{\varphi }^{0}}\ne {{\varphi }_{2}}$. Consider the sequence $\left\{ {{\varphi }^{k}}(u):=A{{\varphi }^{k}}(u),\,\,k=0,1,... \right\}$. Obviously, $1\ge {{\varphi }^{0}}(u)\ge {{\varphi }_{1}}(u)$ and $1\ge {{\varphi }^{0}}(u)\ge {{\varphi }_{2}}(u)$. Since the operator $A$ takes a large function to a large one (is monotone), then, taking into account Lemma 5.1, $1\ge {{\varphi }^{1}}(u)=A{{\varphi }^{0}}(u)\ge A{{\varphi }_{1}}(u)={{\varphi }_{1}}(u)$, $1\ge {{\varphi }^{1}}(u)=A{{\varphi }^{0}}(u)\ge A{{\varphi }_{2}}(u)={{\varphi }_{2}}(u)$ and, thus, ${{\varphi }^{1}}(u)\ge \max \left\{ {{\varphi }_{1}}(u),{{\varphi }_{2}}(u) \right\}={{\varphi }^{0}}(u)$. From this it follows by induction that ${{\varphi }^{k+1}}(u)\ge {{\varphi }^{k}}(u)\ge {{\varphi }^{0}}(u)$ for all $k\ge 0$. The sequence $\left\{ {{\varphi }^{k}}(\cdot ) \right\}$ is monotonically increasing and is bounded above by unity, so it has a pointwise limit $\varphi (u)$, $1\ge \varphi (u)\ge {{\varphi }^{0}}(u)$, which is a solution of equation (5.2) by Lebesgue's theorem on the limit passage under the integral sign. Obviously, the obtained function $\varphi (u)$ is a solution of problem (5.2), (5.3), different from ${{\varphi }_{1}}$ and ${{\varphi }_{2}}$. Thus, two different functions ${{\varphi }_{1}}(u)$ and $\varphi (u)\ge {{\varphi }_{1}}(u)$ are solutions of equation (5.2) and satisfy the conditions of Lemmas 5.4, 5.5 with ${{\varphi }_{*}}\left( u \right)={{\varphi }_{1}}(u)$. But this contradicts the contraction property of the operator $A$ from Lemma 5.4 or 5.5.

Statement 4) of the theorem is proved similarly. Let $1\ge {{\varphi }^{0}}(u)\ge \max \left\{ 0,1-{{e}^{-L(u-{{u}_{*}})}} \right\}$, i.e. ${{\bar{\varphi }}^{0}}(u)\ge {{\varphi }^{0}}(u)\ge {{\underline{\varphi }}^{0}}(u)$. Then ${{\bar{\varphi }}^{1}}(u)\ge {{\varphi }^{1}}(u)\ge \underline{{{\varphi }^{1}}}(u)$ and, moreover, for any $k$ ${{\bar{\varphi }}^{k}}(u)\ge {{\varphi }^{k}}(u)\ge \underline{{{\varphi }^{k}}}(u)$. It follows that there exists a limit $\varphi (u)={{\lim }_{k\to \infty }}{{\varphi }^{k}}(u)\equiv \bar{\varphi }(u)\equiv \underline{\varphi }(u)$. The theorem is proved.
\end{proof}

%\begin{theorem}
{\bf Theorem 5.2.}
    (On uniform convergence and the rate of convergence of the method of successive approximations). Under Assumptions 5.1, the method of successive approximations (101), starting from an initial approximation $\varphi^0(u)$ such that $1\ge \varphi^0(u)\ge \max\{0,1-e^{-L(u-u_*)}\}$, converges uniformly monotonically to the solution of problems (98), (99); moreover, it converges to any $\varepsilon$-neighborhood of the solution of problem (98) at the rate of a geometric progression with a denominator depending on $\varepsilon$.
%\end{theorem}

The statement of the theorem is similar to Theorem 4.2 and follows from Theorem 5.1 and Lemmas 5.4, 5.5.

\section{Numerical experiments}

Let $U(u,t)=u+ct$, $x(z)=z$, and the claim distribution $F(z)$ be a mixture of two exponential distributions:
\begin{equation*}
    F(z)=\varepsilon(1-\exp(-\varepsilon z))+(1-\varepsilon)(1-\exp(-z)),\quad \mu=\int_{\,0}^{\,+\infty}{(1-F(z))dz}=2-\varepsilon.
\end{equation*}
Assume that random claims arrive with probability one at equal deterministic intervals (at the end of the calculation period, say, a month or a year) of length $T$. In this case, $F(z)$ is the distribution function of the sum of insurance claims arriving over a time interval of length $T$. Take $\varepsilon=0.1$, $cT=5$.

Then the iterative process (\ref{(5.5)}) for solving equation (\ref{(5.2)}) can be rewritten in the form:
\begin{equation}\label{(5.6)}
    \varphi^{k+1}(u)=\int\limits_{0}^{u+cT}{\varphi^k(u+cT-z)dF(z)},\quad k=0,1,\ldots
\end{equation}

According to Lemma 5.1, the sequence $\{\bar{\varphi}^k(u),\ k=0,1,\ldots\}$ generated by (\ref{(5.5)}) starting from the initial approximation $\varphi^0(u)=\bar{\varphi}^0(u)\equiv 1$ converges monotonically from above to the solution $\varphi(u)$ of equation (\ref{(5.2)}), and according to Lemma 5.2, the sequence $\{\underline{\varphi}^k(u),\ k=0,1,\ldots\}$ generated by (\ref{(5.5)}) starting from the initial function $\varphi^0(u)=\underline{\varphi}^0(u)=1-e^{-Lu}$ converges monotonically from below to $\varphi(u)$. The course of these iterations is shown in Fig. 5.1. The error of the method $\Delta_k=\sup_{0\le u<+\infty}(\bar{\varphi}^k(u)-\underline{\varphi}^k(u))$ decreases monotonically with increasing number of iterations $k$, as shown in Fig. 5.2. As a result, we have $\varphi(0)\approx 0.801932$, $\Delta_{100}\le 10^{-7}$.

\textbf{Comparison with known approximations.} An alternative to the numerical solution of equations can be the method (Monte Carlo) of repeated simulation of the random risk process for a specific value of initial capital and counting the proportion of ruined trajectories. To assess the accuracy and complexity of the Monte Carlo method, the probability of non-ruin with zero initial capital was estimated by repeated (100000 trajectories) simulation of the random risk process (\ref{(5.1)}). The value of the probability obtained was 0.802452. The standard deviation (0.00137612) was also estimated by repeated (100 times) calculation of this probability. One estimate of the probability required 20 minutes of machine time (Intel Pentium M 1.6). For comparison: the method of successive approximations required 34 iterations and 1 minute to find the probability of non-ruin for all values of initial capital with an accuracy of 0.001 and 100 iterations to achieve an accuracy of $10^{-10}$ (the calculations used such an accuracy of numerical integration and interpolation in (\ref{(5.6)}) that allows finding the solution of equation (\ref{(5.2)}) with an accuracy of $10^{-10}$), which is completely unattainable for Monte Carlo-type methods. To achieve greater accuracy, more accurate numerical integration and interpolation must be used, which will lead to longer calculation times depending on the required accuracy. When finding the solution with an accuracy of $10^{-10}$, we can actually estimate the behavior of the ruin probability for capital $u\le 260$, since for $u\ge 260$, $\psi(u)=1-\phi(u)\le 10^{-10}$.

In practice, various approximations of the ruin probability are also used [35, Section 3, §5]. For comparison with our method, the Cramér-Lundberg bound $\psi_{KL}(u)=e^{-0.0806959 u}$ (Fig. 5.3), the diffusion approximation $\psi_{Dif}(u)=e^{-0.458716 u}$ (Fig. 5.3), the Beekman-Bowers approximation $\psi_{BB}(u)=0.139175\int_{-\infty}^{u}{x^{-0.561746}e^{-0.0553572u}dx}$ (Fig. 5.4), and the De Vylder approximation $\psi_{DV}(u)=0.190611413261459 e^{-0.087436428101587u}$ (Fig. 5.4) were constructed. Comparison of the approximations with the exact solution $\psi(u)$ obtained by the method of successive approximations shows that for this example and $u\in[0,260]$:

\begin{itemize}
    \item the diffusion approximation is completely unsuitable for estimating the ruin probability;
    \item the Cramér-Lundberg bound overestimates the ruin probability by about 5 times;
    \item the Beekman-Bowers approximation gives a relative error of up to 30% for capital up to 50 and greatly overestimates the ruin probability for large initial capitals (Fig. 5.5);
    \item the De Vylder approximation turned out to be the best; it underestimates the ruin probability by a maximum of 70% (Fig. 5.6).
\end{itemize}

\begin{figure}[h]
    \centering
    \rotatebox{90}{Survival probability}
    \includegraphics[width=0.8\textwidth]{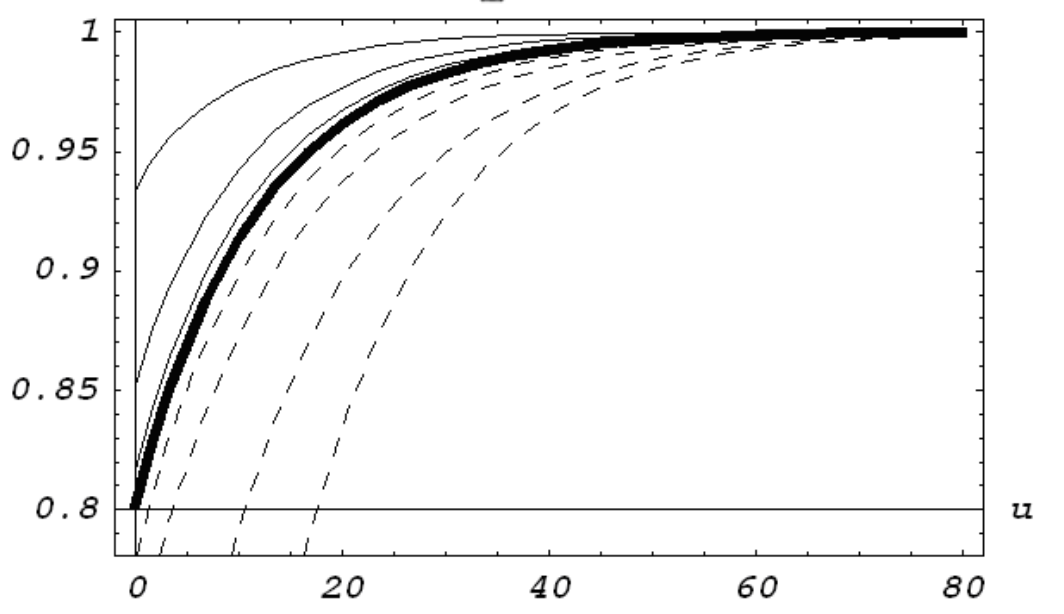}
    \caption{Iteration of $\overline{\Phi}^{k}(u)$ and $\underline{\Phi}^{k}(u)$, $k=2,5,10,50$.}
    \label{fig:5.1}
\end{figure}

\begin{figure}[h]
    \centering
        \rotatebox{90}{Method error}
    \includegraphics[width=0.8\textwidth]{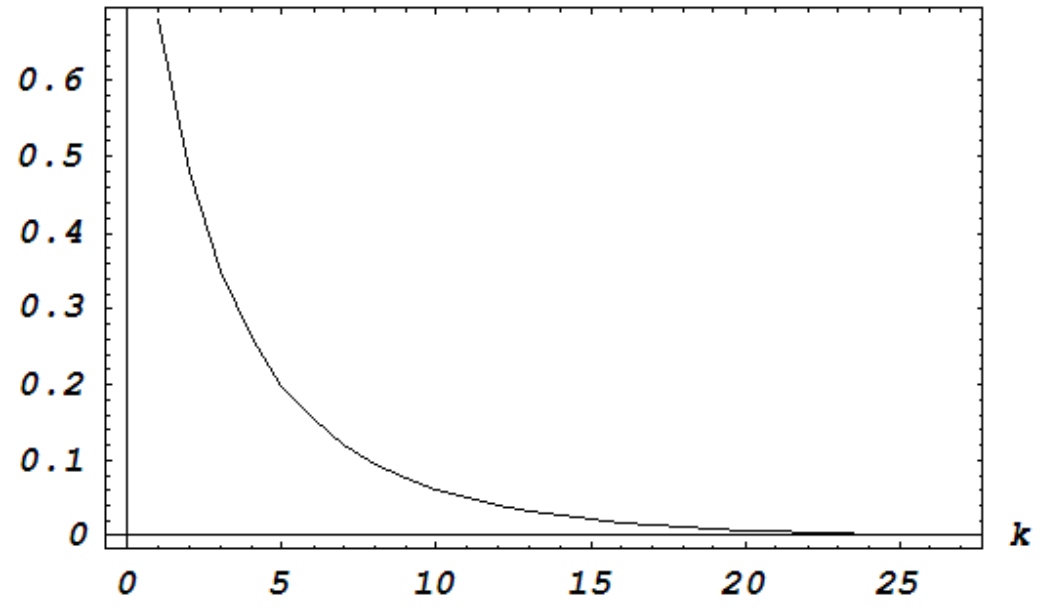}
    \caption{Dependence of the error on the number of iterations.}
    \label{fig:5.2}
\end{figure}

\begin{figure}[h]
    \centering
        \rotatebox{90}{Ruin probability}
    \includegraphics[width=0.8\textwidth]{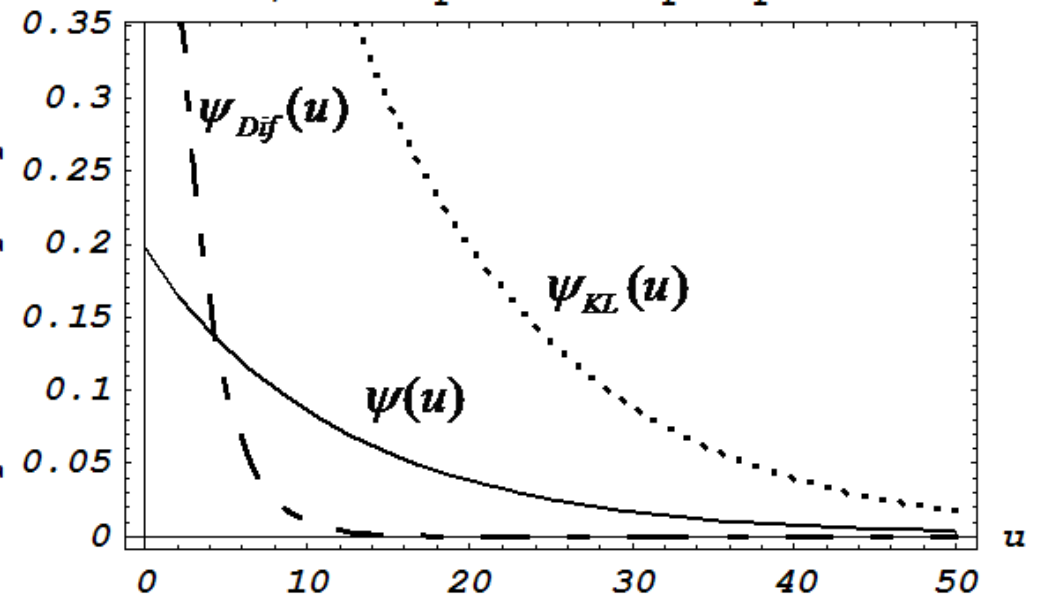}
    \caption{Estimate of the ruin probability.}
    \label{fig:5.3}
\end{figure}

\begin{figure}[h]
    \centering
            \rotatebox{90}{Ruin probability}
    \includegraphics[width=0.8\textwidth]{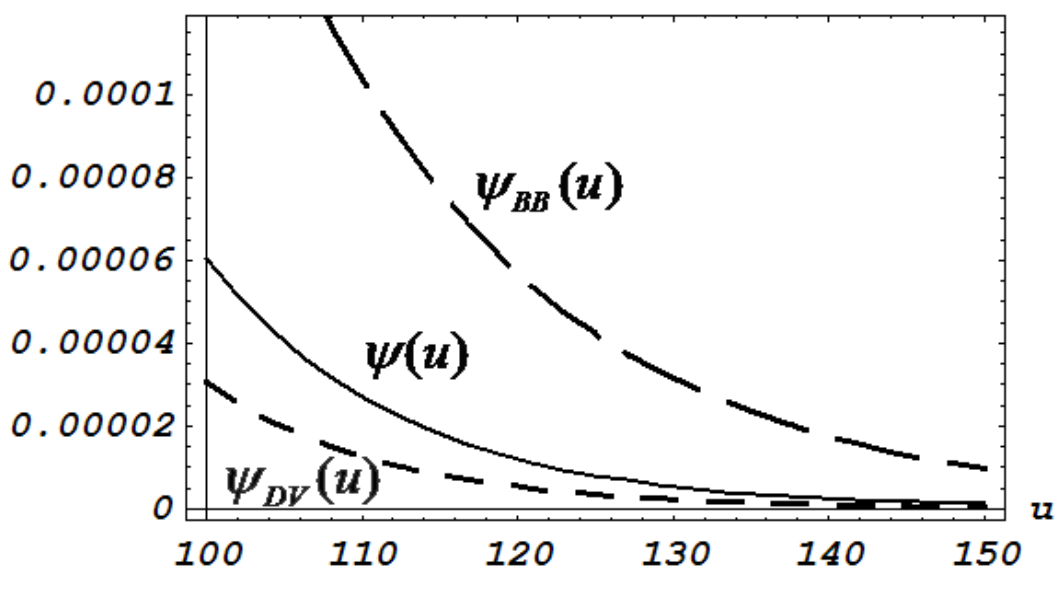}
    \caption{Estimate of the ruin probability.}
    \label{fig:5.4}
\end{figure}

\begin{figure}[h]
    \centering
    \rotatebox{90}{$\psi_{BB}(u)/\psi(u)$}
    \includegraphics[width=0.8\textwidth]{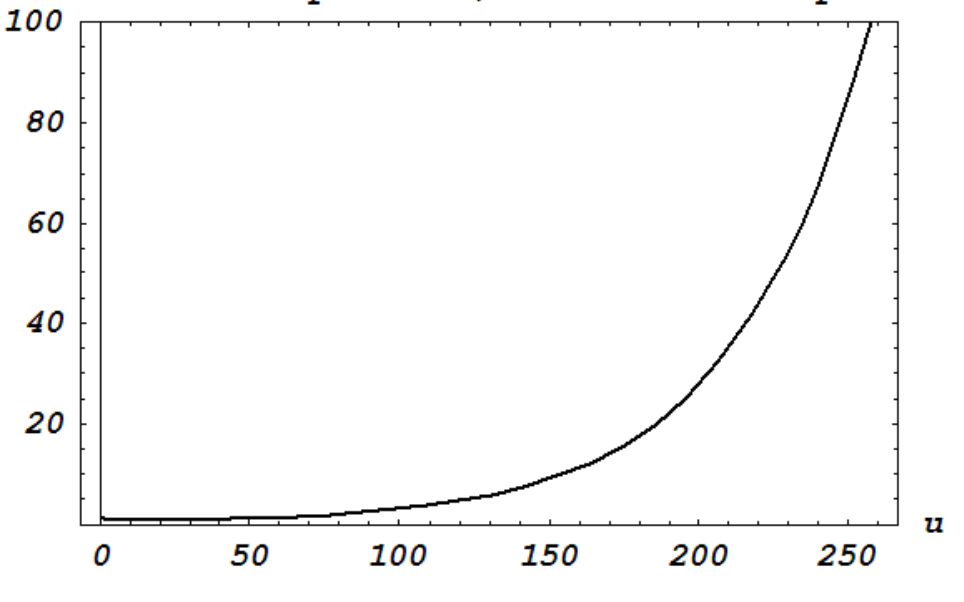}
    \caption{Beekman--Bowers approximation error.}
    \label{fig:5.5}
\end{figure}

\begin{figure}[h]
    \centering
            \rotatebox{90}{Relative error}
    \includegraphics[width=0.8\textwidth]{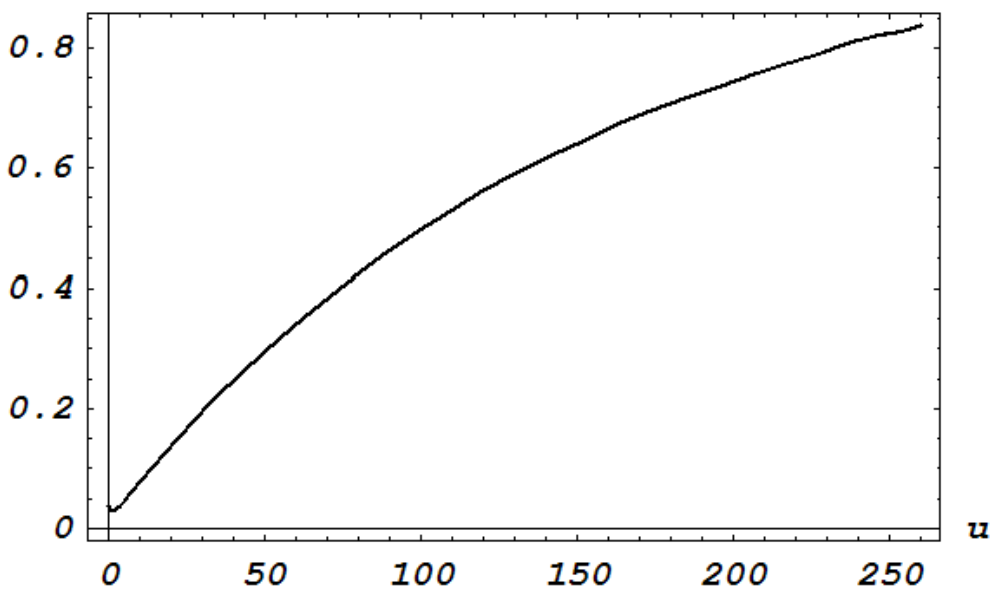}
    \caption{De Vylder approximation error.}
    \label{fig:5.6}
\end{figure}

\section{Conclusions to the section}

In this section, the application of the method of successive approximations (\ref{(5.5)}) for solving general integral equations (\ref{(2.11)}), (\ref{(5.2)}) for the probability of ruin of renewal risk processes (\ref{(2.1)}), (\ref{(5.1)}) describing the random evolution of the capital of an insurance company is justified. The method converges to the solution of equation (\ref{(5.2)}) under the condition that there exists an analogue of the Lundberg constant for the considered process, which is a positive solution of inequality (\ref{(5.4)}). Thus, the existence of a solution to the integral equation (\ref{(5.2)}) taking the value one at infinity is constructively proved. The uniqueness of such a solution is proved under the additional assumption that the claim arrival time or their size can be arbitrarily large. The domain of attraction of the solution is bounded from above by the unit function and from below by the generalized Cramér-Lundberg bound.

\chapter{APPLICATION OF THE METHOD OF SUCCESSIVE APPROXIMATIONS FOR FINDING THE PROBABILITY OF NON-RUIN OF AN INSURANCE COMPANY IN THE PRESENCE OF RANDOM PREMIUMS}

In this section, a risk process is studied that describes the evolution of the capital of an insurance company in the presence of random premiums and claims. Integral equations are derived that describe the probability of non-ruin of the company as a function of the company's initial capital. Necessary and sufficient conditions for the existence and general sufficient conditions for the uniqueness of solutions of these integral equations are established, as well as conditions for the convergence of the method of successive approximations for finding solutions. The results of this section are published in [45, 107, 108].

\section{Description of the model in the presence of stochastic premiums}

In this section, a risk process is considered in which both claims and premiums are stochastic [22]. It describes the stochastic process of customer arrivals to an insurance company more adequately than the classical Cramér-Lundberg model. In [5], for a similar process with independent Poisson streams of random premiums and claims, an integral equation for the probability of non-ruin was derived and some analytical solutions of this equation were found. In this section, we derive more general integral equations for the probability of non-ruin in the case of possibly non-Poisson streams of random premiums and claims and in the presence of deterministic payments, and we justify the method of successive approximations for their numerical solution. Namely, we establish necessary and sufficient conditions for existence and general sufficient conditions for uniqueness of solutions of the equations, prove the convergence of the method of successive approximations, and estimate the rate of convergence.

Consider a random process (risk) $\xi_t$ describing the evolution over time $t$ of the capital of an insurance company and satisfying the stochastic integral equation
\begin{equation}\label{(6.1)}
    \xi_t = u + \int_0^t c(\xi_s) ds + P_t - S_t, \quad t \ge 0,
\end{equation}
where $u$ is the initial capital of the company, $\xi_0 = u$; $c(\xi)$ is the intensity of deterministic payments depending on the current capital (a piecewise continuous function); $S_t = \sum_{k=1}^{N_t^c} z_k$ are the aggregated insurance claims; $P_t = \sum_{k=1}^{N_t^p} p_k$ are the aggregated insurance premiums; $z_k$ are independent random variables (claims) with common distribution function $F(z)$; $p_k$ are independent random variables (premiums) with common distribution function $G(z)$; $N_t^c$ is the number of claims received by time $t$, in particular, it can be a Poisson process with intensity $\alpha_c$; similarly, $N_t^p$ is the number of premiums received by time $t$, in particular, it can be a Poisson process with intensity $\alpha_p$. Consider the probability of ruin of an insurance company $\psi(u) = P\{\exists t \ge 0 : \xi_t < 0\}$ on the infinite time interval $t \in [0, +\infty)$ and the corresponding probability of non-ruin $\varphi(u) = 1 - \psi(u)$ as functions of the initial capital $u \ge 0$.

{\bf Theorem 6.1.}
Let the number of premiums $N_t^p$ and claims $N_t^c$ in (6.1) be given by independent Poisson streams with intensities $\alpha_p$ and $\alpha_c$, respectively. Then the probability of non-ruin $\varphi(u)$ satisfies the integro-differential equation
\begin{equation}\label{(6.2)}
    (\alpha_c + \alpha_p)\varphi(u) - c(u)\varphi'(u) = \alpha_c \int_0^u \varphi(u-z) dF(z) + \alpha_p \int_0^{+\infty} \varphi(u+z) dG(z).
\end{equation}

In [5], Theorem 6.1 was established for the case $c(u) \equiv 0$. Next, a more general equation is derived for the case of non-Poisson streams of claims and premiums.

Consider the situation when there is no deterministic component, i.e., $c(u) \equiv 0$. Then equation (6.2) becomes integral; rewriting it in the case of random premiums, we get:
\begin{equation}\label{(6.3)}
    \varphi(u) = A_1 \varphi(u),
\end{equation}
where the linear integral operator $A_1$ is given by:
\begin{equation}\label{(6.4)}
    A_1 \varphi(u) := \frac{\alpha_c}{\alpha_c + \alpha_p} \int_0^u \varphi(u-z) dF(z) + \frac{\alpha_p}{\alpha_c + \alpha_p} \int_0^{+\infty} \varphi(u+z) dG(z), \quad u \ge 0.
\end{equation}

Consider a more general process constructed as follows. At the initial moment $t = 0$ and at each moment $t'$ of a jump of the process (i.e., arrival of a random premium or claim), two random variables, $\tau_p$ and $\tau_c$, with distribution functions $K_p(t)$ and $K_c(t)$, respectively, are observed. The next jump occurs at the moment $t'' = t' + \min\{\tau_p, \tau_c\}$. If $\tau_p < \tau_c$, then the capital increases by a random premium with distribution function $G(\cdot)$; otherwise, it decreases by a claim with distribution $F(\cdot)$. In the intervals between jumps, the capital changes continuously and deterministically according to the equation: $\frac{d\xi_t}{dt} = c(\xi_t)$.

The described process includes the case of Poisson processes of random premium and claim arrivals.

%\begin{assumption}
{\bf Assumption 6.1.}
    The functions $K_p(t)$ and $K_c(t)$ have no common discontinuity points, and $K_p(t) < 1, K_c(t) < 1$ for all $t \ge 0$.
%\end{assumption}

In the case of Poisson streams of premiums and claims, $K_p(t) = 1 - e^{-\lambda_p t}$, $K_c(t) = 1 - e^{-\lambda_c t}$, and thus Assumption 6.1 is automatically satisfied.

Define the function $U(u,t)$ as the solution of the Cauchy problem: $\frac{dU}{dt} = c(U)$, $U(u,0) = u$, where $c(\cdot)$ is the intensity of the deterministic drift in equation (6.1).

Note that the function $U(u,t)$ is monotonically non-decreasing in $u$.

The following statement holds.

%\begin{theorem}
{\bf Theorem 6.2.}
    In the general case of non-Poisson streams of premiums and claims, the probability of non-ruin $\varphi(u)$ satisfies the integral equation
    \begin{equation}\label{(6.5)}
        \varphi(u) = A_2 \varphi(u),
    \end{equation}
    where
    \begin{align}
        A_2 \varphi(u)
        & = \int_0^\infty \int_0^{U(u,t)} \varphi(U(u,t)-z) dF(z) (1 - K_p(t)) dK_c(t) \nonumber \\
        & \quad + \int_0^\infty \int_0^{+\infty} \varphi(U(u,t)+z) dG(z) (1 - K_c(t)) dK_p(t).
				\label{(6.6)}
    \end{align}
%\end{theorem}
\begin{proof}
    Let the first jump of the process (arrival of a claim or premium) occur at time $t$. For the company not to be ruined, it is necessary that it is not ruined before the first jump, at the moment of the first jump, and throughout the subsequent time. Three cases are possible.

    \begin{enumerate}
        \item[1)] A claim arrives first, which does not ruin the company (the amount is less than $U(u,t)$), after which the company is not ruined on the time interval $(t, +\infty)$. The probability of this event is equal to
        \begin{equation*}
            \int_0^\infty \int_0^{U(u,t)} \varphi(U(u,t)-z) dF(z) (1 - K_p(t)) dK_c(t).
        \end{equation*}

        \item[2)] A premium arrives first, after which the company is not ruined on the time interval $(t, +\infty)$. The probability of this event is equal to
        \begin{equation*}
            \int_0^\infty \int_0^{+\infty} \varphi(U(u,t)+z) dG(z) (1 - K_c(t)) dK_p(t).
        \end{equation*}

        \item[3)] The first claim and the first premium arrive simultaneously. The probability of this event, by Assumption 6.1, is equal to zero.
    \end{enumerate}
    By the law of total probability, the probability of non-ruin of the company can be represented as the sum of the probabilities of these disjoint events.
\end{proof}

Consider the integral equation
\begin{equation}\label{(6.7)}
    \varphi(u) = A \varphi(u),
\end{equation}
where the operator $A$ can be either $A_1$ or $A_2$. Here the function $\varphi(\cdot)$ is assumed to be monotone, and the integrals are understood in the Lebesgue-Stieltjes sense. This is a linear homogeneous integral equation with integrands defined on the unbounded domain $[0, +\infty)$ and non-negative kernels. It always has a trivial (zero) solution. We are interested in a non-decreasing in $u$ solution $\varphi(u)$, $0 \le \varphi(u) \le 1$, satisfying the following boundary condition at infinity:
\begin{equation}\label{(6.8)}
    \varphi(+\infty) = \lim_{u \to +\infty} \varphi(u) = 1.
\end{equation}
This condition means that with unlimited initial capital, the company does not go bankrupt.

\section{Properties of Integral Operators}

In this section, the properties of the operators $A_1, A_2$ from equations (13), (14) are studied. Let us define the (metric) space $H$ of functions $\varphi(u)$ non-decreasing in $u \in [0,+\infty)$ such that $0 \le \varphi(u) \le 1$, with the distance between functions $\varphi_1, \varphi_2 \in H$
\begin{equation*}
    \rho(\varphi_1, \varphi_2) := \sup_{u \ge 0} |\varphi_1(u) - \varphi_2(u)|.
\end{equation*}
Define a partial order on $H$: $\varphi_1 \le \varphi_2$ if $\varphi_1(u) \le \varphi_2(u)$ for any $u \ge 0$.

%\begin{lemma}
{\bf lemma 6.1.}
    The operators $A_1, A_2$ are continuous on $H$ with respect to pointwise convergence.
%\end{lemma}

\begin{proof}
    Let the sequence $\{\varphi^k(u) \in H, k = 0,1,\ldots\}$ converge pointwise to some function $\varphi(u)$. Obviously, $\varphi(u) \in H$. By the Lebesgue dominated convergence theorem, the limit can be taken under the sign of the corresponding integral operator: $\lim_{k \to \infty} A_i \varphi^k = A_i \lim_{k \to \infty} \varphi^k = A_i \varphi$, $i = 1,2$.
\end{proof}

%\begin{lemma}
{\bf Lemma 6.2.}
    The linear integral operator $A_i, i = 1,2$, maps $H$ into $H$, is monotone, i.e., for any $\varphi_1 \le \varphi_2$, we have $A_i \varphi_1 \le A_i \varphi_2$, and is non-expansive (and hence continuous with respect to the metric $\rho$):
    \begin{equation*}
        \rho(A_i \varphi_1, A_i \varphi_2) \le \rho(\varphi_1, \varphi_2) \quad \forall \varphi_1, \varphi_2 \in H.
    \end{equation*}
%\end{lemma}

\begin{proof}
    First, we prove the assertion for the operator $A_1$. The functions $\phi_u^\pm(z) = \varphi(u \pm z)$ are monotone and bounded, so the integrals in the operator are defined. For any $\varphi(u)$, $0 \le \varphi(u) \le 1$, obviously, $A_1 \varphi(u) \ge 0$ and
    \begin{align*}
        A_1 \varphi(u) \le A_1 1 & \le \frac{\alpha_c}{\alpha_c + \alpha_p} \int_0^u dF(z) + \frac{\alpha_p}{\alpha_c + \alpha_p} \int_0^{+\infty} dG(z) \\
        & \le \frac{\alpha_c}{\alpha_c + \alpha_p} F(u) + \frac{\alpha_p}{\alpha_c + \alpha_p} G(+\infty) \le 1.
    \end{align*}
    The non-decrease of the function $A_1 \varphi(u)$ in $u$ follows from the monotonicity and non-negativity of $\varphi(\cdot)$. Thus, $A_1: H \to H$. The monotonicity of the integral operator $A_1$ follows from its linearity and the non-negativity of the kernel. Finally, for any $\varphi_1, \varphi_2$:
    \begin{align*}
        & |A_1 \varphi_1(u) - A_1 \varphi_2(u)| \\
        & \le \frac{\alpha_c}{\alpha_c + \alpha_p} \int_0^u \sup_{x \in [0,u]} |\varphi_1(x) - \varphi_2(x)| dF(z) \\
        & \quad + \frac{\alpha_p}{\alpha_c + \alpha_p} \int_0^{+\infty} \sup_{x \in [u,+\infty)} |\varphi_1(x) - \varphi_2(x)| dG(z) \\
        & \le \frac{\alpha_c}{\alpha_c + \alpha_p} \sup_{x \in [0,u]} |\varphi_1(x) - \varphi_2(x)| F(u) \\
        & \quad + \frac{\alpha_p}{\alpha_c + \alpha_p} \sup_{x \in [u,+\infty)} |\varphi_1(x) - \varphi_2(x)| G(+\infty) \\
        & \le \frac{\alpha_c}{\alpha_c + \alpha_p} \rho(\varphi_1, \varphi_2) + \frac{\alpha_p}{\alpha_c + \alpha_p} \rho(\varphi_1, \varphi_2) = \rho(\varphi_1, \varphi_2).
    \end{align*}
    Hence, $\rho(A_1 \varphi_1, A_1 \varphi_2) \le \rho(\varphi_1, \varphi_2)$, and thus the operator $A_1$ is non-expansive.

    The proof for the operator $A_2$. For any function $\varphi(u)$, $0 \le \varphi(u) \le 1$, obviously, $A_2 \varphi(u) \ge 0$ and
    \begin{align*}
        A_2 \varphi(u) \le A_2 1
        & = \int_0^\infty \int_0^{U(u,t)} dF(z)(1 - K_p(t)) dK_c(t) \\
        & \quad + \int_0^\infty \int_0^{+\infty} dG(z)(1 - K_c(t)) dK_p(t) \\
        & = \int_0^\infty F(U(u,t))(1 - K_p(t)) dK_c(t) \\
        & \quad + \int_0^\infty (1 - K_c(t)) dK_p(t) \\
        & \le \int_0^\infty (1 - K_p(t)) dK_c(t) + \int_0^\infty (1 - K_c(t)) dK_p(t) = 1.
    \end{align*}
    The non-decrease of the function $A_2 \varphi(u)$ in $u$ follows from the monotonicity and non-negativity of $\varphi(\cdot)$ and the monotonicity of $U(\cdot, t)$. Thus, $A_2: H \to H$. The monotonicity of the integral operator $A_2$ follows from its linearity and the non-negativity of the kernel. Finally, for any $\varphi_1, \varphi_2$:
    \begin{align*}
        & |A_2 \varphi_1(u) - A_2 \varphi_2(u)| \\
        & \le \int_0^\infty \int_0^{U(u,t)} |\varphi_1(U(u,t)-z) - \varphi_2(U(u,t)-z)| dF(z)(1 - K_p(t)) dK_c(t) \\
        & \quad + \int_0^\infty \int_0^{+\infty} |\varphi_1(U(u,t)+z) - \varphi_2(U(u,t)+z)| dG(z)(1 - K_c(t)) dK_p(t) \\
        & \le \sup_{u \in [0,+\infty)} |\varphi_1(u) - \varphi_2(u)| \\
        & \quad \times \left[ \int_0^\infty \int_0^{U(u,t)} dF(z)(1 - K_p(t)) dK_c(t) \right. \\
        & \quad + \left. \int_0^\infty \int_0^{+\infty} dG(z)(1 - K_c(t)) dK_p(t) \right] \\
        & \le \sup_{u \in [0,+\infty)} |\varphi_1(u) - \varphi_2(u)|.
    \end{align*}
    Hence, $\rho(A_2 \varphi_1, A_2 \varphi_2) \le \rho(\varphi_1, \varphi_2)$, and thus the operator $A_2$ is non-expansive.
\end{proof}

The following assumptions are the main ones for solving problem (\ref{(6.7)}), (\ref{(6.8)}).

%\begin{assumption}
{\bf Assumption 6.2.}
    For the operator $A$, there exist monotone functions $\varphi_*(u), \varphi^*(u) \in H$ such that
    \begin{enumerate}
        \item[(a)] $\varphi_*(u) \le \varphi^*(u)$ and $A \varphi^*(u) \le \varphi^*(u), \quad A \varphi_*(u) \ge \varphi_*(u)$.
        \item[(b)] $\lim_{x \to +\infty} \varphi_*(u) = \lim_{x \to +\infty} \varphi^*(u) = 1$.
    \end{enumerate}
%\end{assumption}

In view of Lemma 6.2, we can take $\varphi^*(u) \equiv 1$. Another possible function $\varphi^*(u)$ for the operator $A_1$ is given in the following Lemma 3.

%\begin{lemma}
    For any monotone non-decreasing solution of problems (\ref{(6.3)}), (\ref{(6.4)}), (\ref{(6.8)}), from equation (\ref{(6.3)}), (\ref{(6.4)}), the estimate follows
    \begin{equation*}
        \varphi(u) \le \varphi_1^*(u) = \frac{\alpha_p}{\alpha_p + \alpha_c(1 - F(u))}.
    \end{equation*}
%\end{lemma}

%\begin{lemma}
{\bf Lemma 6.3.}
    For the operator $A_1$, the relation $A_1 \varphi_1^*(u) \le \varphi_1^*(u)$ holds.
%\end{lemma}

\begin{proof}
    Obviously, $\varphi_1^*(u)$ is monotonically non-decreasing, $\frac{\alpha_p}{\alpha_c + \alpha_p} \le \varphi_1^*(u) \le 1$ and $\lim_{u \to +\infty} \varphi_1^*(u) = 1$. Then the assertion of the lemma follows from the estimates:
    \begin{align*}
        A_1 \varphi_1^*(u)
        & \le \frac{\alpha_c}{\alpha_c + \alpha_p} \int_0^u \frac{\alpha_p dF(z)}{\alpha_p + \alpha_c(1 - F(u-z))} \\
        & \quad + \frac{\alpha_p}{\alpha_c + \alpha_p} \int_0^{+\infty} \frac{\alpha_p dG(z)}{\alpha_p + \alpha_c(1 - F(u+z))} \\
        & \le \frac{\alpha_c}{(\alpha_c + \alpha_p)} \times \frac{\alpha_p}{\alpha_p + \alpha_c(1 - F(u))} \times F(u) + \frac{\alpha_p}{\alpha_c + \alpha_p} = \varphi_1^*(u).
    \end{align*}
\end{proof}

The following two lemmas give sufficient conditions for Assumption 6.1(a) to hold, concerning the existence of the function $\varphi_*(u)$.

%\begin{lemma}
{\bf Lemma 6.4.}
    Let there exist a positive solution $L_1$ of the inequality
    \begin{equation*}
        \frac{\alpha_c}{\alpha_c + \alpha_p} {\mathbb E} e^{L_1 z_i} + \frac{\alpha_p}{\alpha_c + \alpha_p} {\mathbb E} e^{-L_1 p_i}
        = \frac{\alpha_c}{\alpha_c + \alpha_p} \int_0^{+\infty} e^{L_1 z} dF(z)
        + \frac{\alpha_p}{\alpha_c + \alpha_p} \int_0^{+\infty} e^{-L_1 z} dG(z) \le 1
    \end{equation*}
    and $\lim_{z \to \infty} e^{L_1 z}(1 - F(z)) = 0$. Let $\varphi_*(u) = 1 - e^{-L_1 u}$. Then $\varphi_*(u)$ is a non-decreasing function such that $0 \le \varphi_*(u) \le 1$, $\lim_{u \to +\infty} \varphi_*(u) = 1$, and for the operator $A_1$ the relation
$A_1 \varphi_*(u) \ge \varphi_*(u)$    holds.
%\end{lemma}

\begin{proof}
    The assertion of the lemma follows from the estimates
    \begin{align*}
        A_1 \varphi_*(u)
        & = \frac{\alpha_c}{\alpha_c + \alpha_p} \int_0^u (1 - e^{-L_1(u-z)}) dF(z) 
         + \frac{\alpha_p}{\alpha_c + \alpha_p} \int_0^{+\infty} (1 - e^{-L_1(u+z)}) dG(z) \\
        & = \frac{\alpha_c}{\alpha_c + \alpha_p} \left( F(u) - e^{-L_1 u} \int_0^u e^{L_1 z} dF(z) \right) \\
        & \quad + \frac{\alpha_p}{\alpha_c + \alpha_p} \left( 1 - e^{-L_1 u} \int_0^\infty e^{-L_1 z} dG(z) \right) \\
        & = \frac{\alpha_c}{\alpha_c + \alpha_p} F(u) + \frac{\alpha_p}{\alpha_c + \alpha_p} 
         + \frac{\alpha_c}{\alpha_c + \alpha_p} e^{-L_1 u} \int_u^{+\infty} e^{L_1 z} dF(z) \\
        & \quad - e^{-L_1 u} \left( \frac{\alpha_c}{\alpha_c + \alpha_p} \int_0^\infty e^{L_1 z} dF(z) + \frac{\alpha_p}{\alpha_c + \alpha_p} \int_0^\infty e^{-L_1 z} dG(z) \right) \\
        & \ge 1 - e^{-L_1 u} + \frac{\alpha_c}{\alpha_c + \alpha_p} \left( F(u) - 1 + e^{-L_1 u} \int_u^{+\infty} e^{L_1 z} dF(z) \right) \\
        & = 1 - e^{-L_1 u} + \frac{\alpha_c}{\alpha_c + \alpha_p} \left( F(u) - 1 - \left. e^{-L_1 u} e^{L_1 z} (1 - F(z)) \right|_u^{+\infty} \right. \\
        & \quad + \left. L_1 e^{-L_1 u} \int_u^{+\infty} (1 - F(z)) e^{L_1 z} dz \right) \\
        & = 1 - e^{-L_1 u} + \frac{\alpha_c}{\alpha_c + \alpha_p} L_1 e^{-L_1 u} \int_u^{+\infty} (1 - F(z)) e^{L_1 z} dz \ge 1 - e^{-L_1 u}.
    \end{align*}
\end{proof}

%\begin{lemma}
{\bf Lemma 6.5.}
    Let there exist constants $u_* \ge 0, c_* > 0, L > 0$ such that {\bf (a)} $U(u,t) \ge u + c_* t$ for all $u \ge u_*$ and {\bf (b)} there exists a positive solution $L_2$ of the inequality
    \begin{align}
       & \int_0^{+\infty} e^{L_2 z} dF(z) \cdot \int_0^\infty e^{-c_* L_2 t} (1 - K_p(t)) dK_c(t)\nonumber\\
       & \;\;\;\;\;\;\;\;\;\;+ \int_0^{+\infty} e^{-L_2 z} dG(z) \cdot \int_0^\infty e^{-c_* L_2 t} (1 - K_c(t)) dK_p(t) \le 1\nonumber
    \end{align}
    such that $\lim_{z \to \infty} e^{L_2 z}(1 - F(z)) = 0$. Let $\varphi_*(u) = \max\{0, 1 - e^{-L_2(u - u_*)}\}$. Then $\varphi_*(u)$ is a non-decreasing function such that $0 \le \varphi_*(u) \le 1$, $\lim_{u \to +\infty} \varphi_*(u) = 1$, and for the operator $A_2$ the relation
$A_2 \varphi_*(u) \ge \varphi_*(u)$  holds.
%\end{lemma}

\begin{proof}
    Let $u \ge u_*$. The assertion follows from the estimates
    \begin{align*}
        A_2 \varphi_*(u)
        & = \int_0^\infty \int_0^{U(u,t)} \max\{0, 1 - e^{-L_2(U(u,t) - u_* - z)}\} dF(z) (1 - K_p(t)) dK_c(t) \\
        & \quad + \int_0^\infty \int_0^{+\infty} \max\{0, 1 - e^{-L_2(U(u,t) - u_* + z)}\} dG(z) (1 - K_c(t)) dK_p(t) \\
        & = L_2 \int_0^\infty e^{-L_2(U(u,t) - u_*)} \int_0^{U(u,t) - u_*} e^{L_2 z} F(z) dz (1 - K_p(t)) dK_c(t) \\
        & \quad + \int_0^\infty \left( 1 - e^{-L_2(U(u,t) - u_*)} \int_0^{+\infty} e^{-L_2 z} dG(z) \right) (1 - K_c(t)) dK_p(t) \\
        & \ge \int_0^\infty e^{-L_2(U(u,t) - u_*)} \left( e^{L_2(U(u,t) - u_*)} - \int_0^{+\infty} e^{L_2 z} dF(z) \right) (1 - K_p(t)) dK_c(t) \\
        & \quad + \int_0^\infty \left( 1 - e^{-L_2(U(u,t) - u_*)} \int_0^{+\infty} e^{-L_2 z} dG(z) \right) (1 - K_c(t)) dK_p(t) \\
        & = \int_0^\infty \left( 1 - e^{-L_2(U(u,t) - u_*)} \int_0^{+\infty} e^{L_2 z} dF(z) \right) (1 - K_p(t)) dK_c(t) \\
        & \quad + \int_0^\infty \left( 1 - e^{-L_2(U(u,t) - u_*)} \int_0^{+\infty} e^{-L_2 z} dG(z) \right) (1 - K_c(t)) dK_p(t) \\
        & \ge 1 - e^{L_2 u_*} \left( \int_0^{+\infty} e^{L_2 z} dF(z) \cdot \int_0^\infty e^{-L_2 U(u,t)} (1 - K_p(t)) dK_c(t) \right. \\
        & \quad + \left. \int_0^{+\infty} e^{-L_2 z} dG(z) \cdot \int_0^\infty e^{-L_2 U(u,t)} (1 - K_c(t)) dK_p(t) \right).
    \end{align*}
    For $u \ge u_*$, we have $U(u,t) \ge u + c_* t$, and thus,
    \begin{align*}
        A_2 \varphi_*(u)
        & \ge 1 - e^{-L_2(u - u_*)} \left( \int_0^{+\infty} e^{L_2 z} dF(z) \cdot \int_0^\infty e^{-c_* L_2 t} (1 - K_p(t)) dK_c(t) \right. \\
        & \quad + \left. \int_0^{+\infty} e^{-L_2 z} dG(z) \cdot \int_0^\infty e^{-c_* L_2 t} (1 - K_c(t)) dK_p(t) \right) \ge 1 - e^{-L_2(u - u_*)}.
    \end{align*}
    The lemma is proved.
\end{proof}

%\begin{corollary}
{\bf Corollary 6.1.}
    When $K_c(t) = 1 - e^{-\alpha_c t}$, $K_p(t) = 1 - e^{-\alpha_p t}$, the inequality for finding $L_2$ takes the form:
    \begin{equation*}
        \frac{\alpha_c}{\alpha_c + \alpha_p + c_* L_2} \int_0^{+\infty} e^{L_2 z} dF(z)
        + \frac{\alpha_p}{\alpha_c + \alpha_p + c_* L_2} \int_0^{+\infty} e^{-L_2 z} dG(z) \le 1.
    \end{equation*}
    If $F(\cdot)$ and $G(\cdot)$ are such that
    $\int_0^{+\infty} e^{L_2 z} dF(z) = 1 + \mu_c L_2 + o(L_2)$
    and
    $\int_0^{+\infty} e^{-L_2 z} dG(z) = 1 - \mu_p L_2 + o(L_2)$
    (for example, claims and premiums are bounded), where $\mu_c$ and $\mu_p$ are the mean claims and premiums, then a positive solution $L_2$ of the inequality exists under the condition $\alpha_c \mu_c < \alpha_p \mu_p + c_*$.
%\end{corollary}

Within the framework of Assumption 6.1, define the set $H^*\subset H$ of non-decreasing functions such that $\varphi_*(u) \le \varphi(u) \le \varphi^*(u)$. Obviously, $\varphi_*, \varphi^* \in H^*$.

Under additional assumptions, the operators $A_i$ are contractions on the set $H^*$, but the contraction property is not uniform on $H^*$.

%\begin{assumption}
{\bf Assumption 6.3.}
    $F(z) < 1 \ \forall z \ge 0$.
		
{\bf Assumption 6.4.}
    $G(z) < 1 \ \forall z \ge 0$.
%\end{assumption}

%\begin{lemma}
{\bf Lemma 6.6.}
    (On the contraction coefficient of the operators $A_i, i = 1,2$, when $F(\cdot) < 1$). Let Assumptions 6.2, 6.3 hold for $A_1$ and 6.1–6.3 hold for $A_2$. Then for any $\varepsilon > 0$, there exist numbers $q_i^*(\varepsilon)$, $0 \le q_i^*(\varepsilon) < 1$, $i = 1,2$, such that for any functions $\varphi_1, \varphi_2 \in H^*$ with distance $\rho(A \varphi_1, A \varphi_2) \ge \varepsilon$, the following holds:
    \begin{equation*}
        \rho(A_i \varphi_1, A_i \varphi_2) \le q_i^*(\varepsilon) \cdot \rho(\varphi_1, \varphi_2), \quad i = 1,2.
    \end{equation*}
%\end{lemma}

\begin{proof}
    Fix $\varepsilon > 0$. Find a number $u^*(\varepsilon) \ge 0$ such that $\varphi^*(u) - \varphi_*(u) \le \varepsilon/2$ for all $u \ge u^*(\varepsilon)$. Let the functions $\varphi_1, \varphi_2 \in H^*$ be such that $\rho(A_i \varphi_1, A_i \varphi_2) \ge \varepsilon, i = 1,2$. Obviously, there exists a sequence $\{u^s\}$ such that
    \begin{equation*}
        \lim_{s \to +\infty} |A_i \varphi_1(u^s) - A_i \varphi_2(u^s)| = \rho(A_i \varphi_1, A_i \varphi_2) \ge \varepsilon > 0.
    \end{equation*}
    Since $\varphi_*(u^s) \le A_i \varphi_1(u^s) \le \varphi^*(u^s)$ and $\varphi_*(u^s) \le A_i \varphi_2(u^s) \le \varphi^*(u^s)$, we have
    \begin{equation*}
        |A_i \varphi_1(u^s) - A_i \varphi_2(u^s)| \le \varphi^*(u^s) - \varphi_*(u^s).
    \end{equation*}
    Hence, $u^s \le u^*(\varepsilon)$ for all sufficiently large $s$. Without loss of generality, we can assume that $\lim_{s \to \infty} u^s = u^*$. For the operator $A_1$, the following estimate holds:
    \begin{align*}
        & |A_1 \varphi_1(u^s) - A_1 \varphi_2(u^s)| \\
        & \le \frac{\alpha_c}{\alpha_c + \alpha_p} \int_0^{u^s} |\varphi_1(u^s - z) - \varphi_2(u^s - z)| dF(z) \\
        & \quad + \frac{\alpha_p}{\alpha_c + \alpha_p} \int_0^{+\infty} |\varphi_1(u^s + z) - \varphi_2(u^s + z)| dG(z) \\
        & \le \rho(\varphi_1, \varphi_2) \left( \frac{\alpha_c}{\alpha_c + \alpha_p} F(u^s) + \frac{\alpha_p}{\alpha_c + \alpha_p} \right) \\
        & \le \rho(\varphi_1, \varphi_2) \left( \frac{\alpha_c}{\alpha_c + \alpha_p} F(u^*(\varepsilon)) + \frac{\alpha_p}{\alpha_c + \alpha_p} \right).
    \end{align*}
    Denote
    \begin{equation*}
        q_1^*(\varepsilon) = \left( \frac{\alpha_c}{\alpha_c + \alpha_p} F(u^*(\varepsilon)) + \frac{\alpha_p}{\alpha_c + \alpha_p} \right).
    \end{equation*}
    Since $F(\cdot) < 1$, we have $q_1^*(\varepsilon) < 1$. Passing to the limit in $s$, we obtain the assertion of the lemma for the operator $A_1$.

    For the operator $A_2$, the following estimate holds:
    \begin{align*}
        & |A_2 \varphi_1(u^s) - A_2 \varphi_2(u^s)| \\
        & \le \sup_{u \ge 0} |\varphi_1(u) - \varphi_2(u)| \\
        & \quad \times \left[ \int_0^\infty \int_0^{U(u^s,t)} dF(z)(1 - K_p(t)) dK_c(t) + \int_0^\infty \int_0^\infty dG(z)(1 - K_c(t)) dK_p(t) \right] \\
        & = \rho(\varphi_1, \varphi_2) \left[ \int_0^\infty F(U(u^s,t))(1 - K_p(t)) dK_c(t) + \int_0^\infty (1 - K_c(t)) dK_p(t) \right] \\
        & \le \rho(\varphi_1, \varphi_2) q_2^*(\varepsilon),
    \end{align*}
    where
    \begin{equation*}
        q_2^*(\varepsilon) = \left[ \int_0^\infty F(U(u^*(\varepsilon), t))(1 - K_p(t)) dK_c(t) + \int_0^\infty (1 - K_c(t)) dK_p(t) \right].
    \end{equation*}
    Since $F(\cdot) < 1$, taking into account Assumption 6.1 and the fact that
    \begin{equation*}
        \int_0^\infty (1 - K_p(t)) dK_c(t) + \int_0^\infty (1 - K_c(t)) dK_p(t) = 1,
    \end{equation*}
    we have $q_2^*(\varepsilon) < 1$. Passing to the limit in $s$, we obtain the assertion of the lemma for $A_2$.
\end{proof}

%\begin{lemma}
{\bf Lemma 6.7.}
    (On the non-uniform contraction property of the operators $A_i, i = 1,2$, when $G(\cdot) < 1$). Let Assumptions 6.1, 6.2, 6.4 be satisfied. Then for any $\varepsilon > 0$, there exist numbers $q_i^*(\varepsilon)$, $0 \le q_i^*(\varepsilon) < 1$, $i = 1,2$, such that for any functions $\varphi_1, \varphi_2 \in H^*$ with distance $\rho(A \varphi_1, A \varphi_2) \ge \varepsilon$, the following holds:
    \begin{equation*}
        \rho(A_i \varphi_1, A_i \varphi_2) \le q_i^*(\varepsilon) \cdot \rho(\varphi_1, \varphi_2), \quad i = 1,2.
    \end{equation*}
%\end{lemma}

\begin{proof}
    By Lemma 6.2, $A_i$ is non-expansive; hence, $\varepsilon \le \rho(A_i \varphi_1, A_i \varphi_2) \le \rho(\varphi_1, \varphi_2)$. Define
    \begin{equation*}
        u^*(\varepsilon) = \sup_{v \ge 0} \{ v : \varphi^*(v) - \varphi_*(v) \le \varepsilon/2 \} < +\infty.
    \end{equation*}
    By the definition of the distance, for $\varphi_1, \varphi_2 \in H^*$, there exists a sequence $\{u^s\}$ such that
    \begin{equation}
        \lim_{s \to +\infty} |A_i \varphi_1(u^s) - A_i \varphi_2(u^s)| = \rho(A_i \varphi_1, A_i \varphi_2) \ge \varepsilon > 0.\label{(6.9)}
    \end{equation}
    By Assumption 6.2(a), $\varphi_*(u) \le A_i \varphi_1(u) \le \varphi^*(u), \varphi_*(u) \le A_i \varphi_2(u) \le \varphi^*(u)$; therefore,
    \begin{equation}\label{(6.10)}
        |A_i \varphi_1(u^s) - A_i \varphi_2(u^s)| \le \varphi^*(u^s) - \varphi_*(u^s).
    \end{equation}
    From (\ref{(6.9)}) and (\ref{(6.10)}), it follows that the sequence $\{u^s\}$ is bounded, and $u^s \le u^*(\varepsilon) < +\infty$ for all sufficiently large $s$. Without loss of generality, we can assume that $\lim_{s \to \infty} u^s = u^*$. From the conditions $\varphi_*(u) \le \varphi_1(u) \le \varphi^*(u), \varphi_*(u) \le \varphi_2(u) \le \varphi^*(u)$, it also follows that $|\varphi_1(u) - \varphi_2(u)| \le \varphi^*(u) - \varphi_*(u)$. Then
  \begin{equation}\label{(6.11)}
\begin{aligned}
    & |A_1 \varphi_1(u^s) - A_1 \varphi_2(u^s)| \\
    & \le \frac{\alpha_c}{\alpha_c + \alpha_p} \int_0^{u^s} |\varphi_1(u^s - z) - \varphi_2(u^s - z)| dF(z) \\
    & \quad + \frac{\alpha_p}{\alpha_c + \alpha_p} \int_0^{+\infty} |\varphi_1(u^s + z) - \varphi_2(u^s + z)| dG(z) \\
    & \le \frac{\alpha_c}{\alpha_c + \alpha_p} \rho(\varphi_1, \varphi_2) \\
    & \quad + \frac{\alpha_p}{\alpha_c + \alpha_p} \int_0^{+\infty} \min\{\rho(\varphi_1, \varphi_2), \varphi^*(u^s + z) - \varphi_*(u^s + z)\} dG(z) \\
    & \le \frac{\alpha_c}{\alpha_c + \alpha_p} \rho(\varphi_1, \varphi_2) \\
    & \quad + \frac{\alpha_p}{\alpha_c + \alpha_p} \int_0^{+\infty} \min\{\rho(\varphi_1, \varphi_2), \varphi^*(u^*(\varepsilon) + z) - \varphi_*(z)\} dG(z) \\
    & \le \frac{\alpha_c}{\alpha_c + \alpha_p} \rho(\varphi_1, \varphi_2) \\
    & \quad + \frac{\alpha_p}{\alpha_c + \alpha_p} \rho(\varphi_1, \varphi_2) \\
    & \quad \times \int_0^{+\infty} \min\left\{1, \frac{\varphi^*(u^*(\varepsilon) + z) - \varphi_*(z)}{\rho(\varphi_1, \varphi_2)}\right\} dG(z) \\
    & \le \rho(\varphi_1, \varphi_2) \cdot \left( \frac{\alpha_c}{\alpha_c + \alpha_p} \right. \\
    & \quad + \left. \frac{\alpha_p}{\alpha_c + \alpha_p} \int_0^{+\infty} \min\left\{1, \frac{\varphi^*(u^*(\varepsilon) + z) - \varphi_*(z)}{\varepsilon}\right\} dG(z) \right).
\end{aligned}
\end{equation}    Passing to the limit in $s$, we get
    \begin{equation*}
        \rho(A_1 \varphi_1, A_1 \varphi_2)
        \le \rho(\varphi_1, \varphi_2) \cdot \left( \frac{\alpha_c}{\alpha_c + \alpha_p} + \frac{\alpha_p}{\alpha_c + \alpha_p} \int_0^{+\infty} \min\left\{1, \frac{\varphi^*(u^*(\varepsilon) + z) - \varphi_*(z)}{\varepsilon}\right\} dG(z) \right).
    \end{equation*}
    Since $G(\cdot) < 1$, taking into account Assumption 6.2(b), we get
    \begin{equation*}
        q_1^*(\varepsilon) = \frac{\alpha_c}{\alpha_c + \alpha_p} + \frac{\alpha_p}{\alpha_c + \alpha_p} \int_0^{+\infty} \min\left\{1, \frac{\varphi^*(u^*(\varepsilon) + z) - \varphi_*(z)}{\varepsilon}\right\} dG(z) < 1.
    \end{equation*}
    Hence, $\rho(A_1 \varphi_1, A_1 \varphi_2) \le q_1^*(\varepsilon) \cdot \rho(\varphi_1, \varphi_2)$.

    For the operator $A_2$, the following estimates hold:
		\begin{align*}
        & |A_2 \varphi_1(u^s) - A_2 \varphi_2(u^s)| \\
        & \le \rho(\varphi_1, \varphi_2) \int_0^\infty F(U(u^s,t))(1 - K_p(t)) dK_c(t) \\
				& \quad + \int_0^\infty \int_0^\infty \min\{\rho(\varphi_1, \varphi_2), \varphi^*(U(u^s,t) + z) - \varphi_*(U(u^s,t) + z)\} dG(z)(1 - K_c(t)) dK_p(t) 
				\end{align*}
    \begin{align*}
        & \le \rho(\varphi_1, \varphi_2) \int_0^\infty (1 - K_p(t)) dK_c(t) \\
        & \quad + \int_0^\infty \int_0^\infty \min\{\rho(\varphi_1, \varphi_2), \varphi^*(U(u^*(\varepsilon), t) + z) - \varphi_*(z)\} dG(z)(1 - K_c(t)) dK_p(t) \\
        & \le \rho(\varphi_1, \varphi_2) \cdot \left( \int_0^\infty (1 - K_p(t)) dK_c(t) \right. \\
        & \quad + \left. \int_0^\infty \int_0^\infty \min\left\{1, \frac{\varphi^*(U(u^*(\varepsilon), t) + z) - \varphi_*(z)}{\varepsilon}\right\} dG(z)(1 - K_c(t)) dK_p(t) \right).
    \end{align*}
    Passing to the limit in $s$, we get
    \begin{align*}
        \rho(A_2 \varphi_1, A_2 \varphi_2)
        & \le \rho(\varphi_1, \varphi_2) \cdot \left( \int_0^\infty (1 - K_p(t)) dK_c(t) \right. \\
        & \quad + \left. \int_0^\infty \int_0^\infty \min\left\{1, \frac{\varphi^*(U(u^*(\varepsilon), t) + z) - \varphi_*(z)}{\varepsilon}\right\} dG(z)(1 - K_c(t)) dK_p(t) \right).
    \end{align*}
    Since $0 \le G(\cdot) < 1, \int_0^{+\infty} dG(z) = 1$ and $\varphi^*(U(u^*(\varepsilon), t) + z) - \varphi_*(z) \to 0$ as $z \to +\infty$, the function
    \begin{equation*}
        f_\varepsilon^*(t) = \int_0^\infty \min\left\{1, \frac{\varphi^*(U(u^*(\varepsilon), t) + z) - \varphi_*(z)}{\varepsilon}\right\} dG(z) < 1
    \end{equation*}
    for all $t$. Hence, taking into account Assumption 6.1, we have
    \begin{equation*}
        \int_0^\infty f_\varepsilon^*(t)(1 - K_c(t)) dK_p(t) < \int_0^{+\infty} (1 - K_c(t)) dK_p(t),
    \end{equation*}
    and consequently,
    \begin{equation*}
        q_2^*(\varepsilon) = \int_0^{+\infty} (1 - K_p(t)) dK_c(t) + \int_0^{+\infty} f_\varepsilon^*(t)(1 - K_c(t)) dK_p(t) < 1.
    \end{equation*}
    The lemma is proved.
\end{proof}

\section{Necessary and Sufficient Conditions for the Existence and Uniqueness of a Solution to the Problem}

%\begin{theorem}
{\bf Theorem 6.3.}
    (Necessary and sufficient conditions for the existence and uniqueness of a solution to problems (\ref{(6.7)}), (\ref{(6.8)}). For the existence and uniqueness of a solution $\varphi \in H$ of the operator equation (\ref{(6.7)}) satisfying the boundary condition (\ref{(6.8)}) with operator $A \in \{A_1, A_2\}$, it is necessary and sufficient that Assumption 6.2 be satisfied.
%\end{theorem}

\begin{proof}
    The necessity is obvious; as $\varphi_*(u), \varphi^*(u)$, we can take any solution of problems (\ref{(6.7)}), (\ref{(6.8)}). Let us prove the sufficiency of the conditions for the existence of a solution by constructing a sequence of functions converging to some solution of the problem. Namely, consider the sequence of approximations
    \begin{equation*}
        \{\varphi^{k+1}(u) = A \varphi^k(u), \quad \varphi^0(u) \equiv \varphi^*(u), \quad k = 0,1,\ldots\}.
    \end{equation*}
    By the monotonicity of $\varphi^*(u)$, all functions $\varphi^k(u)$ are non-decreasing in $u$. By the monotonicity of the operator $A$ and the assumption $A \varphi^*(u) \le \varphi^*(u)$, it follows that the sequence $\{\varphi^k(u), k = 0,1,\ldots\}$ decreases monotonically. Since $\varphi^*(u) \ge \varphi_*(u)$, by the monotonicity of the operator $A$ and Assumption 6.2(a), we have $\varphi^1(u) = A \varphi^*(u) \ge A \varphi_*(u) \ge \varphi_*(u)$. Similarly, by induction, we get $\varphi^k(u) \ge \varphi_*(u)$. Thus, the sequence of functions $\{\varphi^k(u)\}$ decreases monotonically and is bounded below by the function $\varphi_*(u)$. Therefore, there exists a limit function $\varphi(u) = \lim_{k \to +\infty} \varphi^k(u)$, which, like all $\varphi^k(u)$, is non-decreasing in $u$, $\varphi^*(u) \ge \varphi(u) \ge \varphi_*(u)$, and thus, by Assumption 6.2(b), $\lim_{u \to +\infty} \varphi(u) = 1$. Passing to the limit in $k$ in the relation $\varphi^{k+1}(u) = A \varphi^k(u)$, by the continuity of the operator $A$ with respect to pointwise convergence, the limit function $\varphi(u)$ satisfies the equation $\varphi = A \varphi$.

    Let us prove the uniqueness of the solution. Let $\varphi_1, \varphi_2 \in H^*$ be two different solutions of equation (\ref{(6.7)}). By the definition of the distance, for $\varphi_1, \varphi_2 \in H^*$, there exists a sequence $\{u^s\}$ such that
    \begin{equation}\label{(6.12)}
        \lim_{s \to +\infty} |\varphi_1(u^s) - \varphi_2(u^s)| = \rho(\varphi_1, \varphi_2) > 0.
    \end{equation}
    Note that
    \begin{equation}\label{(6.13)}
        \rho(A \varphi_1, A \varphi_2) = \lim_{s \to \infty} |A \varphi_1(u^s) - A \varphi_2(u^s)| = \lim_{s \to \infty} |\varphi_1(u^s) - \varphi_2(u^s)| = \rho(\varphi_1, \varphi_2).
    \end{equation}
    From the assumption $\varphi_1, \varphi_2 \in H^*$, it follows that $|\varphi_1(u) - \varphi_2(u)| \le \varphi^*(u) - \varphi_*(u)$ and thus $\lim_{u \to \infty} |\varphi_1(u) - \varphi_2(u)| = 0$. Hence, it follows that the sequence $\{u^s\}$ is bounded. Without loss of generality, we can assume that $\{u^s\}$ has a limit, $\lim_{s \to \infty} u^s = u^*$. Also, without loss of generality, we can assume that $u^*$ is the extreme right (maximal) limit point of the sequences $\{u^s\}$ satisfying condition (\ref{(6.12)}). This means that for any $\varepsilon > 0$, there exists $\delta(\varepsilon) > 0$ such that
    \begin{equation*}
        |\varphi_1(u) - \varphi_2(u)| \le \rho(\varphi_1, \varphi_2) - \delta(\varepsilon)
    \end{equation*}
    for all $u \ge u^* + \varepsilon$.

    Proof for $A_1$. Fix an arbitrary $\varepsilon > 0$. The following estimates hold:
    \begin{align*}
        & |A_1 \varphi_1(u^s) - A_1 \varphi_2(u^s)| \\
        & \le \frac{\alpha_c}{\alpha_c + \alpha_p} \int_0^{u^s} |\varphi_1(u^s - z) - \varphi_2(u^s - z)| dF(z) \\
        & \quad + \frac{\alpha_p}{\alpha_c + \alpha_p} \int_0^{+\infty} |\varphi_1(u^s + z) - \varphi_2(u^s + z)| dG(z) \\
        & \le \frac{\alpha_c}{\alpha_c + \alpha_p} \rho(\varphi_1, \varphi_2) \\
        & \quad + \frac{\alpha_p}{\alpha_c + \alpha_p} \int_0^{2\varepsilon} |\varphi_1(u^s + z) - \varphi_2(u^s + z)| dG(z) \\
        & \quad + \frac{\alpha_p}{\alpha_c + \alpha_p} \int_{2\varepsilon}^{+\infty} |\varphi_1(u^s + z) - \varphi_2(u^s + z)| dG(z) \\
        & \le \frac{\alpha_c}{\alpha_c + \alpha_p} \rho(\varphi_1, \varphi_2) + \frac{\alpha_p}{\alpha_c + \alpha_p} \rho(\varphi_1, \varphi_2) G(2\varepsilon) \\
        & \quad + \frac{\alpha_p}{\alpha_c + \alpha_p} \int_{2\varepsilon}^{+\infty} |\varphi_1(u^s + z) - \varphi_2(u^s + z)| dG(z).
    \end{align*}
    For a given $\varepsilon$, find $\delta(\varepsilon)$ such that
    \begin{equation*}
        |\varphi_1(u) - \varphi_2(u)| \le \rho(\varphi_1, \varphi_2) - \delta(\varepsilon)
    \end{equation*}
    for all $u \ge u^* + \varepsilon$. Since by construction $\lim_{s \to \infty} u^s = u^*$, for all sufficiently large $s$, we have $|u^* - u^s| \le \varepsilon$ and for $z \ge 2\varepsilon$, we have $u^s + z \ge \varepsilon$, $|\varphi_1(u^s + z) - \varphi_2(u^s + z)| \le \rho(\varphi_1, \varphi_2) - \delta(\varepsilon)$. Therefore, the following estimates hold:
    \begin{align*}
        & |A_1 \varphi_1(u^s) - A_1 \varphi_2(u^s)| \\
        & \le \frac{\alpha_c}{\alpha_c + \alpha_p} \rho(\varphi_1, \varphi_2) + \frac{\alpha_p}{\alpha_c + \alpha_p} \rho(\varphi_1, \varphi_2) G(2\varepsilon) \\
        & \quad + \frac{\alpha_p}{\alpha_c + \alpha_p} \int_{2\varepsilon}^{+\infty} (\rho(\varphi_1, \varphi_2) - \delta) dG(z) \\
        & \le \frac{\alpha_c}{\alpha_c + \alpha_p} \rho(\varphi_1, \varphi_2) + \frac{\alpha_p}{\alpha_c + \alpha_p} \rho(\varphi_1, \varphi_2) G(2\varepsilon) \\
        & \quad + \frac{\alpha_p}{\alpha_c + \alpha_p} (\rho(\varphi_1, \varphi_2) - \delta)(1 - G(2\varepsilon)) \\
        & \le \rho(\varphi_1, \varphi_2) - \frac{\alpha_p}{\alpha_c + \alpha_p} \delta(1 - G(2\varepsilon)).
    \end{align*}
    Passing to the limit here as $s \to \infty$, taking into account (\ref{(6.13)}), we get
    \begin{equation*}
        \rho(A \varphi_1, A \varphi_2) \le \rho(\varphi_1, \varphi_2) - \frac{\alpha_p}{\alpha_c + \alpha_p} \delta(\varepsilon) \cdot (1 - G(2\varepsilon)),
    \end{equation*}
    which was required to prove.

    The proof of the assertion of the lemma for the operator $A_2$ is carried out similarly. The following estimates hold:
    \begin{align*}
        & |A_2 \varphi_1(u^s) - A_2 \varphi_2(u^s)| \\
        & \le \int_0^\infty \int_0^{U(u^s,t)} |\varphi_1(U(u^s,t) - z) - \varphi_2(U(u^s,t) - z)| dF(z)(1 - K_p(t)) dK_c(t) \\
        & \quad + \int_0^\infty \int_0^{+\infty} |\varphi_1(U(u^s,t) + z) - \varphi_2(U(u^s,t) + z)| dG(z)(1 - K_c(t)) dK_p(t) \\
        & \le \rho(\varphi_1, \varphi_2) \int_0^\infty F(U(u^s,t))(1 - K_p(t)) dK_c(t) \\
        & \quad + \int_0^\infty \left( \int_0^{2\varepsilon} + \int_{2\varepsilon}^{\infty} \right) |\varphi_1(U(u^s,t) + z) - \varphi_2(U(u^s,t) + z)| dG(z)(1 - K_c(t)) dK_p(t) \\
        & \le \rho(\varphi_1, \varphi_2) \int_0^\infty (1 - K_p(t)) dK_c(t) + \rho(\varphi_1, \varphi_2) G(2\varepsilon) \int_0^\infty (1 - K_c(t)) dK_p(t) \\
        & \quad + \int_0^\infty \int_{2\varepsilon}^{\infty} |\varphi_1(U(u^s,t) + z) - \varphi_2(U(u^s,t) + z)| dG(z)(1 - K_c(t)) dK_p(t).
    \end{align*}
    Since $c(\cdot) \ge 0$, we have $U(u^s,t) \ge u^s$. For all sufficiently large $s$, we have $|u^* - u^s| \le \varepsilon$ and for $z \ge 2\varepsilon$, we have $U(u^s,t) + z \ge u^s + z \ge \varepsilon$, $|\varphi_1(U(u^s,t) + z) - \varphi_2(U(u^s,t) + z)| \le \rho(\varphi_1, \varphi_2) - \delta(\varepsilon)$. Thus, for sufficiently large $s$:
    \begin{align*}
        & |A_2 \varphi_1(u^s) - A_2 \varphi_2(u^s)| \\
        & \le \rho(\varphi_1, \varphi_2) \int_0^\infty (1 - K_p(t)) dK_c(t) \\
        & \quad + \rho(\varphi_1, \varphi_2) G(2\varepsilon) \int_0^\infty (1 - K_c(t)) dK_p(t) \\
        & \quad + \int_0^\infty \int_{2\varepsilon}^{\infty} (\rho(\varphi_1, \varphi_2) - \delta(\varepsilon)) dG(z)(1 - K_c(t)) dK_p(t) \\
        & \le \rho(\varphi_1, \varphi_2) - \delta(\varepsilon)(1 - G(2\varepsilon)) \int_0^\infty (1 - K_c(t)) dK_p(t).
    \end{align*}
    Hence, taking into account (\ref{(6.13)}), the assertion of the lemma follows for the operator $A_2$. The theorem is proved.
\end{proof}

\section{The Method of Successive Approximations for Solving the Problem}

Consider the method of successive approximations for solving problems (\ref{(6.7)}), (\ref{(6.8)}):
\begin{equation}\label{(6.14)}
    \varphi^{k+1}(u) := A \varphi^k(u), \quad k = 0,1,\ldots,
\end{equation}
where $k$ is the iteration number, $0 \le \varphi^0(u) \le 1$; $A \in \{A_1, A_2\}$.

%\begin{theorem}
{\bf Theorem 6.4.}
    (On the convergence of approximations from above and below). Under Assumptions 6.1 and 6.2, the sequence of approximations $\{\varphi^k(u)\}$ constructed according to (\ref{(6.14)}) and starting with $\varphi^0(u) \equiv \varphi^*(u)$ and, in particular, with $\varphi^0(u) \equiv 1$, decreases monotonically and converges pointwise from above to the solution of problems (\ref{(6.7)}), (\ref{(6.8)}), and the sequence of approximations $\{\varphi^k(u)\}$ starting with $\varphi^0(u) \equiv \varphi_*(u)$ increases monotonically and converges pointwise from below to the solution of problems (\ref{(6.7)}), (\ref{(6.8)}).
%\end{theorem}

The assertion of the theorem is a consequence of Lemmas 6.1–6.5 and Theorem 6.3.

%\begin{corollary}
{\bf Corollary 6.2.}
    Under Assumptions 6.1 and 6.2, for any initial approximation $\varphi^0 \in H^*$, the sequence $\{\varphi^k(u), k = 0,1,\ldots\}$ generated by algorithm (\ref{(6.14)}) converges pointwise to the solution of problems (\ref{(6.7)}), (\ref{(6.8)}).
%\end{corollary}

%\begin{theorem}
{\bf Theorem 6.5.}
    (On uniform convergence and the rate of convergence of the method of successive approximations). Under Assumptions 6.1, 6.3 (or 6.4), the method of successive approximations (\ref{(6.14)}) starting from $\varphi^0(u) \in H^*$ converges uniformly monotonically to the solution of problems (\ref{(6.7)}), (\ref{(6.8)}); moreover, it converges to any $\varepsilon$-neighborhood of the solution of problems (\ref{(6.7)}), (\ref{(6.8)}) at the rate of a geometric progression with a denominator depending on $\varepsilon$.
%\end{theorem}

The assertion of the theorem follows from Theorem 6.3 and Lemmas 6.6, 6.7.

\section{Numerical Experiments}

To illustrate the results, take the process (\ref{(6.1)}) in which all claims are deterministic and equal in size, $p_k \equiv p$, and all premiums are deterministic and equal, $z_k \equiv z$, but the moments of claim and premium arrivals are random (Poisson). If, in addition, $z = p = 1$, then an analytical solution of problems (\ref{(6.3)}), (\ref{(6.4)}), \ref{(6.8)} is known [5]:
\begin{equation*}
    \varphi(u) = \varphi([u]) = 1 - \left( \frac{\lambda_c}{\lambda_p} \right)^{[u] + 1}.
\end{equation*}
This same solution can also be found numerically by the method of successive approximations. In this case,
\begin{equation*}
    F(x) = \begin{cases}
        0, & x < z, \\
        1, & x \ge z;
    \end{cases}
    \quad
    G(x) = \begin{cases}
        0, & x < p, \\
        1, & x \ge p;
    \end{cases}
\end{equation*}
and the integral operator takes the form:
\begin{equation*}
    A\varphi(u) = \begin{cases}
        \frac{\alpha_p}{\alpha_c + \alpha_p} \varphi(u + p), & u < z; \\
        \frac{\alpha_p}{\alpha_c + \alpha_p} \varphi(u + p) + \frac{\alpha_c}{\alpha_c + \alpha_p} \varphi(u - z), & u \ge z;
    \end{cases}
\end{equation*}
Choose the process parameters $\lambda_c = 1, \lambda_p = 2, z = 1, p = 1$. Then the solution has the form $\varphi(u) = 1 - (0.5)^{[u] + 1}$ (bold dotted line in Fig. 6.1).

\begin{figure}[h]
    \centering
    \includegraphics[width=0.8\textwidth]{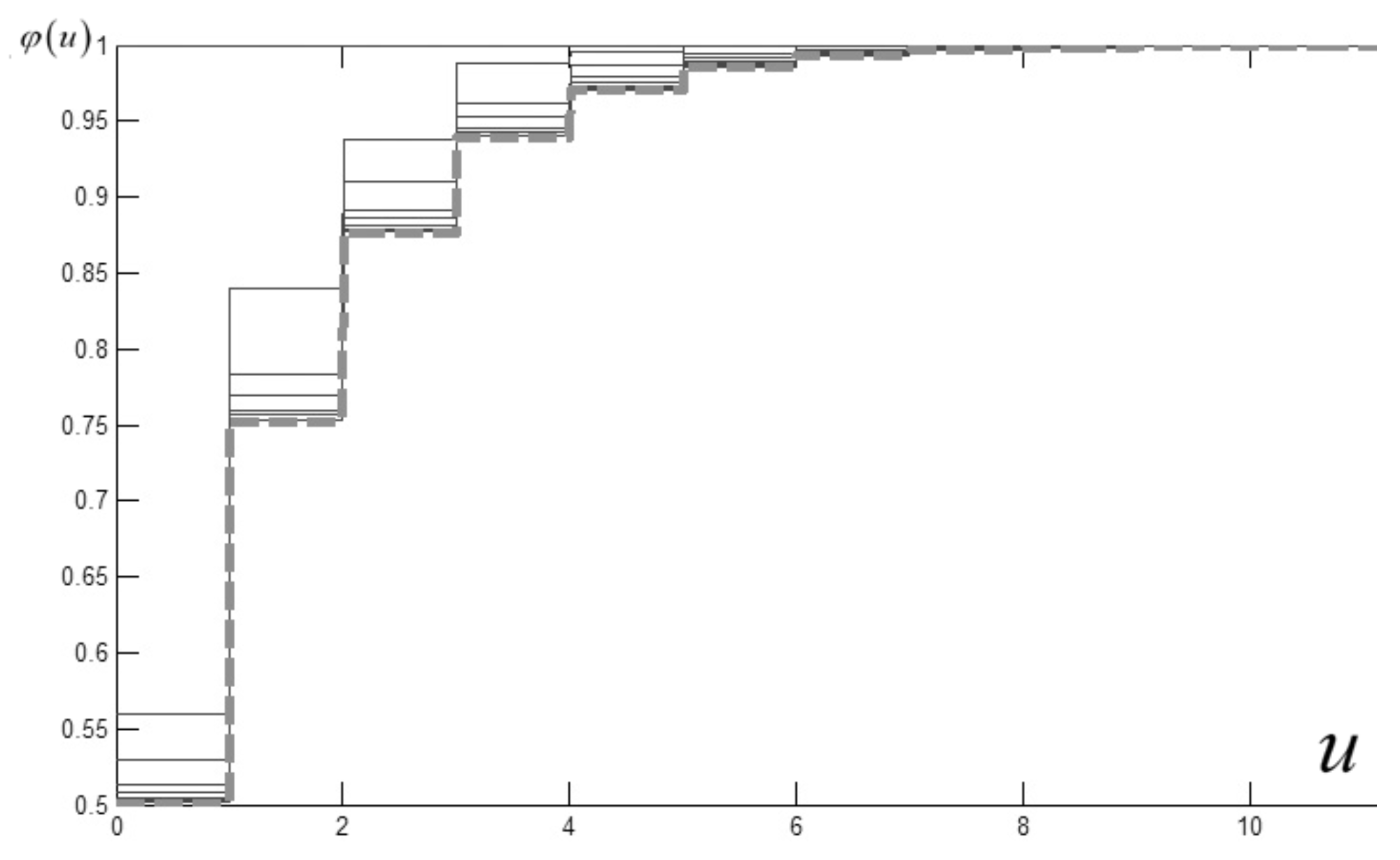}
    \caption{Fig. 6.1}
    %\label{fig:6.1}
\end{figure}

Fig. 6.1 shows the successive approximations (after 5 iterations) of the method to this solution, starting with $\varphi_0(u) \equiv 1$.

\begin{figure}[h]
    \centering
    \includegraphics[width=0.8\textwidth]{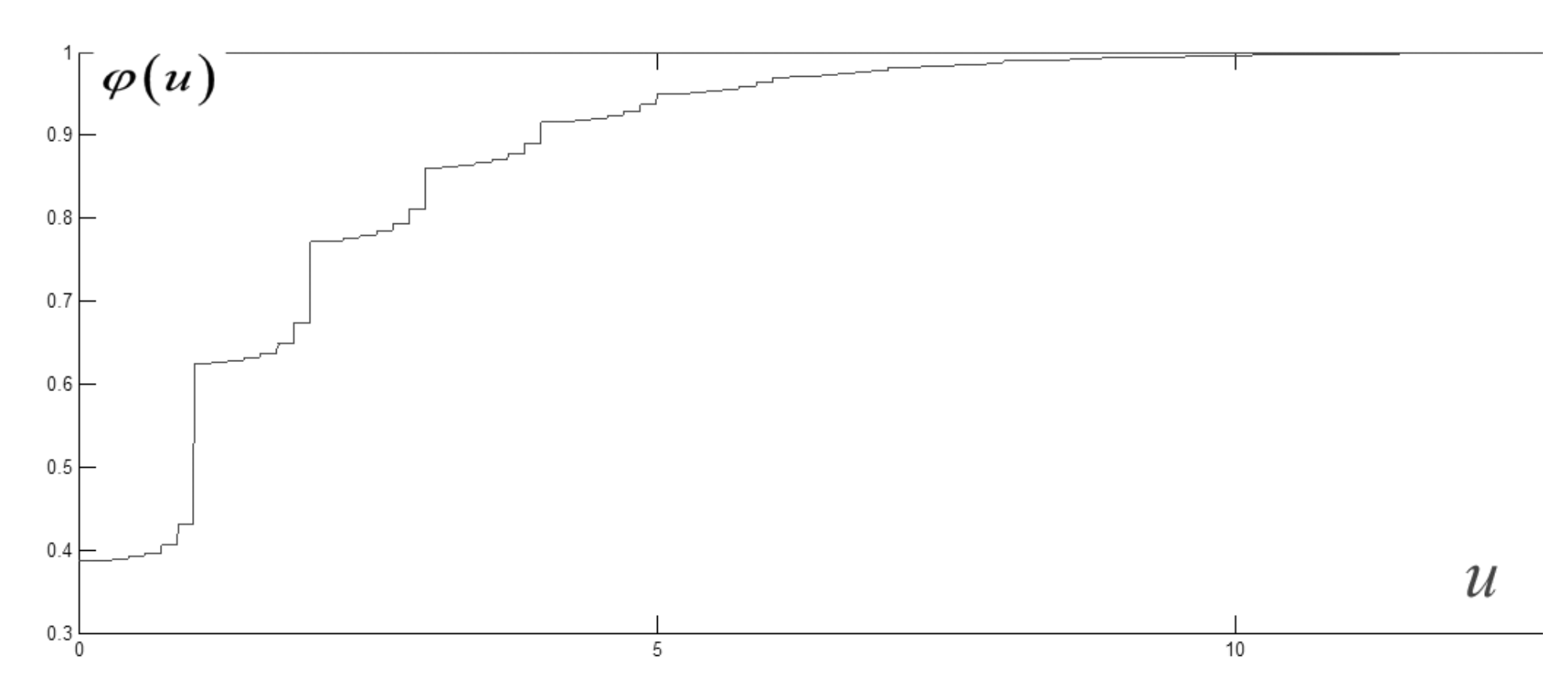}
    \caption{Fig. 6.2}
    %\label{fig:6.2}
\end{figure}

Fig. 6.2 shows the form of the solution obtained by the method of successive approximations after 50 iterations for $\lambda_c = \lambda_p = 1, z = 1, p = \pi - 1$. In this case, an analytical solution is not known.

\chapter{RISK PROCESS IN A MARKOVIAN ENVIRONMENT}

In this section, a general risk process is considered that describes the stochastic evolution of the capital of an insurance company in a Markovian random environment with nonlinear premium income and general insurance contracts. It is shown that the probabilities of ruin of the process (company) as functions of the initial state generally satisfy a certain system of integral equations with boundary conditions at infinity. Sufficient conditions for the existence of a solution to this system are established and the method of successive approximations for calculating the ruin probabilities is justified. When starting from the unit initial approximation, the iterations converge to the solution from above, and when starting from the generalized Cramér-Lundberg exponential bound, from below, and thus, it is always possible to estimate the accuracy of the obtained approximate solution. The performance of the method is verified on a numerical example of a risk process on a Markov chain with two states. The results of the section are published in [40, 43, 105].

\section{Systems of integral and integro-differential equations for the ruin probabilities of a risk process in a Markovian environment}

Consider a random risk process $\xi_t$ in a random (Markovian) environment that can be in one of $i=1,\ldots,n$ states. Then the probability of non-ruin of the process depends on the initial state of the process $(u,i)$, and thus, we have a set of probabilities
\begin{equation*}
    \{\varphi_i(u)=P[\inf_{0\le t<\infty}\xi_t\ge 0 | \xi_0=u, \text{ init. state of env.}=i]\}.
\end{equation*}

Let us derive a system of integral equations for these probabilities of non-ruin. A more general approach for obtaining the distribution of the absolute minimum of a process, as well as the corresponding probabilities of (non-)ruin, for general processes with independent increments on Markov chains is based on matrix integro-differential equations with partial derivatives for the matrix distribution function [18, p. 88]
$\{\varphi_{ij}(t,x)=P[\min_{0\le \tau \le t}\xi(\tau)>x | \xi(0)=0, \text{init. state}=i, \text{final state}=j]\}$.
However, these matrix equations are very complex. They contain partial derivatives with respect to time $t$, unlike similar equations for the classical risk process. To solve these equations, [18] develops analytical methods using matrix characteristic functions. Therefore, it makes sense to derive similar simpler equations directly for the probabilities
\begin{equation*}
    \{\varphi_i(u)=P[\inf_{0\le \tau <\infty}\xi(\tau)>0 | \xi(0)=u, \text{ init. state of chain}=i]\}
\end{equation*}
of non-ruin of the process on an infinite time interval with initial capital $u$. These equations no longer contain partial derivatives with respect to time $t$, unlike the equations in [18, (1.8)], and direct methods, such as the well-known Laplace transform, can be applied to solve them. The effectiveness of this approach is illustrated by an example for a chain with two states.

Consider a more general risk process $\xi(t)$ in a Markovian environment (on a Markov chain) with $m$ states $\{x_1,x_2,\ldots,x_m\}$ (Fig. 7.1):
\begin{align}
    \xi(t) & = u + \left( \sum_{h=0}^{M_t-1} c(h)(\tau_{h+1}-\tau_h) + c(M_t)(t-\tau_{M_t}) \right) \nonumber \\
    & \quad - \left( \sum_{h=0}^{M_t-1} \sum_{k=1}^{N(\tau_{h+1}-\tau_h,\alpha(h))} Y(h,k) + \sum_{k=1}^{N(t-\tau_{M_t},\alpha(M_t))} Y(M_t,k) \right),\label{(7.1)}
\end{align}
where $u$ is the initial capital of the company; $h$ is the ordinal number of the environment state switch; $M_t$ is the number of environment state switches by time $t$; $N(\Delta,\alpha)$ is the number of insurance events during time $\Delta$ with intensity $\alpha$; $\alpha_i$, $\alpha(h)$ are the intensities of insurance events in state $x_i$ and after the $h$-th switch, respectively; $c_i$ and $c(h)$ are the premium income rates in state $x_i$ and after the $h$-th switch, respectively; $\tau_h$ is the $h$-th moment of change of the chain state, $\tau_0=0$; $Y(k,h)$ is the $k$-th claim after the $h$-th switch.

The probability that the chain remains in state $i$ for time $\Delta$ is $e^{-\lambda_i\Delta}$, the probability of the chain transitioning from state $i$ to state $j\ne i$ during time $\Delta$ is $\lambda_i p_{ij}\Delta+o(\Delta)$, where $\sum_{j=1}^m p_{ij}=1$, $p_{ii}=0$.

\begin{figure}[h]
    \centering
    \includegraphics[width=0.8\textwidth]{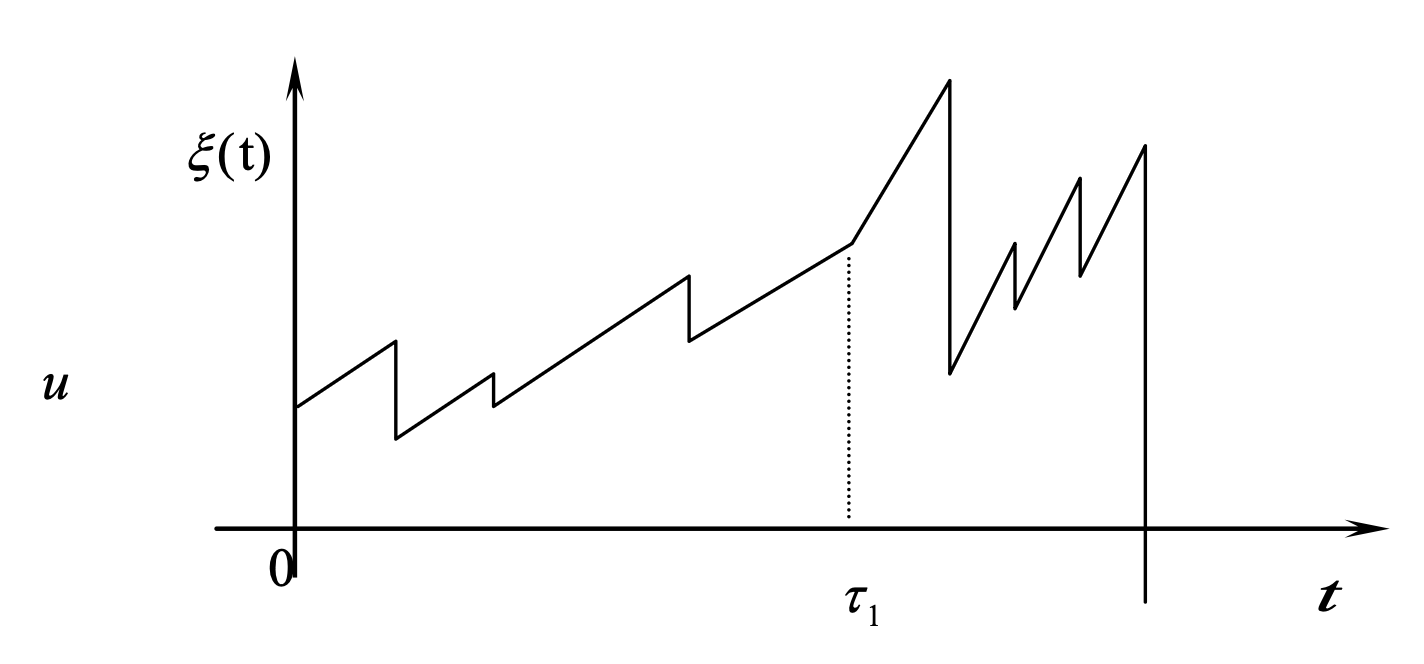}
    \caption{Trajectory of the risk process}
    \label{fig:7.1}
\end{figure}

Consider the probability $\psi_i$ of ruin from the $i$-th initial state of the chain and the probability that ruin does not occur when starting from the $i$-th state, $\varphi_i = 1 - \psi_i$.

We obtain integral equations for the probabilities of non-ruin without assuming differentiability of the functions $\varphi_i$.

%\begin{theorem}
{\bf Theorem 7.1} [40].
    The vector-function $\varphi(u) = \{\varphi_1(u), \varphi_2(u), \ldots, \varphi_m(u)\}$ satisfies the system of integral equations ($i = 1, \ldots, m$):
    \begin{equation}\label{(7.2)}
        \varphi_i(u) = \int_0^\infty e^{-t(\alpha_i + \lambda_i)} \left( \lambda_i \sum_{j=1}^m p_{ij} \varphi_j(u + c_i t) + \alpha_i \int_0^{u + c_i t} \varphi_i(u + c_i t - z) dF_i(z) \right) dt.
    \end{equation}
%\end{theorem}

\begin{proof}
    Let event $L$ be the non-ruin of the company.
    Divide the time axis into intervals of length $\Delta$ and represent the event of non-ruin from the initial state $(x_i, u)$ as the union over $t = 0, 1, \ldots$ of the following disjoint events $A(t), B(t), C(t)$ (see Figs. 7.2, 7.3).

    Event
    \begin{equation*}
        A(t) = A_1(t) \cap A_2(t) \cap \left( \bigcup_{z=0}^{u + a_i(t\Delta + O(\Delta))} (A_3(t,z) \cap A_4(t,z)) \right)
    \end{equation*}
    (see Fig. 7.2), where

    $A_1(t)$: on the interval $[0, (t+1)\Delta)$ there are no transitions (state changes) of the chain (probability of this event is $e^{-\lambda_i (t+1)\Delta}$);

    $A_2(t)$: on the interval $[0, t\Delta)$ no claims arrive (under $A_1(t)$, the conditional probability of this event is $e^{-\alpha_i t\Delta}$);

    $A_3(t,z)$: on the interval $[t\Delta, (t+1)\Delta)$ at the moment $(t\Delta + o(\Delta))$ one claim of size from the interval $[z, z+dz)$ arrives, $z < u + c_i(t\Delta + o(\Delta))$ (under $A_1(t)$, $A_2(t)$, the conditional probability of this event is $(\alpha_i \Delta + o(\Delta)) dF_i(z)$);

    $A_4(t,z)$: ruin does not occur from the new state $(x_i, u + c_i(t\Delta + o(\Delta)) - z)$ after the first claim (probability $\varphi_i(u + c_i(t\Delta + o(\Delta)) - z)$).

    Event
    \begin{equation*}
        B(t) = B_1(t) \cap B_2(t) \cap \left( \bigcup_{j=1, j \ne i}^m (B_3(t,j) \cap B_4(t,j)) \right)
    \end{equation*}
    (see Fig. 7.3), where

    $B_1(t)$: on the interval $[0, t\Delta)$ there are no transitions (state changes) of the chain (probability of this event is $e^{-\lambda_i t\Delta}$);

    $B_2(t)$: on the interval $[0, t\Delta)$ no claims arrive (under $B_1(t)$, the conditional probability of this event is $e^{-\alpha_i t\Delta}$);

    $B_3(t,j)$: on the interval $[t\Delta, (t+1)\Delta)$ at the moment $(t\Delta + o(\Delta))$ the chain transitions from the initial state $x_i$ to state $x_j$ (probability of this event is $(\lambda_i p_{ij} \Delta + o(\Delta))$);

    $B_4(t,j)$: after the transition, ruin does not occur (probability $\varphi_j(u + c_i(t\Delta + o(\Delta)))$).

    Event
    \begin{equation*}
        C(t) = C_1(t) \cap C_2(t) \cap C_3(t),
    \end{equation*}
    where

    $C_1(t)$: on the interval $[0, t\Delta)$ there are no transitions (state changes) of the chain (probability of this event is $e^{-\lambda_i t\Delta}$);

    $C_2(t)$: on the interval $[0, t\Delta)$ no claims arrive (under $C_1(t)$, the conditional probability of this event is $e^{-\alpha_i t\Delta}$);

    $C_3(t)$: on the interval $[t\Delta, (t+1)\Delta)$ more than one event occurs (a state change of the chain or a claim arrival).

    Let us estimate the probability of event $C_3(t)$:
    \begin{align*}
        P\{C_3(t)\} & = 1 - P\{\overline{C_3(t)}\} \\
        & \le 1 - (e^{-\lambda_i \Delta} e^{-\alpha_i \Delta} + e^{-\lambda_i \Delta} \alpha_i \Delta + o(\Delta) + \sum_{j} e^{-\lambda_i \Delta} e^{-\alpha_i \Delta} p_{ij} \lambda_i \Delta + o(\Delta)) = o(\Delta).
    \end{align*}

    Note that $A(n) \cap B(k) \cap C(m) = \varnothing$, $\forall n, k, m \in [0, +\infty)$; $A(n) \cap A(k) = \varnothing$, $B(n) \cap B(k) = \varnothing$, $C(n) \cap C(k) = \varnothing$, $\forall n \ne k$. Then
    \begin{equation*}
        L = \bigcup_t (A(t) \cup B(t) \cup C(t)).
    \end{equation*}

    Let us write the sum over time for the probability of non-ruin:
\begin{equation}\label{(7.3)}
\begin{split}
    \varphi_i(u)
    & = \sum_{t=0}^\infty \left[ e^{-\lambda_i (t+1)\Delta} e^{-\alpha_i t\Delta} \int_0^{u + c_i(t\Delta + o(\Delta))} (\alpha_i \Delta + o(\Delta)) \varphi_i(u + c_i(t\Delta + o(\Delta)) - z) dF_i(z) \right. \\
    & \quad + e^{-\alpha_i t\Delta} e^{-\lambda_i t\Delta} \sum_{j=1}^n (\lambda_i p_{ij} \Delta + o(\Delta)) \varphi_i(u + c_i(t\Delta + o(\Delta))) \\
    & \quad + \left. e^{-\alpha_i t\Delta} e^{-\lambda_i t\Delta} o(\Delta) \right] = \sum_{t=0}^\infty [Q_t].
\end{split}
\end{equation}

    Fix a (sufficiently large) $T$ and represent $\varphi_i$ in the form $\varphi_i(u) = \left( \sum_{t=0}^M + \sum_{t=M+1}^{+\infty} \right) [Q_t]$, let $M \to \infty$, $\Delta \to 0$ such that $\Delta M = T$. Note that $\varphi_i(\cdot)$ are monotonically increasing; hence they are continuous almost everywhere [39, Theorem VIII.1.1, p. 192] and Riemann integrable [39, Theorem V.4.2, p. 129]. Therefore
    \begin{equation}\label{(7.4)}
        \sum_{t=0}^M [Q_t] \to \int_0^T e^{-t(\alpha_i + \lambda_i)} \left( \lambda_i \sum_{j=1}^m p_{ij} \varphi_j(u + c_i t) + \alpha_i \int_0^{u + c_i t} \varphi_i(u + c_i t - z) dF_i(z) \right) dt.
    \end{equation}

    Taking into account that $\varphi_i(\cdot) \le 1$, $F_i(\cdot) \le 1$, and $\sum_{j=1}^m p_{ij} \le 1$, for sufficiently small $\Delta$, such that $o(\Delta) \le (\alpha_i + \lambda_i) \Delta$, we have the estimates:
\begin{equation}\label{(7.5)}
\begin{split}
        \sum_{t=M+1}^{+\infty} [Q_t]
        & \le \sum_{t=M+1}^\infty \left[ e^{-\lambda_i t\Delta} e^{-\alpha_i t\Delta} (\alpha_i \Delta + o(\Delta)) + e^{-\alpha_i t\Delta} e^{-\lambda_i t\Delta} (\lambda_i \Delta + o(\Delta)) \right. \\
        & \quad + \left. e^{-\alpha_i t\Delta} e^{-\lambda_i t\Delta} o(\Delta) \right] \\
        & \le 5 \int_T^{+\infty} (\alpha_i + \lambda_i) e^{-t(\alpha_i + \lambda_i)} dt = 5 e^{-T(\alpha_i + \lambda_i)}.
\end{split}
\end{equation}

    Passing to the limit as $T \to +\infty$ in (\ref{(7.3)}) and taking into account (\ref{(7.4)}), (\ref{(7.5)}) we obtain:
    \begin{equation*}
        \varphi_i(u) = \int_0^\infty e^{-t(\alpha_i + \lambda_i)} \left( \lambda_i \sum_{j=1}^n p_{ij} \varphi_j(u + c_i t) + \alpha_i \int_0^{u + c_i t} \varphi_i(u + c_i t - z) dF_i(z) \right) dt.
    \end{equation*}
\end{proof}

From Theorem 7.1, the following statement follows.

%\begin{corollary}
{\bf Theorem 7.2.}
    The vector-function $\varphi(u) = \{\varphi_1(u), \varphi_2(u), \ldots, \varphi_m(u)\}$ satisfies the system of integro-differential equations:
    \begin{equation}\label{(7.6)}
        \frac{d\varphi_i}{du} = \frac{\alpha_i + \lambda_i}{c_i} \varphi_i(u) - \frac{\alpha_i}{c_i} \int_0^u \varphi_i(u - z) dF_i(z) - \frac{\lambda_i}{c_i} \sum_j p_{ij} \varphi_j(u), \quad i = 1, 2, \ldots, m.
    \end{equation}
%\end{corollary}

\begin{proof}
    Make the change of variables $y_i = u + c_i t$ in equation (\ref{(7.2)}):
    \begin{equation}\label{(7.7)}
        \varphi_i(u) = \frac{1}{c_i} e^{\frac{u(\alpha_i + \lambda_i)}{c_i}} \int_u^\infty e^{-\frac{y_i(\alpha_i + \lambda_i)}{c_i}} \left( \lambda_i \sum_{j=1}^n p_{ij} \varphi_j(y_i) + \alpha_i \int_0^{y_i} \varphi_i(y_i - z) dF_i(z) \right) dy_i.
    \end{equation}
    From this it is clear that $\varphi_i(u)$ are differentiable. Differentiating (\ref{(7.7)}), we obtain:
    \begin{align*}
        \frac{d\varphi_i(u)}{du}
        & = -\frac{1}{c_i} e^{\frac{u(\alpha_i + \lambda_i)}{c_i}} \left( \frac{1}{c_i} e^{-\frac{u(\alpha_i + \lambda_i)}{c_i}} \left( \lambda_i \sum_{j=1}^n p_{ij} \varphi_j(u) + \alpha_i \int_0^u \varphi_i(u - z) dF_i(z) \right) \right) \\
        & \quad + \frac{1}{c_i} \frac{(\alpha_i + \lambda_i)}{c_i} e^{\frac{u(\alpha_i + \lambda_i)}{c_i}} \int_u^\infty e^{-\frac{y_i(\alpha_i + \lambda_i)}{c_i}} \left( \lambda_i \sum_{j=1}^n p_{ij} \varphi_j(y_i) + \alpha_i \int_0^{y_i} \varphi_i(y_i - z) dF_i(z) \right) dy_i \\
        & = -\frac{1}{c_i} \left( \lambda_i \sum_{j=1}^n p_{ij} \varphi_j(u) + \alpha_i \int_0^u \varphi_i(u - z) dF_i(z) \right) + \frac{(\alpha_i + \lambda_i)}{c_i} \varphi_i(u).
    \end{align*}
\end{proof}

%\begin{corollary}
{\bf Theorem 7.3.}
    The functions $\{\varphi_i(u), i = 1, \ldots, m\}$ satisfy the system of integral equations
    \begin{equation}\label{(7.8)}
        \varphi_i(u) = \varphi_i(0) + \frac{\lambda_i}{c_i} \int_0^u \varphi_i(x) dx - \frac{\lambda_i}{c_i} \sum_{j=1}^m p_{ij} \int_0^u \varphi_j(x) dx + \frac{\alpha_i}{c_i} \int_0^u \varphi_i(u - z)(1 - F_i(z)) dz, 
    \end{equation}
 $i = 1, \ldots, m,$   where $\{\varphi_i(0)\}$ are (unknown) constants.
%\end{corollary}

In [91], equations (\ref{(7.6)}), (\ref{(7.8)}) were obtained under the assumption of differentiability of the functions $\varphi_i(u)$ and with $c_i = c > 0$.

Consider now a somewhat more general nonlinear model of a risk process in a Markovian environment than (\ref{(7.1)}).

We assume that in state $i$, the reserves (capital) of the insurance company in the absence of claims increase nonlinearly over time $t$ from the initial level $u$ to the value $U_i(u,t)$, $U_i(u,0) = u$. Similar situations are considered in the literature [65, 91], for example, when the premium income rate $c_i$ depends on the current reserve, $dU_i/dt = c_i(U_i)$, $U_i(0) = u$.

Note that insurance contracts may regulate the amounts $z_i(x)$ of insurance payments for the amount of loss $x$ in environment state $i$, for example, only a certain proportion of the possible loss may be insured, or there may be limits on the amounts of insurance payments. Let the random variable $z$ of claims in any insurance event have distribution function $F_i(z)$ and finite mean $\mu_i$, depending on the environment state $i$.

Then, similarly to Theorem 7.1, we can obtain the following system of integral equations, generalizing (\ref{(7.2)}).

%\begin{theorem}
{\bf Theorem 7.4.}
    The vector-function $\varphi(u) = \{\varphi_1(u), \varphi_2(u), \ldots, \varphi_m(u)\}$ of probabilities of non-ruin $\{\varphi_i(u)\}$ of the process from the initial position $(u, i)$ satisfies the system of integral equations ($i = 1, \ldots, m$):
    \begin{align}
        \varphi_i(u) & = \int_0^\infty e^{-t(\alpha_i + \lambda_i)} \left\{ \lambda_i \sum_{j=1}^m p_{ij} \varphi_j(U_i(u,t)) \right. \nonumber \\
        & \quad + \left. \alpha_i \int_0^{U_i(u,t)} \varphi_i(U_i(u,t) - z) dF_i(z) \right\} dt, \quad i = 1, \ldots, m.\label{(7.9)}
    \end{align}
%\end{theorem}

The systems of equations (\ref{(7.2)}) - (\ref{(7.9)}) are linear. We are interested in the solution satisfying the natural boundary conditions
\begin{equation}
    \lim_{u \to +\infty} \varphi_i(u) = 1, \quad i = 1, \ldots, m,\label{(7.10)}
\end{equation}
meaning that with an infinitely large initial capital, the process does not go bankrupt. Such a solution may not exist, for example, if the process goes bankrupt from any initial state, as is the case when the average claims per unit time are greater than the corresponding premiums. Then the solution is trivial, $\{\varphi_i(u) \equiv 0, i = 1, \ldots, m\}$.

\section{Solution of the system by the Laplace transform method}

The solution of the system of integral equations (\ref{(7.2)}) or its consequences (\ref{(7.6)}), (\ref{(7.8)}) with boundary conditions (\ref{(7.10)}) is a nontrivial problem [18, 65, 88]. In [40], however, it is shown that the system (\ref{(7.8)}), (\ref{(7.10)}) in some cases can be solved using the standard Laplace transform method. Applying the Laplace transform to equation (\ref{(7.8)}) and using
\begin{equation*}
    \int_0^{+\infty} e^{-\nu t} \int_0^t \varphi_i(u) du dt = \frac{\bar{\varphi}_i(\nu)}{\nu}, \quad
    \int_0^{+\infty} e^{-\nu t} \int_0^t \varphi_i(t - z)(1 - F_i(z)) dz dt = \bar{\varphi}_i(\nu) \frac{1 - \bar{F}_i(\nu)}{\nu},
\end{equation*}
where $\bar{\varphi}_i(\nu) = \int_0^{+\infty} e^{-\nu t} \varphi_i(t) dt$, $\bar{F}_i(\nu) = \int_0^{+\infty} e^{-\nu t} dF_i(t)$, we have the system of linear equations
\begin{equation}\label{(7.11)}
    \bar{\varphi}_i(\nu) = \frac{\varphi_i(0)}{\nu} + \frac{\lambda_i}{c_i} \frac{\bar{\varphi}_i(\nu)}{\nu} - \frac{\lambda_i}{c_i} \sum_{j=1}^m p_{ij} \frac{\bar{\varphi}_j(\nu)}{\nu} + \frac{\alpha_i}{c_i} \bar{\varphi}_i(\nu) \frac{1 - \bar{F}_i(\nu)}{\nu}, \quad i = 1, \ldots, m.
\end{equation}

Represent the system of linear equations (\ref{(7.11)}) in matrix form
\begin{equation}\label{(7.12)}
    \left( I - \frac{1}{\nu} L + \frac{1}{\nu} LP - \frac{1}{\nu} HD(\nu) \right) \bar{\varphi}(\nu) = \frac{\varphi(0)}{\nu}.
\end{equation}
Here $\bar{\varphi}(\nu) = (\bar{\varphi}_1(\nu), \ldots, \bar{\varphi}_m(\nu))'$; $I$ is the identity diagonal matrix, $L$ is the diagonal matrix with elements $\lambda_i/c_i$; $P = \{p_{ij}, p_{ii} = 0\}$; $D(\nu)$ is the diagonal matrix with elements $(1 - \bar{F}_i(\nu))$; $H$ is the diagonal matrix with elements $\alpha_i/c_i$; $\varphi(0) = \{\varphi_1(0), \ldots, \varphi_m(0)\}'$ is the column vector of unknown constants. Consequently,
\begin{equation*}
    \bar{\varphi}(\nu) = (\nu I - L(I - P) - HD(\nu))^{-1} \varphi(0).
\end{equation*}
Having found the inverse Laplace transform $\Phi(u)$ of the matrix $(\nu I - L(I - P) - HD(\nu))^{-1}$, from this we can compute the function $\varphi(u) = \Phi(u) \varphi(0)$ up to unknown constants $\varphi(0) = (\varphi_1(0), \ldots, \varphi_m(0))'$, which, in turn, can be found from the boundary conditions (\ref{(7.10)}): $\lim_{u \to +\infty} \Phi(u) \varphi(0) = \lim_{u \to +\infty} \varphi(u) = 1$, or verify that suitable constants do not exist and the system (\ref{(7.8)}), (\ref{(7.10)}) has no solution.

If the components $\bar{F}_i(\nu)$ of the matrix $D(\nu)$ are fractional-rational functions of $\nu$ (as in the example), then the components of the matrix $(\nu I - L(I - P) - HD(\nu))^{-1}$ are fractional-rational functions of $\nu$ and the inverse Laplace transform $\Phi(u)$ certainly exists.

\section{Solution of the system by the method of successive approximations}

In [41], the probability of non-ruin $\varphi(u)$ of the classical risk process was found from the integral equation (\ref{eqn:1.2}) by the Picard method of successive approximations:
\begin{equation*}
    \varphi^{k+1}(u) = 1 - \frac{\alpha \mu}{c} + \frac{\alpha}{c} \int_0^u \varphi^k(u - z)[1 - F(z)] dz, \quad k = 0, 1, \ldots,
\end{equation*}
where $\varphi^0(u)$ is some initial function, $1 - \alpha \mu / c \le \varphi^0(u) \le 1$. In this section, this iterative approach is applied to solve the system of integral equations (\ref{(7.9)}) with boundary conditions (\ref{(7.10)}).

Consider the following method of successive approximations for finding the functions $\{\varphi_1(u), \ldots, \varphi_m(u)\}$:
\begin{align}
    \varphi_i^{k+1}(u) & = \int_0^\infty e^{-t(\alpha_i + \lambda_i)} \left\{ \lambda_i \sum_{j=1}^m p_{ij} \varphi_j^k(U_i(u,t)) \right. \nonumber \\
    & \quad + \left. \alpha_i \int_0^{U_i(u,t)} \varphi_i^k(U_i(u,t) - z) dF_i(z) \right\} dt, \quad i = 1, \ldots, m,\label{(7.13)}
\end{align}
where $\{\varphi_i^0(u)\}$ are some initial functions such that $0 \le \varphi_i^0(u) \le 1$, $\varphi_i^0(+\infty) = 1$.

%\begin{assumption}
{\bf Assumtion 7.1.}
    Let in each state $i$ of the environment
    \begin{enumerate}
        \item[(a)] the reserve growth function $U_i(u,t)$ be monotonically increasing in $u$ and $t$, $\lim_{t \to +\infty} U_i(u,t) = +\infty$, and $U_i(u,t) \ge u + c_i t$ for $u \ge u_*$ with some constant $c_i \ge 0$;
        \item[(b)] $\alpha_i \mu_i / c_i < 1$, and thus, there exists a constant common to all states (an analogue of the Lundberg constant) $L > 0$ such that
        \begin{equation*}
            (\alpha_i / c_i) \int_0^{+\infty} e^{Lz} [1 - F_i(z)] dz \le 1.
        \end{equation*}
    \end{enumerate}
%\end{assumption}

Denote by $\{\bar{\varphi}^k(u)\}$ the sequence of approximations (\ref{(7.13)}) starting from the initial approximation $\varphi_i^0(u) \equiv 1$, and by $\{\underline{\varphi}^k(u)\}$ the sequence of approximations (\ref{(7.13)}) starting from the initial approximation $\varphi_i^0(u) = \varphi_*(u) = \max\{0, 1 - e^{-L(u - u_*)}\}$, $i = 1, \ldots, m$.

%\begin{lemma}
{\bf Lemma 7.1.}
    (On the convergence of the sequence of approximations from above). Let Assumption 7.1(a) be satisfied. Then in each environment state:
    \begin{enumerate}
        \item[(i)] the function $\bar{\varphi}_i^k(u)$ is monotonically non-decreasing in $u \ge 0$ and $0 \le \bar{\varphi}_i^k(u) \le 1$, $k \ge 0$;
        \item[(ii)] for each $u$, the sequence of values $\{\bar{\varphi}_i^k(u) \ge 0, k = 0, 1, \ldots\}$ decreases monotonically with increasing $k$;
        \item[(iii)] the limit function $\bar{\varphi}_i(u) = \lim_{k \to \infty} \bar{\varphi}_i^k(u)$ is the $i$-th component of the solution of the system of equations (\ref{(7.9)}), $\bar{\varphi}_i(u)$ is monotone (non-decreasing) in $u$ and $0 \le \bar{\varphi}_i(u) \le 1$;
        \item[(iv)] under the additional Assumption 7.1(b), for all $k$, $\bar{\varphi}_i^k(u) \ge \varphi_*(u) = \max\{0, 1 - e^{-L(u - u_*)}\}$, and thus $1 \ge \bar{\varphi}_i(u) \ge \varphi_*(u) = \max\{0, 1 - e^{-L(u - u_*)}\}$.
    \end{enumerate}
%\end{lemma}

\begin{proof}
    (i) The first assertion is proved by induction. For $k = 0$, it is obviously true. If $0 \le \bar{\varphi}_i^k(u) \le 1$, $i = 1, \ldots, m$, then
    \begin{align*}
        0 \le \bar{\varphi}_i^{k+1}(u) & \le \int_0^\infty e^{-t(\alpha_i + \lambda_i)} \left\{ \lambda_i \sum_{j=1}^m p_{ij} + \alpha_i \int_0^{U_i(u,t)} dF_i(z) \right\} dt \\
        & \le \int_0^\infty e^{-t(\alpha_i + \lambda_i)} \left\{ \lambda_i + \alpha_i \int_0^\infty dF_i(z) \right\} dt = 1.
    \end{align*}
    If all $\bar{\varphi}_j^k(u), j = 1, \ldots, m$, are monotone in $u$, then by the monotonicity of $U_i(u,t)$, the function $\bar{\varphi}_i^{k+1}(u)$ is also monotone in $u$.

    (ii) The second assertion follows by induction from the relations $\bar{\varphi}_i^1(u) \le 1 = \bar{\varphi}_i^0(u)$, $i = 1, \ldots, m$, and
    \begin{align*}
        \bar{\varphi}_i^{k+1}(u) - \bar{\varphi}_i^k(u)
        & = \int_0^\infty e^{-t(\alpha_i + \lambda_i)} \left\{ \lambda_i \sum_{j=1}^m p_{ij} (\bar{\varphi}_j^k(U_i(u,t)) - \bar{\varphi}_j^{k-1}(U_i(u,t))) \right. \\
        & \quad + \left. \alpha_i \int_0^{U_i(u,t)} (\bar{\varphi}_i^k(U_i(u,t) - z) - \bar{\varphi}_i^{k-1}(U_i(u,t) - z)) dF_i(z) \right\} dt. %\quad i = 1, \ldots, m.
    \end{align*}

    (iii) For any $u$, the sequence $\{\bar{\varphi}_i^k(u)\}$ decreases monotonically and is bounded below; therefore, there exists a limit $\bar{\varphi}_i(u) = \lim_{k \to \infty} \bar{\varphi}_i^k(u)$. The limit function $\bar{\varphi}_i(u)$, like all $\bar{\varphi}_i^k(u)$, is monotone in $u$ and satisfies $0 \le \bar{\varphi}_i(u) \le 1$. Passing to the limit in (\ref{(7.13)}) as $k \to \infty$, by the Lebesgue dominated convergence theorem, we obtain that the limit functions $\bar{\varphi}_i(u), i = 1, \ldots, m$, satisfy the system of equations (\ref{(7.9)}).

    (iv) Obviously, $\bar{\varphi}_i^0(u) \equiv 1 \ge \max\{0, 1 - e^{-L(u - u_*)}\}$. Let us show that if $\bar{\varphi}_i^k(u) \ge \max\{0, 1 - e^{-L(u - u_*)}\}$, $i = 1, \ldots, m$, then $\bar{\varphi}_i^{k+1}(u) \ge \max\{0, 1 - e^{-L(u - u_*)}\}$. Using the monotonicity of all $\bar{\varphi}_i^k$ and the fact that $\bar{\varphi}_j^k(U_i(u,t)) \ge \max\{0, 1 - e^{-L(U_i(u,t) - u_*)}\}$, $\bar{\varphi}_i^k(U_i(u,t) - z) \ge \max\{0, 1 - e^{-L(U_i(u,t) - u_* - z)}\}$, we get
    \begin{align*}
        \bar{\varphi}_i^{k+1}(u)
        & \ge \int_0^\infty e^{-t(\alpha_i + \lambda_i)} \left\{ \lambda_i \sum_{j=1}^m p_{ij} \max\{0, 1 - e^{-L(U_i(u,t) - u_*)}\} \right. \\
        & \quad + \left. \alpha_i \int_0^{U_i(u,t)} \max\{0, 1 - e^{-L(U_i(u,t) - u_* - z)}\} dF_i(z) \right\} dt, \quad i = 1, \ldots, m.
    \end{align*}
    Let $u \ge u_*$. Integrating the inner integral by parts and using $F_i(0) = 0$, we get
    \begin{align*}
        & \int_0^{U_i(u,t)} \max\{0, 1 - e^{-L(U_i(u,t) - u_* - z)}\} dF_i(z) \\
        & = \int_0^{\max\{0, U_i(u,t) - u_*\}} (1 - e^{-L(U_i(u,t) - u_* - z)}) dF_i(z) \\
        & = e^{-L(U_i(u,t) - u_*)} L \int_0^{U_i(u,t) - u_*} e^{Lz} F_i(z) dz.
    \end{align*}
    Thus,
    \begin{align*}
        \bar{\varphi}_i^{k+1}(u)
        & \ge \int_0^\infty e^{-t(\alpha_i + \lambda_i)} \left\{ \lambda_i \sum_{j=1}^m p_{ij} (1 - e^{-L(U_i(u,t) - u_*)}) \right. \\
        & \quad + \left. \alpha_i e^{-L(U_i(u,t) - u_*)} L \int_0^{U_i(u,t) - u_*} e^{Lz} F_i(z) dz \right\} dt, \quad i = 1, \ldots, m.
    \end{align*}
    Transform
    \begin{align*}
        L \int_0^{U_i(u,t)} e^{Lz} F_i(z) dz
        & = L \int_0^{U_i(u,t) - u_*} e^{Lz} (1 - (1 - F_i(z))) dz \\
        & = e^{L(U_i(u,t) - u_*)} - 1 - L \int_0^{U_i(u,t) - u_*} e^{Lz} (1 - F_i(z)) dz \\
        & \ge e^{L(U_i(u,t) - u_*)} - 1 - L \int_0^{+\infty} e^{Lz} (1 - F_i(z)) dz.
    \end{align*}
    Then
    \begin{align*}
        \bar{\varphi}_i^{k+1}(u)
        & \ge \int_0^\infty e^{-t(\alpha_i + \lambda_i)} \left\{ \lambda_i \sum_{j=1}^m p_{ij} (1 - e^{-L(U_i(u,t) - u_*)}) \right. \\
        & \quad + \left. \alpha_i e^{-L(U_i(u,t) - u_*)} \left( e^{L(U_i(u,t) - u_*)} - 1 - L \int_0^{+\infty} e^{Lz} (1 - F_i(z)) dz \right) \right\} dt \\
        & = \lambda_i \int_0^\infty e^{-t(\alpha_i + \lambda_i)} dt - \lambda_i \int_0^\infty e^{-t(\alpha_i + \lambda_i)} \sum_{j=1}^m p_{ij} e^{-L(U_i(u,t) - u_*)} dt \\
        & \quad + \alpha_i \int_0^\infty e^{-t(\alpha_i + \lambda_i)} dt \\
        & \quad - \alpha_i \left( 1 + L \int_0^{+\infty} e^{Lz} (1 - F_i(z)) dz \right) \cdot \int_0^\infty e^{-t(\alpha_i + \lambda_i) - L(U_i(u,t) - u_*)} dt.
    \end{align*}
    Taking into account that $\int_0^{+\infty} e^{Lz} (1 - F_i(z)) dz \le c_i / \alpha_i$ and $U_i(u,t) \ge u + c_i t$, we can write
    \begin{align*}
        \bar{\varphi}_i^{k+1}(u)
        & \ge 1 - \lambda_i e^{-L(u - u_*)} \sum_{j=1}^m p_{ij} \int_0^\infty e^{-t(\alpha_i + \lambda_i + L c_i)} dt \\
        & \quad - e^{-L(u - u_*)} \alpha_i (1 + L c_i / \alpha_i) \int_0^\infty e^{-t(\alpha_i + \lambda_i + L c_i)} dt \\
        & = 1 - e^{-L(u - u_*)} \left( \frac{\lambda_i}{(\alpha_i + \lambda_i + L c_i)} \sum_{j=1}^m p_{ij} + \frac{(\alpha_i + L c_i)}{(\alpha_i + \lambda_i + L c_i)} \right) = 1 - e^{-L(u - u_*)}.
    \end{align*}
    The lemma is proved.
\end{proof}

%\begin{lemma}
{\bf Lemma 7.2.}
    (On the convergence of the sequence of approximations from below). Let Assumptions 7.1(a), (b) be satisfied. Then in each environment state:
    \begin{enumerate}
        \item[(i)] each function $\underline{\varphi}_i^k(u)$ is monotonically non-decreasing (increasing) in $u \ge 0$ and $0 \le \underline{\varphi}_i^k(u) \le 1$, $k \ge 0$;
        \item[(ii)] the sequence of values $\{\underline{\varphi}_i^k(u) \ge 0, k = 0, 1, \ldots\}$ increases monotonically with increasing $k$, and thus for all $k$, $\underline{\varphi}_i^k(u) \ge \varphi_*(u) = \max\{0, 1 - e^{-L(u - u_*)}\}$;
        \item[(iii)] the limit function $\underline{\varphi}_i(u) = \lim_{k \to \infty} \underline{\varphi}_i^k(u)$ is the $i$-th component of the solution of the system of equations (\ref{(7.9)}), $\underline{\varphi}_i(u)$ is monotone (non-decreasing) in $u$ and $1 \ge \underline{\varphi}_i(u) \ge \varphi_*(u) = \max\{0, 1 - e^{-L(u - u_*)}\}$.
    \end{enumerate}
%\end{lemma}

\begin{proof}
    The proof of (i) is the same as in Lemma 7.1. From the proof of (iv) of Lemma 7.1, it follows that $\underline{\varphi}_i^1(u) \ge \max\{0, 1 - e^{-L(u - u_*)}\} = \underline{\varphi}_i^0(u)$. Hence, similarly to the proof of (ii) of Lemma 7.1, it follows that $\underline{\varphi}_i^{k+1}(u) \ge \underline{\varphi}_i^k(u) \ge \underline{\varphi}_i^0(u)$ for any $k \ge 0$ and $u \ge 0$. The proof of (iii) is analogous to the proof of assertion (iv) of Lemma 7.1.
\end{proof}

%\begin{lemma}
{\bf Corollary 7.1.}
    If the system of equations (\ref{(7.9)}) has a unique monotone solution $\{\varphi_1(u), \ldots, \varphi_m(u)\}$, then under Assumption 7.1(a), $\bar{\varphi}_i(u) = \varphi_i(u)$, and under Assumptions 7.1(a), (b), $\bar{\varphi}_i(u) = \underline{\varphi}_i(u) = \varphi_i(u)$, and for the obtained solution $\{\varphi_i(u)\}$, the boundary condition $\varphi_i(+\infty) = 1$, $i = 1, \ldots, m$, is satisfied.
%\end{lemma}

%\begin{theorem}
{\bf Theorem 7.5.}
    (On the convergence of the method of successive approximations). Assume that Assumptions 7.1(a), (b) are satisfied. Then
    \begin{itemize}
        \item[1)] system (\ref{(7.9)}), (\ref{(7.10)}) has a solution with monotone components $\varphi_i(u)$ such that $1 \ge \varphi_i(u) \ge \varphi_*(u) = \max\{0, 1 - e^{-L(u - u_*)}\}$;
        \item[2)] such a solution is unique;
        \item[3)] the sequence of approximations (\ref{(7.13)}) starting with $\{\varphi_i^0(u) \equiv 1\}$ converges monotonically from above, and starting with $\{\varphi_i^0(u) = \varphi_*(u) = \max\{0, 1 - e^{-L(u - u_*)}\}\}$, converges monotonically from below to the solution $\{\varphi_i(u)\}$ for each $i$;
        \item[4)] for any initial approximation $\{\varphi_i^0(u)\}$ such that $1 \ge \varphi_i^0(u) \ge \varphi_*(u) = \max\{0, 1 - e^{-L(u - u_*)}\}$, the corresponding sequence of approximations (\ref{(7.13)}) converges pointwise to the solution $\{\varphi_i(u)\}$.
    \end{itemize}
%\end{theorem}

\begin{proof}
    1) From Lemmas 7.1 and 7.2, it follows that the limit functions $\{\bar{\varphi}_1(u), \ldots, \bar{\varphi}_m(u)\}$ and $\{\underline{\varphi}_1(u), \ldots, \underline{\varphi}_m(u)\}$ are solutions of the system of equations (\ref{(7.9)}) and satisfy the conditions $1 \ge \bar{\varphi}_i(u) \ge \max\{0, 1 - e^{-L(u - u_*)}\}$, $1 \ge \underline{\varphi}_i(u) \ge \max\{0, 1 - e^{-L(u - u_*)}\}$. Thus, the first assertion on the existence of a solution of the system (\ref{(7.9)}), (\ref{(7.10)}) is proved.

    Let us show that $\bar{\varphi}_i(u) \equiv \underline{\varphi}_i(u)$, $i = 1, \ldots, m$. Since $\bar{\varphi}_i^0(u) \equiv 1 \ge \max\{0, 1 - e^{-L(u - u_*)}\} = \underline{\varphi}_i^0(u)$, similarly to the proof of Lemma 7.1(ii), it follows that $\bar{\varphi}_i^k(u) \ge \underline{\varphi}_i^k(u)$ for all $u \ge 0$, $k \ge 0$, and thus $\bar{\varphi}_i(u) \ge \underline{\varphi}_i(u)$ for all $u \ge 0$. By the linearity of the system (\ref{(7.9)}), the difference $\bar{\varphi}_i(u) - \underline{\varphi}_i(u)$ also satisfies the system of equations:
 
 \begin{equation}\label{(7.14)}
\begin{split}
        \bar{\varphi}_i(u) - \underline{\varphi}_i(u)
        & = \int_0^\infty e^{-t(\alpha_i + \lambda_i)} \left\{ \lambda_i \sum_{j=1}^m p_{ij} (\bar{\varphi}_j(U_i(u,t)) - \underline{\varphi}_j(U_i(u,t))) \right. \\
        & \quad + \left. \alpha_i \int_0^{U_i(u,t)} (\bar{\varphi}_i(U_i(u,t) - z) - \underline{\varphi}_i(U_i(u,t) - z)) dF_i(z) \right\} dt, \\
				\quad i = 1, \ldots, m.
  \end{split}
 \end{equation}

    Let us show that $\bar{\varphi}_i(u) - \underline{\varphi}_i(u) \equiv 0$, $i = 1, \ldots, m$. Suppose to the contrary that $\sup_{u \ge 0} [\bar{\varphi}_l(u) - \underline{\varphi}_l(u)] = \varepsilon > 0$ for some $l$ and $\sup_{u \ge 0} [\bar{\varphi}_j(u) - \underline{\varphi}_j(u)] \le \varepsilon$ for all $j$. Since $\bar{\varphi}_j(u) - \underline{\varphi}_j(u) \le \min\{1, e^{-L(u - u_*)}\}$ for all $j$, it follows that there exists a bounded sequence $\{u^s \le v, s = 0, 1, \ldots\}$ such that $\lim_{s \to \infty} u^s = u \le v$, $\lim_{s \to \infty} [\bar{\varphi}_l(u^s) - \underline{\varphi}_l(u^s)] = \varepsilon > 0$ and $\bar{\varphi}_l(u^s) - \underline{\varphi}_l(u^s) \le \varepsilon$. Moreover, $\bar{\varphi}_l(u^s) - \underline{\varphi}_l(u^s) \le \min\{1, \varepsilon, e^{-L(u^s - u_*)}\}$. From relation (\ref{(7.14)}) for $i = l$, $u = u^s$, it follows that
  \begin{equation}\label{(7.15)}
\begin{split}
        \bar{\varphi}_l(u^s) - \underline{\varphi}_l(u^s)
        & = \int_0^\infty e^{-t(\alpha_l + \lambda_l)} \left\{ \lambda_l \sum_{j=1}^m p_{lj} (\bar{\varphi}_j(U_l(u^s,t)) - \underline{\varphi}_j(U_l(u^s,t))) \right. \\
        & \quad + \left. \alpha_l \int_0^{U_l(u^s,t)} (\bar{\varphi}_l(U_l(u^s,t) - z) - \underline{\varphi}_l(U_l(u^s,t) - z)) dF_l(z) \right\} dt \\
        & \le \int_0^{+\infty} e^{-t(\alpha_l + \lambda_l)} \left\{ \lambda_l \varepsilon + \alpha_l \int_0^{U_l(u^s,t)} \min\{1, \varepsilon, e^{-L(U_l(u^s,t) - u_* - z)}\} dF_l(z) \right\} dt \\
        & \le \varepsilon (\alpha_l + \lambda_l) \int_0^\infty e^{-t(\alpha_l + \lambda_l)} F_l(U_l(u^s,t)) dt.
  \end{split}
 \end{equation}
    Two cases are possible: $F_l(z) < 1$ for all $z \ge 0$ and $F_l(z) = 1$ for all $z \ge \bar{z}_l$. Consider the first case. Passing to the limit in (\ref{(7.15)}) as $s \to \infty$ and using the monotonicity of the distribution function $F_l(z)$, we get
    \begin{equation*}
        1 \le \int_0^{+\infty} (\alpha_l + \lambda_l) e^{-t(\alpha_l + \lambda_l)} F_l(U_l(v,t)) dt.
    \end{equation*}
    But by the continuity of $U_i(u,\cdot)$, we have $F_l(U_l(v,t)) < 1$ for all $t \ge 0$. Hence, we obtain a contradiction that $\int_0^{+\infty} (\alpha_l + \lambda_l) e^{-t(\alpha_l + \lambda_l)} F_l(U_l(v,t)) dt < 1$. Thus, in this case, $\bar{\varphi}_j(u) \equiv \underline{\varphi}_j(u)$, $j = 1, \ldots, m$.

    Consider the second case, when $F_l(z) = 1$ for all $z \ge \bar{z}_l$. From inequality (\ref{(7.15)}), it follows that
 \begin{equation*}
\begin{split}
        \bar{\varphi}_l(u^s) - \underline{\varphi}_l(u^s)
        & \le \int_0^{+\infty} e^{-t(\alpha_l + \lambda_l)} \left\{ \lambda_l \varepsilon + \alpha_l \int_0^{U_l(u^s,t)} \min\{1, \varepsilon, e^{-L(U_l(u^s,t) - u_* - z)}\} dF_l(z) \right\} dt \\
        & \le \int_0^{+\infty} e^{-t(\alpha_l + \lambda_l)} \left\{ \lambda_l \varepsilon + \alpha_l \min\{1, \varepsilon, e^{-L(U_l(u^s,t) - u_* - \bar{z}_l)}\} \int_0^{\min\{\bar{z}_l, U_l(u^s,t)\}} dF_l(z) \right\} dt \\
        & \le \frac{\varepsilon \lambda_l}{(\alpha_l + \lambda_l)} + \alpha_l \int_0^{+\infty} e^{-t(\alpha_l + \lambda_l)} \min\{1, \varepsilon, e^{-L(U_l(u^s,t) - u_* - \bar{z}_l)}\} dt.
        \end{split}
 \end{equation*}

    Passing to the limit here as $s \to \infty$, we get
    \begin{equation*}
        \varepsilon \le \frac{\varepsilon \lambda_l}{(\alpha_l + \lambda_l)} + \alpha_l \int_0^{+\infty} e^{-t(\alpha_l + \lambda_l)} \min\{1, \varepsilon, e^{-L(U_l(u,t) - u_* - \bar{z}_l)}\} dt
    \end{equation*}
    or
    \begin{align*}
     &   \varepsilon \le (\alpha_l + \lambda_l) \int_0^\infty e^{-t(\alpha_l + \lambda_l)} \min\{1, \varepsilon, e^{-L(U_l(u,t) - u_* - \bar{z}_l)}\} dt \\
		&\le (\alpha_l + \lambda_l) \int_0^\infty e^{-t(\alpha_l + \lambda_l)} \min\{1, \varepsilon, e^{-L(U_l(0,t) - u_* - \bar{z}_l)}\} dt.
    \end{align*}
    But since $U_l(0,t) \to +\infty$, we have $e^{-L(U_l(0,t) - u_* - \bar{z}_l)} \to 0$ as $t \to +\infty$ and $\min\{\varepsilon, e^{-L(U_l(0,t) - u_* - \bar{z}_l)}\} < \varepsilon$ for all sufficiently large $t$. Hence, it follows that
    \begin{equation*}
        (\alpha_l + \lambda_l) \int_0^{+\infty} e^{-t(\alpha_l + \lambda_l)} \min\{1, \varepsilon, e^{-L(U_l(0,t) - u_* - \bar{z}_l)}\} dt < \varepsilon.
    \end{equation*}
    The obtained contradiction proves that in the case under consideration, we also have $\bar{\varphi}_j(u) \equiv \underline{\varphi}_j(u)$, $j = 1, \ldots, m$.

    Let us now prove assertion 2) of the theorem on the uniqueness of the solution $\{1 \ge \varphi_i(u) \ge \max\{0, 1 - e^{-L(u - u_*)}\}\}$. Since $\bar{\varphi}_j^0(u) \equiv 1 \ge \varphi_j(u) \ge \max\{0, 1 - e^{-L(u - u_*)}\} = \underline{\varphi}_j^0(u)$ for all $j$, by the linearity of equation (7.13), $\bar{\varphi}_i^1(u) \ge \varphi_i(u) \ge \underline{\varphi}_i^1(u)$ and for any $k$, $\bar{\varphi}_i^k(u) \ge \varphi_i(u) \ge \underline{\varphi}_i^k(u)$ for all $i$. Hence, passing to the limit as $k \to \infty$, taking into account the proved identity $\bar{\varphi}_i(u) \equiv \underline{\varphi}_i(u)$, we get $\bar{\varphi}_i(u) \equiv \varphi_i(u) \equiv \underline{\varphi}_i(u)$, $i = 1, \ldots, m$.

    Assertion 3) follows from Lemmas 7.1 and 7.2.

    Assertion 4) is proved similarly to assertion 2). Let $1 \ge \varphi_i^0(u) \ge \max\{0, 1 - e^{-L(u - u_*)}\}$, i.e., $\bar{\varphi}_i^0(u) \ge \varphi_i^0(u) \ge \underline{\varphi}_i^0(u)$, $i = 1, \ldots, m$. Then for the sequence of approximations $\{\varphi_i^k(u), i = 1, \ldots, m\}$ generated by method (\ref{(7.13)}) from the initial approximation $\{\varphi_i^0(u), i = 1, \ldots, m\}$, for any $k$, we have $\bar{\varphi}_i^k(u) \ge \varphi_i^k(u) \ge \underline{\varphi}_i^k(u)$, $i = 1, \ldots, m$. Hence, it follows that for any $i$, there exists a limit $\varphi_i(u) = \lim_{k \to \infty} \varphi_i^k(u) \equiv \bar{\varphi}_i(u) \equiv \underline{\varphi}_i(u)$. The theorem is proved.
\end{proof}

\section{Example}

Consider an example of a risk process on a Markov chain with two states $i = 1, 2$ and transition matrix
\begin{equation*}
    P = \begin{pmatrix}
        0 & 1 \\
        1 & 0
    \end{pmatrix}.
\end{equation*}
It is assumed that only the premium income rates depend on the state, $c_1 = 2$, $c_2 = 22/21$. All state change intensities $\lambda_i = 1$, claim arrival intensities $\alpha_i = 1$, and claims are exponentially distributed with $F_i(x) = 1 - e^{-x}$, $x \ge 0$, $i = 1, 2$.

The Laplace transform of the function $F_i(x) = 1 - e^{-x}$ is $\bar{F}_i(\nu) = 1/(1 + \nu)$. Consequently, the matrices that appear in equation (\ref{(7.12)}) have the form
\begin{equation*}
    L = H = \begin{pmatrix}
        1/2 & 0 \\
        0 & 21/22
    \end{pmatrix}, \quad
    D(\nu) = \begin{pmatrix}
        \nu/(\nu + 1) & 0 \\
        0 & \nu/(\nu + 1)
    \end{pmatrix}.
\end{equation*}
Using formula (\ref{(7.11)}), we write:
\begin{align*}
    \bar{\varphi}(\nu)
    & = \begin{pmatrix}
        \frac{-1 + 2\nu^2}{2(1 + \nu)} & \frac{1}{2} \\
        \frac{21}{22} & -\frac{21}{22} + \nu - \frac{21\nu}{22(1 + \nu)}
    \end{pmatrix}^{-1} \cdot \begin{pmatrix}
        \varphi_1(0) \\
        \varphi_2(0)
    \end{pmatrix} \nonumber \\
    & = \begin{pmatrix}
        \frac{42 + 82\nu + 4\nu^2 - 44\nu^3}{22\nu + 85\nu^2 + 40\nu^3 - 44\nu^4} & \frac{22(1 + \nu)^2}{22\nu + 85\nu^2 + 40\nu^3 - 44\nu^4} \\
        \frac{42(1 + \nu)^2}{22\nu + 85\nu^2 + 40\nu^3 - 44\nu^4} & \frac{22(1 + \nu)(2\nu^2 - 1)}{22\nu + 85\nu^2 + 40\nu^3 - 44\nu^4}
    \end{pmatrix} \cdot \begin{pmatrix}
        \varphi_1(0) \\
        \varphi_2(0)
    \end{pmatrix}.
\end{align*}
After inverse Laplace transform, we get
\begin{equation*}
    \varphi(u) \approx \begin{pmatrix}
        a_{11}(u) & a_{12}(u) \\
        a_{21}(u) & a_{22}(u)
    \end{pmatrix} \begin{pmatrix}
        \varphi_1(0) \\
        \varphi_2(0)
    \end{pmatrix},
\end{equation*}
\begin{align*}
    a_{11}(u) & \approx 1.90909 - 0.8290 e^{-0.7634u} - 1.1124 e^{-0.3275u} + 0.2862 e^{2u}, \\
    a_{12}(u) & \approx 1 + 0.0304 e^{-0.7634u} - 0.6806 e^{-0.3275u} + 0.2862 e^{2u}, \\
    a_{21}(u) & \approx 1.90909 - 0.0581 e^{-0.7634u} - 1.2993 e^{-0.3275u} + 0.6678 e^{2u}, \\
    a_{22}(u) & \approx 1 + 0.0213 e^{-0.7634u} - 0.7949 e^{-0.3275u} + 0.8163 e^{2u}.
\end{align*}
Determine the unknown constants $\varphi_1(0)$ and $\varphi_2(0)$ from the conditions (\ref{(7.10)}) that $\varphi_i(+\infty) = 1$, $i = 1, 2$. We get $\varphi_1(0) = 0.3667$, $\varphi_2(0) = 0.3$. Thus,
\begin{equation}\label{(7.16)}
\begin{split}
    \varphi_1(u) & \approx 1 - 0.0213 e^{-0.7634u} - 0.6121 e^{-0.3275u}, \\
    \varphi_2(u) & \approx 1 + 0.0149 e^{-0.7634u} - 0.7149 e^{-0.3275u}.
\end{split}
\end{equation}

Now we find the probabilities of non-ruin of the risk process from Example 1 by the method of successive approximations (\ref{(7.13)}).

The iterative process (\ref{(7.13)}) for the matrix from Example 1 can be rewritten in the form:
\begin{equation}\label{(7.17)}
\begin{aligned}
    \varphi_1^{k+1}(u) & = \frac{1}{c_1} \int_u^{+\infty} e^{-\frac{\alpha_1 + \lambda_1}{c_1}(y - u)} \left\{ \lambda_1 \varphi_2^k(y) + \int_0^y \varphi_1^k(y - z) dF_1(z) \right\} dy, \\
    \varphi_2^{k+1}(u) & = \frac{1}{c_2} \int_u^{+\infty} e^{-\frac{\alpha_2 + \lambda_2}{c_2}(y - u)} \left\{ \lambda_2 \varphi_1^k(y) + \int_0^y \varphi_2^k(y - z) dF_2(z) \right\} dy, \\
		 k = 0, 1, \ldots
\end{aligned}
\end{equation}

The functional form of the iterations, the presence of double integrals and infinite limits of integration, as well as the possible large number of iterations, require an economical numerical implementation of the process (\ref{(7.17)}) and consideration of rounding errors.

If $\varphi_i^0(u) = \bar{\varphi}_i^0(u) \equiv 1, i = 1, 2$, then, according to Lemma 7.1, the sequence of approximations (\ref{(7.17)}) decreases monotonically and converges from above to the solution of the corresponding system of equations (\ref{(7.9)}). The course of iterations for $\varphi_1^k(u) = \bar{\varphi}_1^k(u)$ is shown in Fig. 7.2.

\begin{figure}[h]
    \centering
    \includegraphics[width=0.8\textwidth]{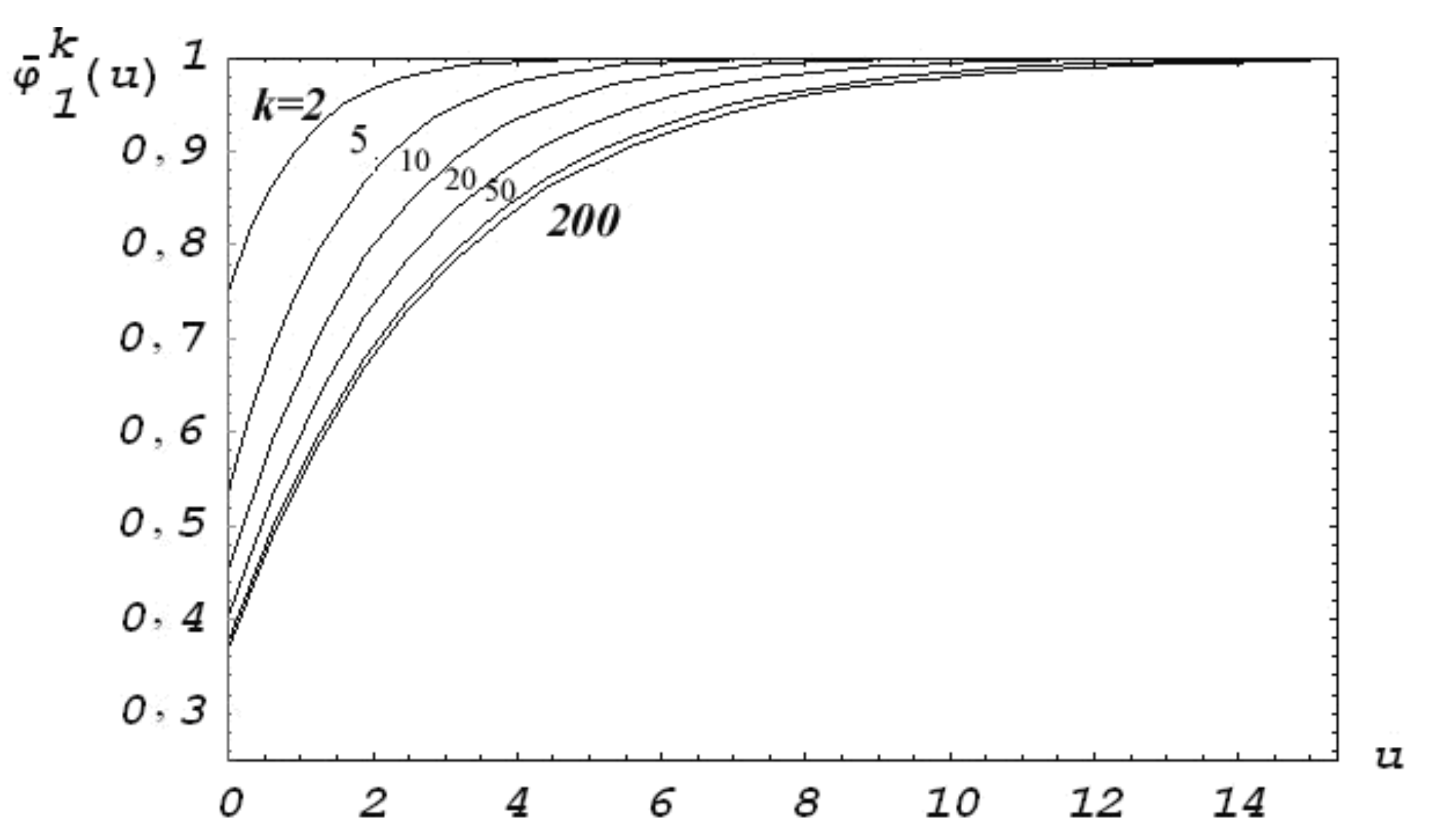}
    \caption{Iterations from above for $\varphi_1^k(u)$}
    \label{fig:7.2}
\end{figure}

If $\varphi_i^0(u) = \underline{\varphi}_i^0(u) = 1 - e^{-Lu}, i = 1, 2$, where $L = 1/22$ is the common Lundberg constant, then according to Lemma 7.2, the sequence of approximations (\ref{(7.17)}) increases monotonically and converges from below to the solution of the system (\ref{(7.9)}), (\ref{(7.10)}); the course of the corresponding iterations for $\varphi_1^k(u) = \underline{\varphi}_1^k(u)$ is shown in Fig. 7.3.

\begin{figure}[h]
    \centering
    \includegraphics[width=0.8\textwidth]{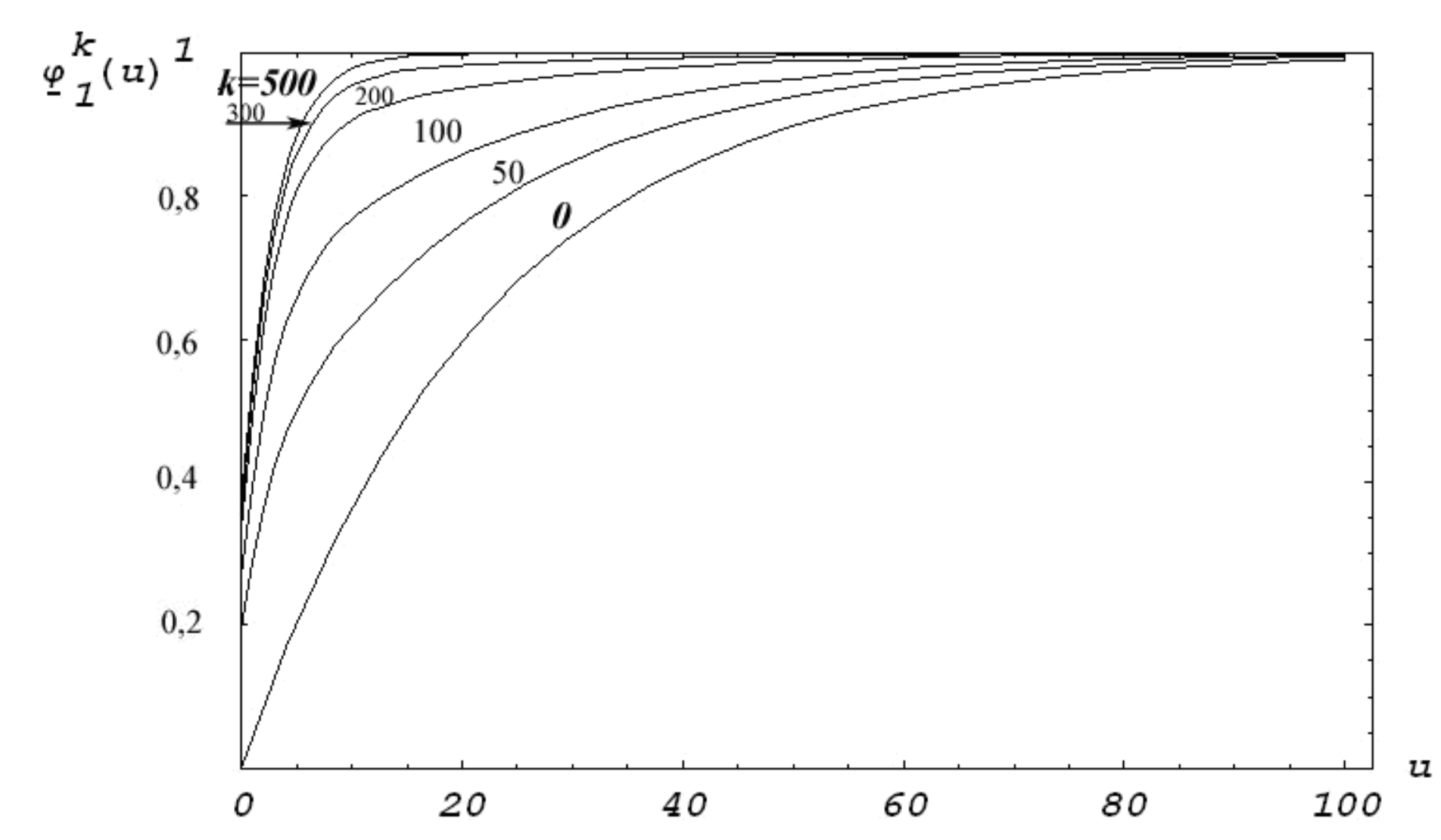}
    \caption{Iterations from below for $\varphi_1^k(u)$}
    \label{fig:7.3}
\end{figure}

Fig. 7.4 graphically shows the error of the method $\Delta_1(u) = \bar{\varphi}_1^{200}(u) - \underline{\varphi}_1^{500}(u)$ as a function of $u$, and Fig. 7.5 shows the dependence of the error $\Delta_1^k = \sup_{0 \le u < \infty} (\bar{\varphi}_1^k(u) - \underline{\varphi}_1^k(u))$ on the number of iterations $k$ (for approximations from above – 1; from below – 2). For the approximations $\varphi_2^k(u)$, the picture of convergence to the solution is similar.

\begin{figure}[h]
    \centering
    \includegraphics[width=0.8\textwidth]{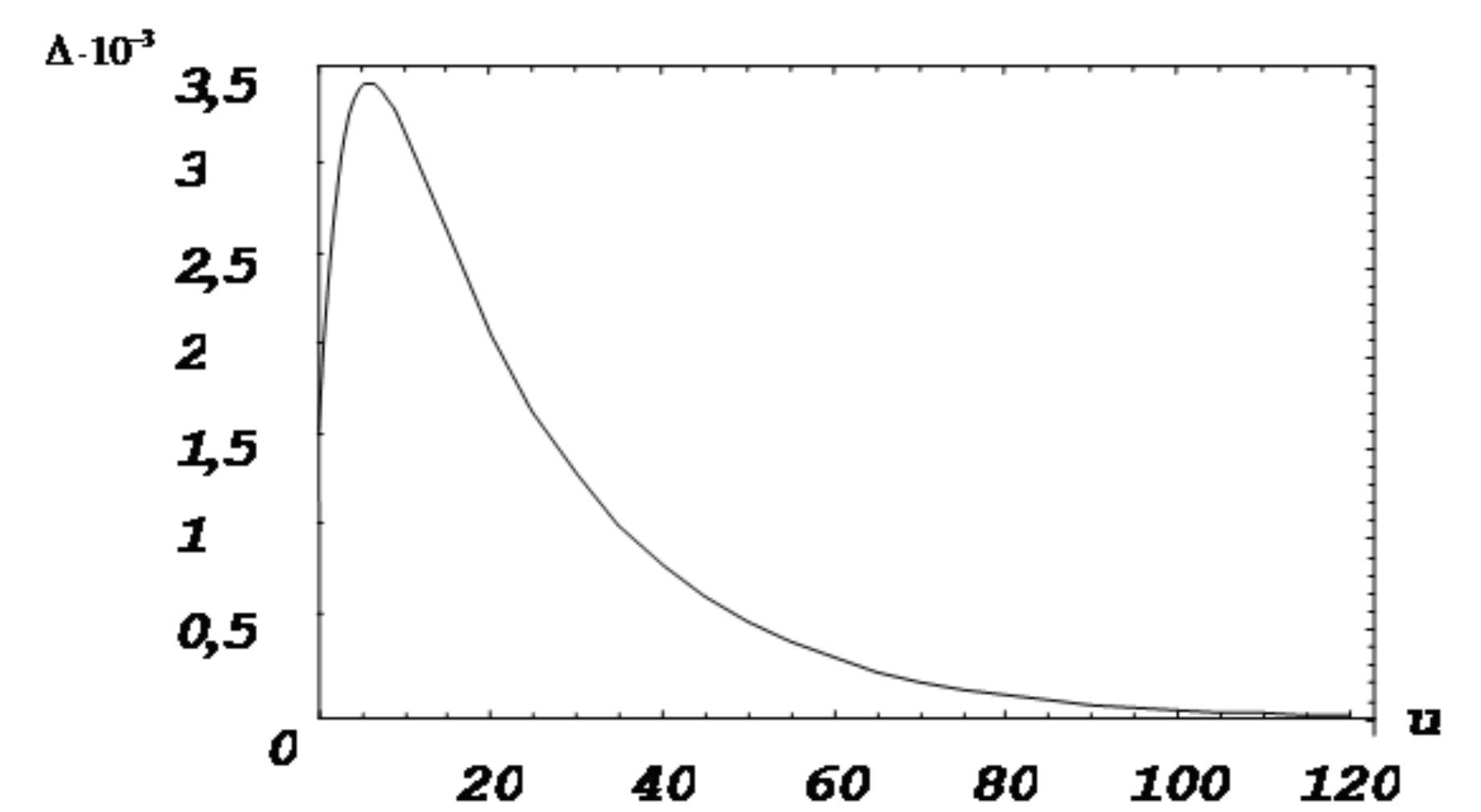}
    \caption{Error of the method $\Delta_1(u)$}
    \label{fig:7.4}
\end{figure}

\begin{figure}[h]
    \centering
    \includegraphics[width=0.8\textwidth]{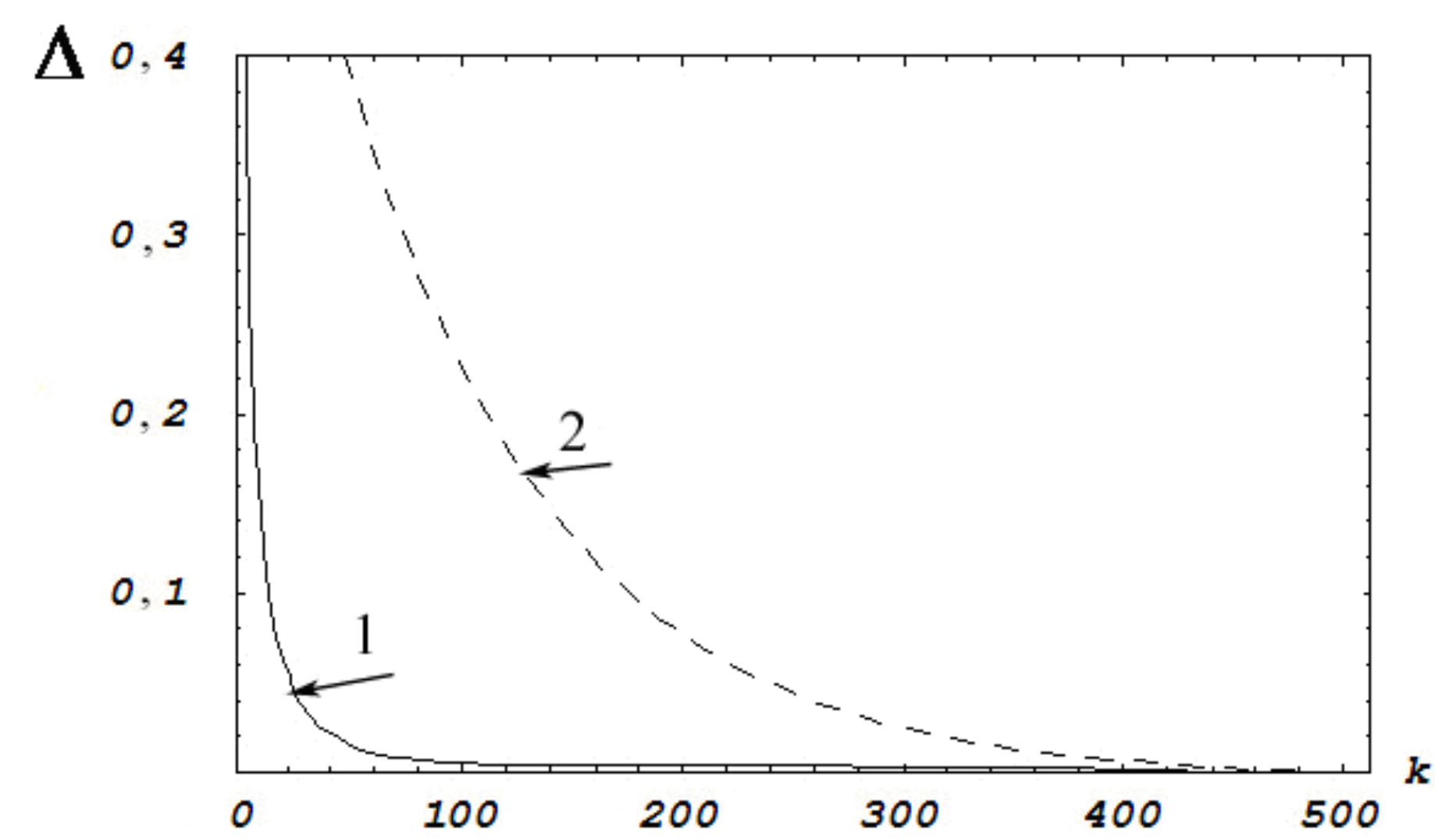}
    \caption{Error $\Delta_1^k$ as a function of the number of iterations}
    \label{fig:7.5}
\end{figure}

The approximate solutions $\{\bar{\varphi}_i^{200}(u), \underline{\varphi}_i^{500}(u), i = 1, 2\}$ coincide with the solution obtained in Example 1 with an accuracy of order 0.003, and the difference is caused by errors in the numerical implementation of the algorithm.

\section{Conclusions to the section}

In this section, a general risk process is considered that describes the stochastic evolution of the capital of an insurance company in a Markovian random environment with nonlinear premium income and general insurance contracts. It is shown that the probabilities of ruin of the process (company) as functions of the initial state generally satisfy a certain system of integral equations with boundary conditions at infinity. Sufficient conditions for the existence of a solution to this system are established and the method of successive approximations for calculating the ruin probabilities is justified. When starting from the unit initial approximation, the iterations converge to the solution from above, and when starting from the exponential Lundberg bound, from below, and thus, it is always possible to estimate the accuracy of the obtained approximate solution. The performance of the method is verified on a numerical example of a risk process on a Markov chain with two states.

\chapter*{CONCLUSIONS}
\addcontentsline{toc}{chapter}{CONCLUSIONS}

The dissertation obtained the following main results.
\begin{enumerate}
    \item Integral equations for the probability of (non-)ruin of various generalizations of the classical risk process are obtained, in particular, the process with deterministic premiums depending on current capital, with random premiums, with non-Poisson streams of premiums and claims, and for the risk process in a random Markovian environment.
    \item Necessary and sufficient conditions for existence and general sufficient conditions for existence and uniqueness of solutions of the considered integral equations and systems of equations of insurance mathematics are established.
    \item The method of successive approximations for the numerical or analytical solution of the considered integral equations of insurance mathematics is theoretically justified.
    \item A methodology for estimating the accuracy of the method of successive approximations for solving integral equations of insurance mathematics by constructing approximations to the solution from below and above is developed.
    \item The proposed method of successive approximations for solving integral equations of insurance mathematics has been tested on a number of numerical examples, and its comparison with the Monte Carlo method and with approximate estimates of solutions has been carried out.
\end{enumerate}

The method of successive approximations for solving integral equations of insurance mathematics developed in the dissertation allows:
\begin{itemize}
    \item increasing the accuracy of actuarial calculations;
    \item calculating (within the framework of the chosen mathematical model of the company) the probability of ruin of an insurance company with any predetermined accuracy;
    \item verifying the accuracy and applicability of various known approximate formulas for the ruin probability;
    \item improving the accuracy of empirical approximations iteratively if necessary;
    \item estimating the accuracy and adjusting the parameters of the Monte Carlo method when calculating the ruin probability by simulation modeling.
\end{itemize}

\bibliographystyle{plain}
%\bibliography{references}

\end{document}